\documentclass[a4paper]{amsart}
\usepackage{amsfonts}
\usepackage{amssymb}
\usepackage{euscript}
\usepackage{xspace}
\usepackage{enumerate}

\newcommand{\LPmn}{\PrMn{\LP}matricially normed\xspace}

\newcommand{\MAX}{\mathop{\mathrm{MAX}}}
\newcommand{\MIN}{\mathop{\mathrm{MIN}}}

\newcommand{\Ast}{${}^\ast$\xspace}
\newcommand{\bAst}{${}^{(\ast)}$\xspace}
\newcommand{\bast}{^{(\ast)}}

\newcommand{\SA}{{\mathrm{sa}}}         
\newcommand{\Cay}{{\EuScript{C}}}       
\newcommand{\Md}{{\mathfrak{M}}}        
\newcommand{\pM}{{\mathfrak{m}}}        
\newcommand{\Rd}{{\mathfrak{R}}}        
\newcommand{\Bd}{{\mathfrak{B}}}        
\newcommand{\NS}[1][V]{{\mathsf{#1}}} 
\newcommand{\ns}[1][v]{{\mathsf{#1}}} 
\newcommand{\Ball}{{\mathop\mathrm{Ball}}}
\newcommand{\clBall}{\operatorname{\cl{Ball}}}
\newcommand{\Cone}{{\mathop\mathrm{Cone}}}

\newcommand{\OA}[1][A]{{\mathsf{#1}}} 
\newcommand{\UA}{{\mathsf{U}}}        
\newcommand{\FA}{{\mathsf{F}}}        
\newcommand{\QA}{{\mathsf{Q}}}        
\newcommand{\Shift}{{\mathsf{Shift}}} 
\newcommand{\Toep}{{\mathrm{Toeplitz}}} 
\newcommand{\mO}[1]{{\mathbf{#1}}}      
\newcommand{\Fant}{{\EuScript{F}}}      
\newcommand{\U}{\mathrm{U}}             
\newcommand{\PSU}{\mathrm{PSU}}         
\newcommand{\dist}{\mathrm{dist}}       

\newcommand{\blank}{{\llcorner\!\lrcorner}}
\newcommand{\PrM}[2][60]{${#2}$\hbox{-}\penalty#1\hspace*{0pt}}
\newcommand{\PrMn}[2][10000]{${#2}$\hbox{-}\penalty#1\hspace*{0pt}}

\newcommand{\diag}{\mathop{\mathrm{diag}}}

\newcommand{\Aut}{\mathop\mathrm{Aut}}

\newcommand{\Bound}{\mathbb{B}}

\newcommand{\Caratheodory}{Ca\-ra\-th\'e\-o\-do\-ry\xspace}
\newcommand{\Cech}{\v{C}ech\xspace}
\newcommand{\Fantappie}{Fan\-tap\-pi\`e\xspace}
\newcommand{\Hausdorff}{Haus\-dorff\xspace}
\newcommand{\Hoermander}{H\"or\-man\-der\xspace}

\newcommand{\Lempert}{Lem\-pert\xspace}

\newcommand{\Moebius}{M\"o\-bi\-us\xspace}
\newcommand{\Nevanlinna}{Ne\-van\-lin\-na\xspace}
\newcommand{\Poincare}{Poin\-ca\-r\'e\xspace}
\newcommand{\Reiffen}{Reiffen\xspace}

\newcommand{\SzNagy}{Sz\H{o}\-ke\-fal\-vi-Nagy\xspace}
\newcommand{\Vigue}{Vi\-gu\'e\xspace}

\newcommand{\Artanh}{\mathop\mathrm{arc\,tanh}}

\newcommand{\CstarAlgebra}{C\PrMn{^\ast}al\-ge\-bra\xspace}

\newcommand{\5}[2]{\langle{{#1},{#2}}\rangle}

\newcommand{\cl}[1]{{\overline{#1}}}
\newcommand{\bd}{\partial}

\newcommand{\C}{{\mathbb{C}}}
\newcommand{\R}{{\mathbb{R}}}
\newcommand{\D}{{\mathbb{D}}}
\newcommand{\bdD}{{\bd\D}}
\newcommand{\clD}{{\cl{\mathbb{D}}}}

\newcommand{\N}{{\mathbb{N}}}

\newcommand{\conj}[1]{\overline{#1}}

\newcommand{\HO}{\mathcal{O}}
\newcommand{\clHO}{\cl{\mathcal{O}}}
\newcommand{\BH}{\mathrm{H}^\infty}

\newcommand{\NBC}{\mathrm{C}}
\newcommand{\CINF}{\mathrm{C}^\infty}
\newcommand{\LP}[1][\infty]{\mathrm{L}^{#1}}

\newcommand{\Pol}{\mathcal{P}}
\newcommand{\clPol}{\cl{\mathcal{P}}}
\newcommand{\Rat}{\mathcal{R}}
\newcommand{\clRat}{\cl{\mathcal{R}}}

\newcommand{\Shf}[1]{\mathbf{U}_{#1}}

\newcommand{\Tg}{\mathbf{T}}
\newcommand{\CTg}{\mathbf{T}^\ast}

\newcommand{\co}{\mathop\mathrm{co}}

\newcommand{\lin}{\mathop\mathrm{lin}}

\newcommand{\Ran}{\mathop\mathrm{Ran}}
\newcommand{\Ker}{\mathop\mathrm{Ker}}
\newcommand{\Spec}{\mathop\mathrm{Spec}\nolimits}

\newcommand{\RE}{\mathop\mathrm{Re}}

\newcommand{\cbn}[1]{\|{#1}\|_{\mathrm{cb}}}

\newcommand{\ide}{\mathsf{I}}

\newcommand{\ID}[1][{}]{\mathrm{id}_{{#1}}}

\newcommand{\Mat}[1][n]{\mathbb{M}_{#1}}
\newcommand{\HilS}[1][H]{\EuScript{#1}}

\newcommand{\Unze}[1]{#1^{+}}
\newcommand{\const}{\mathrm{const}}

\newtheorem{theorem}{Theorem}[section]
\newtheorem{proposition}[theorem]{Proposition}
\newtheorem{corollary}[theorem]{Corollary}
\newtheorem{lemma}[theorem]{Lemma}
\newtheorem{claim}[theorem]{Claim}
\theoremstyle{definition}
\newtheorem{definition}{Definition}[section]
\theoremstyle{remark}
\newtheorem{example}{Example}[section]
\newtheorem{remark}[theorem]{Remark}
\newtheorem{problem}{Problem}[section]
\newtheorem{conjecture}[theorem]{Conjecture}

\begin{document}

\title[Quotients of operator algebras]{Finite dimensional quotients of
 commutative operator algebras}

\author{Ralf Meyer}

\address{ Mathematisches Institut\\
	  Westf. Wilhelms-Universit\"at M\"unster\\
	  Einsteinstr.~62\\
	  48149 M\"unster\\
	  Germany }

\email{rameyer@math.uni-muenster.de}

\begin{abstract}
  It is shown that the matrix normed structure of a non-unital operator algebra
  determines that of its unitization.  This makes the study of certain unital
  operator algebras much easier and provides several interesting
  counterexamples.

  Every two-dimensional, unital operator algebra is completely isometrically
  isomorphic to an algebra of $2\times2$-matrices, and every contractive
  homomorphism between two such algebras is completely contractive.  This is
  used to define analogues for commutative, unital operator algebras of the
  Carath\'eodory distance and the Carath\'eodory-Reiffen metric on complex
  manifolds.

  There exists an isometric, completely contractive map between
  three-di\-mensional Q-algebras that is not completely isometric.
  Moreover, for every strongly pseudoconvex domain~$\mathfrak{M}$, the
  algebra $\mathrm{H}^\infty(\mathfrak{M})$ has a contractive
  representation by $3\times3$-matrices that is not completely
  contractive.

  Recently, Arveson introduced the $d$-shift as a model for $d$-contractions.
  Completely isometric representations of quotients of the operator algebra
  generated by the $d$-shift are computed explicitly.  For $d=1$, this gives a
  version of Nevanlinna-Pick theory.  It happens that every quotient of finite
  dimension~$r$ has a completely isometric representation by $r\times
  r$-matrices.  Finally, the class of operator algebras with this property is
  investigated.
\end{abstract}


\maketitle

\section{Introduction}
\label{sec:Intro}

An operator algebra~$\OA$ is just a subalgebra of $\Bound(\HilS)$, the bounded
operators on a Hilbert space~$\HilS$.  The operator norm on~$\Bound(\HilS)$
gives rise to a norm on~$\OA$.  Moreover, $\OA\otimes\Mat\subset
\Bound(\HilS\otimes\C^n)$ in a natural way, where~$\Mat$ denotes the algebra of
\PrM{n\times n}matrices with the usual C\Ast-norm.  Thus every operator algebra
comes with natural norms on all tensor products $\OA\otimes\Mat$.  The main
interest of this article lies on this additional structure.

It is the framework for the \emph{model theory} of (commuting)
operators on Hilbert space.  The starting point of model theory was
the \SzNagy dilation theorem \cite{Riesz-Nagy:90},
\cite{Nagy-Foias:70}, which asserts that for any contraction
$T\in\Bound(\HilS)$, there is an essentially unique unitary
operator~$U$ on a Hilbert space~$\HilS[D]$ containing~$\HilS$ such
that $T^n=P_{\HilS} U^nP_{\HilS}$ for all $n\in\N$.  Here~$P_{\HilS}$
denotes the projection onto the subspace~$\HilS$.  The unitary~$U$ is
called a \emph{(power) dilation} of~$T$.  This allows to apply the
rich theory of unitary operators to the study of contractions.

Until recently, attempts at a generalization of this result to 
\emph{(multi)operators}, i.e.\ \PrMn{d}tuples of commuting operators
$\mO{T}=(T_1,\ldots,T_d)$ on a common Hilbert space, have been rather
unsuccessful.  It was soon discovered by And\^o~\cite{Ando:63} that
two commuting contractions still have a unitary dilation, which,
however, is no longer unique.  But Parrot~\cite{Parrott:70} gave a
counterexample of three commuting contractions that do not have a
unitary dilation.  At the same time, Arveson \cite{Arveson:69},
\cite{Arveson:72} found necessary and sufficient conditions for the
existence of dilations in terms of the matrix normed structure
described above.  An accessible account of these classical results is
Paulsen's monograph~\cite{Paulsen:86}.

In~\cite{Arveson:97}, finally an interesting model theory for multi\PrM{}operators
is developed.  A \emph{\PrMn{d}contraction} is a multi\PrM{}operator
$\mO{T}=(T_1,\ldots,T_d)$ such that
\begin{equation}  \label{equ:dcontraction}
\|T_1\xi_1+\dots+T_d\xi_d\|^2 \le \|\xi_1\|^2 + \dots + \|\xi_d\|^2
\end{equation}
for all $\xi_1,\ldots,\xi_n\in\HilS$.  An equivalent condition is that
the \PrM{d\times d}matrix
\begin{equation}  \label{equ:dcCond}
\begin{pmatrix}
        T_1 & T_2 & \dots & T_d \\
        0   & 0   & \dots & 0 \\
	0   & 0   & \dots & 0 \\
        \hdotsfor{4}
\end{pmatrix}
\end{equation}
is a contraction.  Hence the innocent-looking
condition~\eqref{equ:dcontraction} already involves the matrix normed structure
on $\Bound(\HilS)\otimes\Mat[d]$.  Indeed, the norm on the operator algebra
generated by~$\mO{T}$, in general, does not contain enough information to
determine whether~$\mO{T}$ is a \PrMn{d}contraction.

A particular \PrMn{d}contraction is the \PrMn{d}shift
$\mO{S}=(S_1,\ldots,S_d)$, acting on the Hilbert space~$H^2_d$, which will be
described in greater detail below.  The main result of~\cite{Arveson:97} is
that every \PrMn{d}contraction has an essentially unique dilation with
additional properties to a multi\PrM{}operator of the form $n\cdot \mO{S}\oplus
\mO{Z}$, where $n\cdot \mO{S}$ stands for the direct sum of~$n$ copies
of~$\mO{S}$ acting on~$(H^2_d)^n$ and~$\mO{Z}$ is a normal multi\PrM{}operator
with spectrum contained in the boundary of the standard Euclidean unit ball
$\D_d\subset\C^d$.  Moreover, every operator that has such a dilation is a
\PrMn{d}contraction.  The normal part~$\mO{Z}$ is often missing, e.g.\
if the matrix in~\eqref{equ:dcCond} has norm strictly less than~$1$.

For $d=1$, the \PrMn{1}shift is just the usual unilateral shift, and the above
dilation is the von Neumann-Wold decomposition of an isometry.  This is almost
as good as a unitary dilation, and actually what is needed in several
applications of dilation theory, e.g.~\cite{Agler:90}.  For $d>1$, however, the
\PrMn{d}shift is no longer subnormal\footnote{A \emph{subnormal
(multi)operator} is the restriction of a normal (multi)operator to an invariant
subspace.}.  This is the reason why the model theory for \PrMn{d}contractions
was discovered so late.

Let~$\Omega$ be a compact space.  A \emph{uniform algebra} on~$\Omega$
is a closed unital subalgebra of $\NBC(\Omega)$ that separates the
points of~$\Omega$.  More generally, a \emph{function algebra}
on~$\Omega$ is a subalgebra of $\NBC(\Omega)$.  A function
algebra~$\FA$ comes with a natural matrix normed structure, viewing
elements of $\FA\otimes\Mat$ as functions~$f$ from~$\Omega$ to~$\Mat$
with norm $\|f\|_\infty=\sup_{\omega\in\Omega} \|f(\omega)\|$.  It
follows easily from spectral theory that every function algebra~$\FA$
on~$\Omega$ is also a function algebra on its spectrum $\Spec(\FA)$.
Thus the space~$\Omega$ is not very important.  The operator algebra
generated by a subnormal operator is always (completely isometric to)
a function algebra.

Function algebras are interesting in their own right because they
arise in complex analysis.  A typical example of a uniform algebra is
the algebra $\BH(\Md)$ of bounded holomorphic functions on a complex
manifold~$\Md$.  However, its matrix normed structure has not been of
great use in complex analysis so far.  Some results in complex
analysis, e.g.\ \Lempert's theorem, can be proved quite naturally
using dilation theory (see~\cite{Agler:90}), but there are also more
elementary proofs \cite{Meyer:96}, \cite{Meyer:97} of those results.

If~$\OA$ is a unital operator algebra and $\ide\subset\OA$ is an ideal, there
is a natural matrix normed structure on the quotient algebra $\QA=\OA/\ide$.
Somewhat surprisingly, the resulting object is again an \emph{abstract operator
algebra}, i.e.\ it can be represented \emph{completely isometrically} on a
Hilbert space~$\HilS$ \cite{Blecher-Ruan-Sinclair:90}.  A representation
$\rho\colon \QA\to\Bound(\HilS)$ is called completely isometric if all the maps
$$
        \rho_{(n)}
  =     \rho\otimes\ID[\Mat]\colon
         \QA\otimes\Mat[n]\to \Bound(\HilS)\otimes\Mat[n]
        \cong\Bound(\HilS\otimes\C^n)
$$
are isometric.  The proof due to Blecher, Ruan, and Sinclair uses an axiomatic
characterization of unital operator algebras and is not constructive.  Indeed,
completely isometric representations of quotients are often quite hard to find.

A commutative operator algebra has a rich ideal structure and a lot of finite
dimensional quotients.  It might be expected that these quotients are simpler
than the original object.  This is true in the sense that all information
contained in the quotient is already contained in the original object.  But in
fact, taking quotients instead often brings the hidden complexity of an
operator algebra to the surface.

Function algebras appear to be rather simple objects, and this is certainly
true, say, from the point of view of spectral theory.  But \emph{Q\PrM{}algebras},
i.e.\ quotients of function algebras, are among the most complicated operator
algebras.  Let~$\FA$ be a function algebra on~$\Omega$ and let~$\ide$ be the
ideal
\begin{equation}  \label{equ:ideal}
        \ide
=       \ide(\omega_1,\dots,\omega_n)
=       \{ f\in\FA\mid f(\omega_1)=\dots=f(\omega_n)=0 \}
\end{equation}
with distinct points $\omega_1,\dots,\omega_d\in\Omega$.  Then an element
of $\QA=\FA/\ide$ is determined by its function values at the points
$\omega_1,\dots,\omega_d$.  Write~$[f]$ for the projection of $f\in\FA$
in~$\QA$.  The norm of an element in $\FA\otimes\Mat$ can, in principle, be
computed from its range.  The norm on~$\QA$, however, depends on the existence
of a solution of an interpolation problem with prescribed range.

Solving interpolation problems with prescribed range is one of the most
difficult problems in complex analysis.  If there are just two points, i.e.\
$d=1$, this amounts to computing the \Caratheodory distance of $\omega_1$
and~$\omega_2$, but this is possible only in very few special cases.  For
$d>3$, quotients of the disk algebra $\Pol(\clD)$ can be computed explicitly.
It is also possible to compute the norm on quotients of $\Pol(\clD^2)$
\cite{Agler:??}, but there seems to be no theory for other examples.  Thus
quotients of function algebras usually cannot be computed.  Moreover, in those
cases where they can, they tend to have as few completely contractive
representations as possible (Theorem~\ref{the:ACQClassify} and
Theorem~\ref{the:approxCTg}).

The most basic requirement on a model theory is a simple criterion which
operators can be modeled.  Formally, this can be translated into a criterion
to decide whether a given representation of a certain operator algebra
(generated by the model multi\PrM{}operator) is completely contractive.  Thus
taking quotients is very relevant to model theory: The (completely)
contractive representations of the quotient $\OA/\ide$ are precisely the
(completely) contractive representations of~$\OA$ whose kernel
contains~$\ide$.  Thus a criterion to decide whether a representation of~$\OA$
is completely contractive automatically applies to representations
of~$\OA/\ide$.

In the opinion of the author, a good model theory should actually have the
stronger property that completely isometric representations of all quotients
can be computed explicitly.  This criterion is not met by function algebras
(with the exception of the disk algebra).  However, completely isometric
representations of the operator algebra $\Shift_d$ generated by the
\PrMn{d}shift~$\mO{S}$ can be computed explicitly.  The last part of this
article is concerned with the theory of the quotients of $\Shift_d$.

The completely isometric representations of quotients of $\Shift_d$ can be
easily written down explicitly, but the proof that they are indeed completely
isometric is formal and not constructive.  It is based on the fact that the
quotient is again an abstract operator algebra.

It turns out that every quotient of $\Shift_d$ of finite dimension~$r$ has a
completely isometric representation by \PrMn{r\times r}matrices.  A unital,
commutative operator algebra with this property is said to have \emph{minimal
quotient complexity}.  The author conjectures that this property already
essentially characterizes the quotients of $\Shift_d$.  More precisely, the
conjecture is that a finite dimensional, indecomposable operator algebra of
minimal quotient complexity is either a quotient of $\Shift_d$ of the transpose
of such a quotient.

\goodbreak

The structure of this article is as follows:

Section~\ref{sec:Notation} contains general terminology used in this article.

In Section~\ref{sec:Unitization}, the matrix normed structure on the
unitization of an operator algebra~$\OA$ is shown to be determined by the
matrix normed structure of~$\OA$ and not to depend on the choice of a
completely isometric representation.  Hence for many purposes a unital operator
algebra can be replaced by a \PrMn{1}codimensional ideal~$\ide$.

This is particularly useful if the multiplication on the \PrMn{1}codimensional
ideal~$\ide$ is the zero map.  This means that~$\ide$ essentially is just a
linear space of operators on a Hilbert space.  Such a space with its matrix
normed structure is called an \emph{operator (vector) space}, and every
operator space~$\NS$, endowed with the zero multiplication, occurs as a maximal
ideal of a unital abstract operator algebra~$\Unze{\NS}$.  This abstract
operator algebra is uniquely determined and called the \emph{trivial
unitization} of~$\NS$.  A linear map $\rho\colon \NS_1\to\NS_2$ can be extended
uniquely to a unital homomorphism $\Unze\rho\colon
\Unze{\NS_1}\to\Unze{\NS_2}$.  Then $\|\rho\|_{(n)}=\|\Unze\rho\|_{(n)}$, and
thus $\cbn{\Unze\rho}=\cbn{\rho}$.  Here
$$
	\|\rho\|_{(n)}=\|\rho_{(n)}\|
\qquad\text{and}\qquad
	\cbn{\rho}=\sup_{n\in\N} \|\rho_{(n)}\|.
$$
Consequently, the representation theory of a trivial unitization is precisely
as well\PrM{}behaved, or pathological, as the linear representation theory of
the underlying operator space.  This is the basis for several counterexamples.

Finally, the unitization technique yields a slight refinement of the theorem of
Smith that a \PrMn{d}contractive linear mapping into $\Mat[d]$ is automatically
completely contractive. If~$\OA$ is a unital, commutative operator algebra,
then a \PrMn{d-1}contractive, unital homomorphism $\OA\to\Mat[d]$ is completely
contractive.  This generalizes Agler's discovery that every contractive, unital
homomorphism $\OA\to\Mat[2]$ is completely contractive.

In Section~\ref{sec:twodimOpalg}, two\PrM{}dimensional, unital operator
algebras are studied.  The unitization technique reduces the classification of
these operator algebras to that of one\PrM{}dimensional operator algebras,
which is rather trivial.  Every two-dimensional, unital operator algebra has a
completely isometric representation by \PrM{2\times2}matrices.  The norms and
\emph{complete norms} of all algebraic isomorphisms between
two\PrM{}dimensional operator algebras are computed.  It turns out that, for
any such automorphism, $\|\rho\|=\cbn{\rho}$.  By taking quotients, this
immediately generalizes to unital homomorphisms between unital operator
algebras which have rank two as linear maps.  This generalizes results
previously known about representations of function algebras by
\PrM{2\times2}matrices (\cite{Agler:90}, \cite{Salinas:91}, \cite{Chu:92},
\cite{Fu-Russo:95}, \cite{Meyer:96}).

In Section~\ref{sec:CaraDist}, the classification of two\PrM{}dimensional
operator algebras is used to define analogues of the \Caratheodory
pseudodistance and the \Caratheodory-\Reiffen pseudometric
\cite{Jarnicki-Pflug:93} for (commutative) unital operator algebras~$\OA$.
Since \Caratheodory has not much to do with these objects, they are called
\emph{quotient distance} and \emph{quotient metric}.  The quotient distance is
a distance on the spectrum of~$\OA$ that is defined precisely as in the case of
complex manifolds.  Essentially, it describes the equivalence class of the
two\PrM{}dimensional, unital operator algebra $\OA/\ide(\omega_1,\omega_2)$.
The tangent space of~$\OA$ consists of pairs $(\omega;d)$, where
$\omega\in\Spec(\OA)$ and~$d$ is a derivation of~$\OA$ at~$\omega$, and the
quotient metric is just the norm of this linear functional.  These definitions
also make sense for noncommutative operator algebras.  But since elements of
the form $[f,g]$ are in the kernel of all characters and derivations at some
point of the spectrum, the spectrum and the tangent space ignore any
noncommutativity of~$\OA$.

For every model theory, there should be an explicit criterion which
\PrM{2\times2}matrices can be modeled.  Hence the quotient distance and metric
can be computed explicitly for the operator algebras involved.

An application of the quotient distance is a simple criterion when a finite
dimensional, commutative, unital operator algebra is \emph{decomposable}, i.e.\
when it decomposes into a non-trivial orthogonal sum of two ideals: This
happens if and only if some two\PrM{}dimensional quotient is isometric to
$\NBC(\{0,1\})$, i.e.\ iff some two\PrM{}dimensional quotient can be decomposed
orthogonally.

In Section~\ref{sec:Transposition}, the usual transposition operation on~$\Mat$
is defined for abstract operator algebras.  This operation is isometric, but
usually not completely isometric.  However, transposition is completely
isometric for Q\PrM{}algebras and thus serves as a simple criterion to show
that an operator algebra is not a Q\PrM{}algebra.  Moreover, transposition
provides the easiest examples of homomorphisms between two\PrMn{}dimensional,
non-unital operator algebras that are not completely isometric.  The transpose
of the operator algebra $\Shift_d$ is interesting for theoretical purposes
because it models the adjoints of \PrMn{d}contractions and has similar formal
properties as~$\Shift_d$.

If~$\OA$ is a unital operator algebra and $\omega\in\Spec(\OA)$, define
$\ide(\omega)$ as in~\eqref{equ:ideal} and $\ide(\omega)^2$ as the closure of
$\ide(\omega)\cdot \ide(\omega)$.  Then define $\OA(\omega)=\OA/\ide(\omega)^2$
and the \emph{cotangent space} $\CTg_\omega\OA=\ide(\omega)/\ide(\omega)^2$
of~$\OA$ at~$\omega$.  The tangent space with the quotient metric is its normed
dual.  Of course, $\OA(\omega)$ is the trivial unitization of $\CTg_\omega\OA$.
Section~\ref{sec:Cotangent} contains several counterexamples of badly behaved
cotangent spaces of function algebras.

For certain function algebras, the tangent and cotangent spaces were already
introduced by Paulsen in~\cite{Paulsen:92}.  He was interested in determining
when it happens that every contractive representation of $\Rat(\cl{\Md})$ is
completely contractive, where $\Rat(\cl{\Md})$ denotes the algebra of rational
functions without singularities in~$\cl{\Md}$, considered as a subalgebra of
$\NBC(\cl{\Md})$.  If~$\Md$ is a balanced domain, it is easy to show that the
cotangent space of $\Rat(\cl{\Md})$ at zero is completely equivalent to
$\MIN(\NS)$, where~$\NS$ is the normed space whose unit ball is the polar
of~$\Md$ and $\MIN(\NS)$ denotes the minimal \LPmn structure on~$\NS$.  Thus it
is very rare that every contractive linear representation
of~$\CTg_0\Rat(\cl{\Md})$ is completely contractive: This happens iff
$\MIN(\NS)=\MAX(\NS)$.  Paulsen shows in~\cite{Paulsen:92} that this cannot
hold for $\dim \NS\ge5$.  Thus $\Rat(\cl{\Md})$ has a contractive
representation that is not completely contractive whenever $\dim\Md\ge5$.

In Section~\ref{sec:Cotangent}, Paulsen's negative result is extended to
bounded, strongly pseudoconvex domains: If~$\Md$ is a bounded, strongly
pseudoconvex domain in $\C^d$, $d\neq1$, then $\BH(\Md)$ has a contractive
representation by \PrM{3\times3}matrices that is not \PrMn{2}contractive.  This
uses that the matrix normed structure on the cotangent space can be computed
approximately near the boundary and approaches $\MIN(\ell^2_d)$.  Moreover, the
result of \Lempert (\cite{Lempert:88}) that there exist domains not
biholomorphic to~$\D_2$ but with tangent space at some point isometric
to~$\ell^2_2$ provides an example of an isometric, completely contractive
homomorphism between three\PrM{}dimensional Q\PrM{}algebras that is not
\PrM{2}isometric.  Using quite different techniques, such a homomorphism has
recently been obtained by Paulsen in~\cite{Paulsen:97}.

This adds to the evidence that Q\PrM{}algebras are rather complicated operator
algebras.  In Section~\ref{sec:ComplexAnalysis}, two reasons will be given why
the operator algebraic viewpoint should not be expected to give deep results in
complex analysis.  First, the matrix normed structure on Q\PrM{}algebras does
not distinguish between certain nice and pathological objects.  Secondly, there
is more structure on a Q\PrM{}algebra than the matrix normed structure: A
Q\PrM{}algebra~$\QA$ comes with natural norms on $\QA\otimes\NS$ for all normed
spaces~$\NS$.  The relevance of this structure for interpolation theory is
discussed in~\cite{Meyer:97a}.  From the point of view of interpolation theory,
the restriction to the spaces~$\Mat$ is artificial.

Section~\ref{sec:modelD} starts with a brief account of Arveson's model theory
for \PrMn{d}contractions.  There are two other approaches to Arveson's Hilbert
space~$H^2_d$ with interesting consequences.  Both use the explicit form
$u_x(z)=(1-\5{z}{x})^{-1}$, $x,z\in\Ball(\ell^2_d)$, of the reproducing kernel
for the Hilbert space~$H^2_d$.  That this is a reproducing kernel means
essentially that $\5{f}{u_x}=f(x)$ for all $f\in H^2_d$, and this computation
is done by Arveson.  But he does not explore the connections this opens to
other areas of mathematics.

Up to a constant, $u_x(z)$ is the $1/(d+1)$st power of the Bergman kernel for
the ball, so that~$H^2_d$ is a ``twisted Bergman space''.  These spaces have
been studied by harmonic analysts because they carry a natural projective
representation of the automorphism group $\PSU(d,1)$ of the ball.  Especially,
there is a natural projective representation of $\PSU(d,1)$ on~$H^2_d$.
Consequently, the closed operator algebra~$\Shift_d$ generated by the
\PrMn{d}shift retains all symmetries of the ball, i.e.\ $f\mapsto f\circ a$ is
a completely isometric isomorphism of~$\Shift_d$ for all $a\in\Aut(\D_d)$.
This is the main consequence of the twisted Bergman space picture of~$H^2_d$.
Moreover, some less important results of Arveson are special cases of results
about general twisted Bergman spaces on symmetric domains.
In~\cite{Bagchi-Misra:96}, Bagchi and Misra compute the spectrum of $\Shift_d$
and show that the \PrMn{d}shift is not subnormal.  However, they fail to single
out~$H^2_d$ as a case of special interest because their main goal is to
determine when the analogue of the \PrMn{d}shift on a twisted Bergman space is
subnormal and when the closed unit ball is a (complete) spectral set for it.

The reproducing kernel of~$H^2_d$ is also related to the \emph{\Fantappie
transform} for \PrMn{\C}convex domains.  For the special case of~$\D_d$, the
\Fantappie transform is a linear bijection $\Fant\colon
\HO(\D_d)'\to\HO(\cl{\D_d})$ and thus induces a canonical bilinear form on
$\HO(\cl{\D_d})\subset\HO(\D_d)$.  In order to make this form sesquilinear,
some complex conjugate signs must be added.  The resulting conjugate \Fantappie
transform~$\conj{\Fant}$ is given by $\bigl(\conj{\Fant}(l)\bigr)(x)=l(u_x)$
for $l\in\HO(\D_d)'$ and is a bijection onto the space of coanalytic functions
on a neighborhood of~$\cl{\D_d}$.  Thus $\5{f}{g}= \bigl(\conj{\Fant}^{-1}
(\conj{g})\bigr)(f)$ defines a sesquilinear form on $\HO(\cl{\D_d})$.  This is
nothing but the inner product on~$H^2_d$.  Thus if $l\in\HO(\D_d)'$, then inner
products in~$H^2_d$ involving $\conj\Fant(l)$ can be computed easily by
$\5{f}{\conj\Fant(l)}= l(f)$.

In Section~\ref{sec:QuotientMultiplier}, completely isometric representations
of quotients of the algebra $\Shift_d$ generated by Arveson's \PrMn{d}shift are
computed explicitly.  It turns out that, up to a self-adjoint part coming from
the boundary~$\bd\D_d$, the representation of $\Shift_d/\ide$ on $H^2_d\ominus
\ide\cdot H^2_d$ is completely isometric.  The proof that this representation
is indeed completely isometric starts with any completely isometric
representation and shows that it is a quotient of several copies of
$H^2_d\ominus \ide H^2_d$.  Thus it makes essential use of the fact that
quotients of operator algebras are again abstract operator algebras.
Therefore, the proof is not as elementary as it may seem at first glance.
Another disadvantage of the proof is that it is not constructive: Given $F\in
(\Shift_d/\ide)_{(n)}$, it does not produce a representative
$\hat{F}\in(\Shift_d)_{(n)}$ with equal norm.

A special case of particular interest is if the ideal is of the form
$\ide(x_1,\dots,x_m)$ for distinct points $x_1,\dots,x_m\in\D_d$.  There exists
$F\in (\Shift_d)_{(n)}$ with prescribed values $F(x_j)=y_j$ and positive,
invertible real part if and only if the block matrix~$P(F)$ with entries
$$
P(F)_{ij}= \frac{y_i+y_j^\ast}{1-\5{x_i}{x_j}}\in\Mat[n]
$$
is positive definite and invertible.  A similar criterion allows to check
whether $\|F\|_{(n)}<1$.  These formulas are direct generalizations of the
existence criteria of \Nevanlinna-Pick theory, which is the special case $d=1$
of the above.  However, an important difference between $d=1$ and $d>1$ is that
transposition is no longer completely isometric for $d>1$.  If there exists
$F\in(\Shift_d)_{(n)}$ with $F(x_j)=y_j$ and $\|F\|_{(n)}<1$, there need not
exist $F\in(\Shift_d)_{(n)}$ with $F(x_j)=y_j^t$ and $\|F\|_{(n)}<1$.

The computation of two\PrM{}dimensional quotients of $\Shift_d$ yields the
quotient distance and metric for $\Shift_d$.  The result turns out to be the
usual \Caratheodory distance and metric for~$\D_d$.  This computation is
greatly facilitated by the invariance of $\Shift_d$ under automorphisms of the
ball.

An important property of quotients of $\Shift_d$ is that they have completely
isometric representations by \PrM{r\times r}matrices if they have
dimension~$r$.  Unital, commutative operator algebras with this property are
said to have \emph{minimal quotient complexity}.  This is a very rare property.
Indeed, the author conjectures that the only finite dimensional operator
algebras with this property are orthogonal direct sums of quotients of
$\Shift_d$ and $\Shift_d^t$.  As a first step towards proving this conjecture,
the case of trivial unitizations is studied.  We obtain that cotangent spaces
of closed operator algebras of minimal complexity are completely isometric to
$\Bound(\C,\HilS)$ or $\Bound(\HilS,\C)$ for some Hilbert space~$\HilS$ with
the obvious matrix normed structure.  Especially, the only possibilities for
finite dimensional cotangent spaces are $\CTg_0\Shift_d$ and $\CTg_0\Shift_d^t$
with $d\in\N$.  In particular, if a function algebra~$\FA$ has minimal
complexity, its cotangent spaces can have dimension at most~$1$.  This excludes
algebras like $\Pol(\Md)$, $\Rat(\Md)$, etc., if $\dim\Md\ge 2$.

\subsection{Notation}
\label{sec:Notation}

If $\HilS$ is a Hilbert space and $\HilS[K]\subset\HilS$, then
$\HilS\ominus\HilS[K]$ denotes the orthogonal complement of~$\HilS[K]$
in~$\HilS$.  If $L\colon \NS_1\to\NS_2$ is a linear map, $\Ker L$ and $\Ran L$
denote its kernel and range, respectively.  For $p\in[1,\infty]$, $n\in\N$,
let~$\ell^p_n$ be the \PrMn{n}dimensional \PrMn{\ell^p}space.  Let~$\Mat$ be
the \PrM{n\times n}matrices with the usual C\PrM{^\ast}norm and let~$\Mat[n,m]$
be the \PrM{n\times m}matrices, normed as operators from~$\ell^2_m$
to~$\ell^2_n$.

Let~$\D$ be the open unit disk in~$\C$ and~$\bdD$ its boundary, the unit
circle.  If~$\NS$ is a normed space, let $\Ball(\NS)$ be its open unit ball and
$\clBall(\NS)$ its closed unit ball.  Let $\C^\ast=\C\setminus\{0\}$.

If~$\Omega$ is a compact space, $\NBC(\Omega)$ stands for the \CstarAlgebra of
\PrMn{\C}valued continuous functions on~$\Omega$.  If~$\Md$ is a complex
manifold, $\BH(\Md)$ denotes the algebra of bounded holomorphic functions
on~$\Md$ and $\HO(\Md,\Rd)$ the algebra of holomorphic maps from~$\Md$
to~$\Rd$.  As usual, $\HO(\Md)=\HO(\Md,\C)$.  If $K\subset\C^n$ is a compact
set, let $\HO(K)$ be the algebra of functions holomorphic in a neighborhood
of~$K$ and $\clHO(K)$ its closure in $\NBC(K)$.  Similar conventions apply to
the algebras $\Pol(K)$ of polynomials and $\Rat(K)$ of rational functions
without singularities in~$K$.  View them as unital subalgebras of $\NBC(K)$ and
write $\clPol(K)$ and~$\clRat(K)$ for their closures.

The \emph{\Caratheodory{}\Ast pseudodistance} on a complex manifold~$\Md$ is
defined by
\begin{equation} \label{equ:defcM}
       c^\ast_\Md(\omega_1,\omega_2)
=      \sup\{ |f(\omega_2)| \mid f\in\HO(\Md,\D),\ f(\omega_1)=0 \}
\end{equation}
for $\omega_1,\omega_2\in\Md$.  The \emph{\Caratheodory-\Reiffen pseudometric}
on~$\Md$ is defined by
\begin{equation} \label{equ:defgM}
       \gamma_\Md(\omega,d)
=      \sup\{ |d(f)| \mid f\in \HO(\Md,\D),\ f(\omega)=0 \}
\end{equation}
for $(\omega,d)\in\Tg{\Md}$, i.e.\ $d\colon\CINF(\Md)\to\C$ is a derivation at
$\omega\in\Md$ ($\CINF(\Md)$ is the algebra of smooth functions $\Md\to\C$).
The \emph{\Caratheodory pseudodistance~$c_\Md$} (without~$\ast$) is related
to~$c_\Md^\ast$ by $c_\Md^\ast = {\tanh}\circ c_\Md$.  If $\Md=\D$, these
definitions yield the \emph{\Moebius distance}
$$
        m(\lambda_1,\lambda_2)
  =     c_\D^\ast(\lambda_1,\lambda_2)
  =     \left| \frac{\lambda_1-\lambda_2}{1-\lambda_1\conj{\lambda_2}} \right|,
$$
the \emph{\Poincare distance} $p=c_\D$, and the \emph{\Poincare metric}
$$
\gamma(\lambda,l) = \frac{|l|}{1-|\lambda|^2},
$$
where $(\lambda,l)$ stands for the derivation $f\mapsto lf'(\lambda)$.
See~\cite{Jarnicki-Pflug:93} for these distances and metrics.

\goodbreak

If $\NS\subset\Bound(\HilS)$ is an operator space, let
$$
\NS_{(n)}= \NS\otimes\Mat\subset\Bound(\HilS\otimes\C^n)
$$
Write~$\|\blank\|_{(n)}$ for the norm on~$\NS_{(n)}$ coming from this
representation of~$\NS_{(n)}$.

Let $\phi\colon \NS_1\to \NS_2$ be a linear map between operator spaces.
Then~$\phi$ induces linear maps $\phi_{(n)}=\phi\otimes\ID[\Mat]\colon
\NS_1\otimes\Mat\to \NS_2\otimes\Mat$.  $\phi$ is called
\emph{\PrMn{n}contractive} if~$\phi_{(n)}$ is contractive, and \emph{completely
contractive} if all the maps~$\phi_{(n)}$, $n\in\N$, are contractive.
Similarly, $\phi$ is called \emph{completely isometric} if all the
maps~$\phi_{(n)}$ are isometric, and \emph{\PrMn{n}isometric} if~$\phi_{(n)}$
is isometric.  Let
$$
\cbn{\phi} = \sup_{n\in\N} \|\phi_{(n)}\|
$$
and call~$\phi$ \emph{completely bounded} if $\cbn{\phi}<\infty$.  An isometric
map is \emph{not} required to be surjective.  A surjective, completely
isometric map is called a \emph{complete equivalence}.  If both $\NS_1$
and~$\NS_2$ are operator algebras, usually only homomorphisms are considered,
and a complete equivalence of operator algebras is a completely isometric
isomorphism.  As a matter of convention, \emph{an abstract operator
algebra~$\OA$ is called unital only if it has a unit~$e$ with $\|e\|=1$}
because if $\|e\|\neq1$ it can have no ``unital'' completely isometric
representations.

A matrix normed space or algebra satisfying certain axioms
\cite{Blecher-Ruan-Sinclair:90} is called an \emph{\LPmn space} or an
\emph{\LPmn algebra}.  The point of these axioms is that a matrix normed space
has a completely isometric (linear) representation on a Hilbert space if and
only if it is an \LPmn space, and a \emph{unital} matrix normed algebra (with
unit of norm~$1$!)  has a completely isometric (unital, multiplicative)
representation on a Hilbert space if and only if it is an \LPmn algebra.  Thus
we prefer to call them \emph{abstract operator spaces} and \emph{abstract
operator algebras}.  A not necessarily unital \LPmn algebra is only called an
abstract operator algebra if it has a completely isometric representation on a
Hilbert space.

There is a natural way to define the \emph{(orthogonal) direct sum} of operator
spaces $\NS_j\subset\HilS_j$, $j=1,2$, as $\NS_1\oplus\NS_2\subset
\Bound(\HilS_1\oplus\HilS_2)$.  Thus
\begin{equation}  \label{equ:orthoSum}
        \|(\ns_1,\ns_2)\|_{(n)}
=     \max \{\|\ns_1\|_{(n)},\|\ns_2\|_{(n)} \}
\end{equation}
for $\ns_j\in \NS_j\otimes\Mat$, $j=1,2$.  Equation~\eqref{equ:orthoSum} shows
that the resulting operator space does not depend on the chosen completely
isometric representations of~$\NS_j$, $j=1,2$.  All this remains true for
operator algebras.

Another natural construction is the \emph{quotient operator space structure}.
If~$\NS_1$ is an abstract operator space and~$\NS_2$ is a closed subspace,
define a matrix normed structure on the quotient space $\NS_1/\NS_2$ by
identifying
$$
(\NS_1/\NS_2)\otimes\Mat\cong (\NS_1\otimes\Mat)/(\NS_2\otimes\Mat)
$$
and taking the quotient norm on the latter space.  This gives an abstract
operator space \cite{Blecher-Ruan-Sinclair:90} because the result satisfies the
axioms for an \LPmn space.  Moreover, if~$\OA$ is a unital \LPmn algebra (with
unit of norm~$1$), and $\ide\subset \OA$ is a proper ideal, then $\OA/\ide$ is
again a unital \LPmn algebra and thus an abstract operator algebra.  This is
the only place where the formalism of \LPmn algebras is needed in this article.

The natural projection $\pi\colon \NS_1\to \NS_1/\NS_2$ is a \emph{complete
quotient map}, i.e.\ every map~$\pi_{(n)}$ is a quotient map.  This implies
that if $\rho\colon \NS_1/\NS_2\to\NS_3$ is any linear map, then
$\|\rho_{(n)}\|= \|(\rho\circ\pi)_{(n)}\|$ for all $n\in\N$ and in particular
$\cbn{\rho}=\cbn{\rho\circ\pi}$.

Further references for matricially normed spaces include
\cite{Blecher-Ruan-Sinclair:90}, \cite{Effros-Ruan:88}, \cite{Paulsen:92}, and
\cite{Paulsen:86}.

\section{Unitization of operator algebras}
\label{sec:Unitization}

Let $\OA\subset\Bound(\HilS)$ be an operator algebra with
$\ID[\HilS]\notin\OA$.  The goal of this section is to express the matrix
normed structure of its unitization $\Unze{\OA}=\OA\oplus \C\cdot
\ID[\HilS]$\footnote{The algebra $\OA$ may well be unital, i.e.\ have a unit of
norm~$1$, so that the notation $\Unze{\OA}$ is preferable to~$\tilde{\OA}$.}
in terms of the matrix normed structure of~$\OA$.  Thus it is independent of
the chosen representation.  The basic idea is that the open unit balls of
$\Mat$ are symmetric domains and especially have a transitive automorphism
group.  This easily implies the results of this section for uniform algebras
and Q\PrM{}algebras (without any further explicit computations!).  However, in
order to get statements for arbitrary operator algebras, some calculations are
necessary.

\begin{lemma}  \label{lem:onedimOpalg}
  Every \PrMn{1}dimensional, unital abstract operator algebra is completely
  equivalent to~$\C$ with its usual operator algebra structure.
\end{lemma}

\begin{proof}
  Let~$\OA$ be a \PrMn{1}dimensional, unital abstract operator algebra and let
  $\rho\colon \OA\to\Bound(\HilS)$ be any unital, completely isometric
  representation.  If~$1_{\OA}$ is the identity element of~$\OA$, then
  $\rho(1_{\OA})=\ID[\HilS]$ and this determines~$\rho$.  Now the lemma follows
  from the identity $\|S\otimes T\|=\|S\|\cdot\|T\|$.
\end{proof}

Define $\Cay(X)= (1-X)/(1+X)$ for all $X\in\Bound(\HilS)$ such that $1+X$ is
invertible (here $1=\ID[\HilS]$).  This differs from the usual Cayley transform
that maps a self-adjoint operator to a unitary only by some factors of~$i$ and
is indeed an analogue of the Cayley transform for skew-adjoint operators.  Let
$\Bd=\Ball\bigl(\Bound(\HilS)\bigr)$ and
$$
      \Bd_+
=     \{ X\in\Bound(\HilS) \mid \text{$\RE X$ positive and invertible} \}.
$$

\begin{lemma}  \label{lem:Cayley}
  $\Cay(X)$ is well\PrM{}defined for $X\in\Bd\cup\Bd_+$.  $\Cay$ maps~$\Bd$
  bijectively onto~$\Bd_+$ and is its own inverse, i.e.\ 
  $\Cay\circ\Cay(X)=X$ for $X\in\Bd\cup\Bd_+$.  More generally, $1+\Cay(X)$ is
  invertible whenever $\Cay(X)$ is defined and then $\Cay\circ\Cay(X)=X$.
\end{lemma}

\begin{proof}
  It is easy to see that elements of $1+\Bd$ and $1+\Bd_+$ are invertible.  The
  identity $\Cay\circ\Cay(X)=X$ is easy to check, whenever the left side is
  well\PrM{}defined.  Thus it only remains to show $\Cay(\Bd)\subset \Bd_+$ and
  $\Cay(\Bd_+)\subset\Bd$.
  \begin{align*}
      2\RE \Cay(X)
  &=  (1+X)^{-1}\bigl( (1-X)(1+X^\ast) + (1+X)(1-X^\ast) \bigr)(1+X^\ast)^{-1}
\\&=  2(1+X)^{-1}( 1 -XX^\ast )(1+X^\ast)^{-1}
  \end{align*}
  shows that $\RE \Cay(X)$ is positive and invertible for $X\in\Bd$, i.e.\
  $\Cay(\Bd)\subset\Bd_+$.  A similar calculation shows that for $X\in\Bd_+$,
  $$
  1-\Cay(X)\Cay(X)^\ast = 4(1+X)^{-1}\RE(X) (1+X^\ast)^{-1}
  $$
  is positive and invertible, so that $\Cay(X)\in\Bd$.
\end{proof}

\begin{theorem}  \label{the:posCone}
  Let~$\OA$ be a closed, unital abstract operator algebra and $X\in\OA_{(n)}$.
  Then the following are equivalent:

  \begin{enumerate}[(i)]%

  \item $X$ is of the form~$\Cay(Y)$ for some $Y\in\Ball(\OA_{(n)})$;

  \item $1+X$ is invertible in~$\OA_{(n)}$ and $\Cay(X)\in\Ball(\OA_{(n)})$;

  \item $\rho_{(n)}(X)$ has invertible and positive real part for all
    \PrMn{n}contractive, unital representations $\rho\colon\OA\to\Bound(\HilS)$;

  \item $\rho_{(n)}(X)$ has invertible and positive real part for some
    \PrMn{n}isometric, unital representation $\rho\colon\OA\to\Bound(\HilS)$.

  \end{enumerate}
\end{theorem}

\begin{definition}
  The set of elements satisfying one of these equivalent conditions is called
  the \emph{(positive) cone} $\Cone(\OA_{(n)})$ of~$\OA_{(n)}$.
\end{definition}

\begin{proof}
  Replacing~$\OA$ by~$\OA_{(n)}$, if necessary, it can be assumed that $n=1$,
  i.e.\ $X\in\OA$.  Since~$\OA$ is closed, all elements of $1+\Ball(\OA)$ are
  invertible in~$\OA$, so that~$\Cay(Y)$ is defined and lies in~$\OA$ for all
  $Y\in\Ball(\OA)$.  Thus the equivalence of (i) and~(ii) follows as in
  Lemma~\ref{lem:Cayley}.  Moreover, (iii) trivially implies~(iv).

  If~$\rho$ is some contractive ($n=1$) representation and $X=\Cay(Y)$ with
  $Y\in\Ball(\OA)$, then $\rho(X)=\Cay\bigl(\rho(Y)\bigr)$ because~$\rho$ is a
  homomorphism, and this lies in~$\Bd_+$ by Lemma~\ref{lem:Cayley}.  Hence~(i)
  implies~(iii), so that it remains to show that~(iv) implies~(ii).

  Therefore, take any isometric unital representation
  $\rho\colon\OA\to\Bound(\HilS)$ and $X\in\OA$ such that~$\rho(X)$ has
  positive and invertible real part.  It is not difficult to see that
  $A\in\Bound(\HilS)$ is invertible if its real part is positive and
  invertible.  Thus $\rho(X)+\lambda\ID[\HilS]$ is invertible in
  $\Bound(\HilS)$ for all $\lambda\in\C$ with $\RE\lambda\ge0$.  Therefore, the
  spectrum of $1+\rho(X)$ is a compact subset of $\{\lambda\in\C\mid
  \RE\lambda>1\}$.  The function $\lambda\mapsto 1/\lambda$ can be approximated
  uniformly on a neighborhood of this compact set by polynomials.  Thus the
  inverse of $1+\rho(X)$ lies in~$\rho(\OA)$, so that $1+X$ is invertible
  in~$\OA$ and not just in $\Bound(\HilS)$.  Hence~$\Cay(X)$ is a
  well\PrM{}defined element of~$\OA$.  Moreover,
  $\rho\bigl(\Cay(X)\bigr)=\Cay\bigl(\rho(X)\bigr)$.  Since
  $\|\Cay\bigl(\rho(X)\bigr)\|<1$ by Lemma~\ref{lem:Cayley} and~$\rho$ is
  isometric, $\|\Cay(X)\|<1$.  Thus~(iv) implies~(ii).
\end{proof}

Due to this correspondence, the matrix normed structure of a closed unital
operator algebra $\OA\subset\Bound(\HilS)$ can equally well be described by its
positive cone.  This is especially advantageous if~$\HilS$ does not come with
an orthonormal basis but just with a \emph{frame}.  A set of vectors
$(\xi_j)_{j\in J}$ in a Hilbert space~$\HilS$ is a frame if there exist numbers
$A,B\in(0,\infty)$ such that, for all $\eta\in\HilS$,
$$
A\|\eta\|^2 \le \sum_{j\in J} |\5{\xi_j}{\eta}|^2 \le B\|\eta\|^2.
$$
Thus a frame need not be linearly independent.  With every frame one can
associate a bounded linear map $S\colon \HilS\to\ell^2(J)$ mapping~$\eta$ to
$(\5{\eta}{\xi_j})_{j\in J}$.  The assumptions guarantee that~$S$ is bounded
and that $S^\ast S$ is invertible.

\begin{proposition}  \label{pro:frame}
  Let~$\HilS$ be a Hilbert space, $(\xi_j)_{j\in J}$ a frame, and
  $T\in\Bound(\HilS)$.  Then~$T$ has positive real part iff $S^\ast
  TS\in\Bound\bigl(\ell^2(J)\bigr)$ has positive real part.  This happens iff
  the matrix $\tilde{T}\in\Bound\bigl(\ell^2(J)\bigr)$ with entries
  $$
	\tilde{T}_{ij}=\5{T\xi_j}{\xi_i}+\5{\xi_j}{T\xi_i}
  $$
  is positive definite.

  If~$S$ is invertible, then~$T$ has positive and invertible real part iff
  $S^\ast TS$ has positive and invertible real part.
\end{proposition}

\begin{proof}
  If~$T$ has positive real part, then so has $STS^\ast\in
  \Bound\bigl(\ell^2(J)\bigr)$ because
  \begin{equation}  \label{equ:frameTT}
  2\RE(STS^\ast) = STS^\ast + ST^\ast S^\ast = 2S\RE(T)S^\ast.
  \end{equation}
  Conversely, if $STS^\ast$ has positive real part then so has $S^\ast STS^\ast
  S\in\Bound(\HilS)$ and hence~$T$ because $S^\ast S$ is invertible.  If~$S$ is
  invertible then $STS^\ast$ has invertible real part iff~$T$ has
  by~\eqref{equ:frameTT}.

  Let $\{e_j\}_{j\in J}$ be the canonical orthonormal basis of~$\ell^2(J)$.
  Then $S^\ast e_j=\xi_j$ and thus
  $$
	\5{2\RE(STS^\ast) e_j}{e_i}
  =	\5{2S\RE(T)S^\ast e_j}{e_i}
  =	\5{2\RE(T)\xi_j}{\xi_i}
  =	\5{T\xi_j}{\xi_i}+\5{\xi_j}{T\xi_i}.
  $$
  Hence the matrix~$\tilde{T}$ comes from the operator $2\RE(STS^\ast)$.
\end{proof}

The situation of Proposition~\ref{pro:frame} will occur in
Section~\ref{sec:QuotientMultiplier}: The Hilbert spaces on which quotients of
$\Shift_d$ are represented do not come with natural orthonormal bases, but
linearly independent frames are rather easy to obtain.

It is easy to write down automorphisms of $\Cone(\OA_{(n)})$: If $A\in\Mat$ is
invertible and $B\in\Mat$ satisfies $\RE B=0$, then $\Phi_{A,B}\colon X\mapsto
AXA^\ast + B$ defines a bijection from~$\Cone(\OA_{(n)})$ onto itself, with
inverse $X\mapsto A^{-1}X(A^{-1})^\ast -A^{-1}B(A^{-1})^\ast$.  It is easy to
see that these maps really map $\Cone(\OA_{(n)})$ into itself using the
characterization~(iv) of Theorem~\ref{the:posCone}.  Consequently, the map
$\Cay\circ\Phi_{A,B}\circ\Cay$ gives a bijection
$\Ball(\OA_{(n)})\to\Ball(\OA_{(n)})$.

\begin{theorem}  \label{the:unitizationWellDefined}
  Let $\OA\subset\Bound(\HilS)$ be a closed operator algebra not containing
  $\ID[\HilS]$.  Then the matrix normed structure of~$\Unze{\OA}$ is uniquely
  determined by the matrix normed structure of~$\OA$: The open unit ball of
  $(\Unze{\OA})_{(n)}$ consists precisely of the elements $\Psi(X)$ where $X\in
  \Ball(\OA_{(n)})$ and $\Psi=\Cay\circ\Phi_{A,B}\circ\Cay$ for $A,B\in\Mat$,
  $A$ invertible and $\RE B = 0$.

  If $\rho\colon \OA\to\Bound(\HilS')$ is a completely
  isometric representation and $\ID[\HilS']\notin \rho(\OA)$, then the
  unitization $\Unze{\rho}\colon \Unze{\OA}\to\Bound(\HilS')$, defined by
  $\Unze{\rho}|_{\OA}=\rho$, $\Unze{\rho}(1_{\Unze{\OA}})=\ID[\HilS']$, is also
  completely isometric.

  More generally, let $\rho\colon \OA\to\Bound(\HilS')$ be a representation
  with $\ID[\HilS']\notin \rho(\OA)$.  Then~$\Unze{\rho}$ is
  \PrMn{n}contractive if and only if~$\rho$ is \PrMn{n}contractive.  The
  map~$(\Unze{\rho})_{(n)}$ is a quotient map onto its image if and only
  if~$\rho_{(n)}$ is a quotient map onto its image.
\end{theorem}

\begin{proof}
  All elements $\Psi(X)$ with $X\in\Ball(\OA_{(n)})$ and~$\Psi$ as above lie in
  $\Ball\bigl((\Unze{\OA})_{(n)}\bigr)$.  Conversely, let~$X$ in
  $\Ball\bigl((\Unze{\OA})_{(n)}\bigr)$.  The map $\Unze{\OA}\to
  \Unze{\OA}/\OA$ is a completely contractive homomorphism.  Since~$\Unze{\OA}$
  is unital, so is $\Unze{\OA}/\OA$, so that $\Unze{\OA}/\OA$ is completely
  equivalent to~$\C$ with the usual matrix normed structure by
  Lemma~\ref{lem:onedimOpalg}.  This yields a completely contractive unital
  homomorphism $\omega\colon\Unze{\OA}\to\C$.  Thus $X_0=\omega_{(n)}(X)$ lies
  in $\Ball(\Mat)$.
  
  For suitable $A,B\in\Mat$, the corresponding
  $\Psi=\Cay\circ\Phi_{A,B}\circ\Cay\colon \Ball(\Mat)\to\Ball(\Mat)$ satisfies
  $\Psi(0)=X_0$.  Let~$\Psi$ also denote the map
  $$
\Cay\circ\Phi_{A,B}\circ\Cay\colon \Ball\bigl((\Unze\OA)_{(n)}\bigr)
	\to\Ball\bigl((\Unze\OA)_{(n)}\bigr).
  $$
  Then $\omega_{(n)}\circ\Psi^{-1}=\Psi^{-1}\circ\omega_{(n)}$ and
  $\omega_{(n)}\bigl(\Psi^{-1}(X)\bigr)=0$, i.e.\ $\Psi^{-1}(X)\in \OA_{(n)}$.
  Since also $\|\Psi^{-1}(X)\|<1$ and $X=\Psi\bigl(\Psi^{-1}(X)\bigr)$, all
  elements of $\Ball\bigl((\Unze{\OA})_{(n)}\bigr)$ are of the form $\Psi(Y)$
  for some $Y\in\Ball(\OA_{(n)})$.

  The definition of~$\Psi(X)$ is purely algebraic and does not use the matrix
  normed structure of~$\Unze{\OA}$, but only that certain elements
  of~$\Unze{\OA}$, namely those of the form $1+\Phi_{A,B}\circ\Cay(X)$ with
  appropriate $X,A,B$, are invertible.  Thus any completely isometric
  representation of~$\OA$ whose image does not contain the identity, yields the
  same unit ball for~$(\Unze{\OA})_{(n)}$.  Hence the unital extension of such
  a completely isometric representation of~$\OA$ is completely isometric.

  The map~$\rho_{(n)}$ is contractive iff it maps $\Ball(\OA_{(n)})$ into
  $\Ball\bigl(\Bound(\HilS'\otimes\C^n)\bigr)$, and similarly
  for~$\Unze{\rho}_{(n)}$.  Thus contractiveness of~$\Unze{\rho}_{(n)}$
  trivially implies contractiveness of~$\rho_{(n)}$.  The converse follows
  because $\Unze{\rho}_{(n)}\bigl(\Psi(X)\bigr)=\Psi\bigl(\rho_{(n)}(X)\bigr)$
  for all $X\in\Ball(\OA_{(n)})$ and the right side lies in
  $\Ball\bigl(\Bound(\HilS'\otimes\C^n)\bigr)$ if $\|\rho_{(n)}(X)\|<1$.

  The statement about quotient maps follows in the same way, using
  that~$\rho_{(n)}$ is a quotient map onto its image iff it maps
  $\Ball(\OA_{(n)})$ onto $\Ball\bigl(\Bound(\HilS'\otimes\C^n)\bigr)\cap
  \rho_{(n)}(\OA_{(n)})$ and the same statement for~$\Unze{\rho}_{(n)}$.
\end{proof}

For the special case of quotients of the function algebra $\BH(\D)$, this
corresponds to well\PrM{}known facts of Nevanlinna-Pick theory.
Theorem~\ref{the:posCone} means that interpolation problems from the unit disk
to~$\Mat$ with values in $\Ball(\Mat)$ and $\Cone(\Mat)$ are equivalent.
Theorem~\ref{the:unitizationWellDefined} means that one can always restrict to
the case where the origin is mapped to $0\in\Ball(\Mat)$ or $1\in\Cone(\Mat)$.

If $\OA\subset\Bound(\HilS)$ does contain $\ID[\HilS]$, consider
$\OA\subset\Bound(\HilS\oplus\C)$ in the obvious way and define
$\Unze{\OA}\subset\Bound(\HilS\oplus\C)$.
Theorem~\ref{the:unitizationWellDefined} shows that~$\Unze{\OA}$ is a
well\PrM{}defined abstract operator algebra if~$\OA$ is an abstract operator algebra
(it does not matter whether~$\OA$ is closed), i.e.\ it does not depend on the
choice of a completely isometric representation.  Every homomorphism
$\rho\colon \OA_1\to\OA_2$ can be extended uniquely to a unital homomorphism
$\Unze{\rho}\colon \Unze{\OA_1}\to \Unze{\OA_2}$.

\begin{corollary}  \label{cor:unitizationFunctor}
  $\OA\to\Unze{\OA}$ is a functor from the category of abstract operator
  algebras with completely contractive homomorphisms as morphisms to the
  category of unital abstract operator algebras with unital, completely
  contractive homomorphisms as morphisms.  Furthermore, it maps complete
  quotient maps to complete quotient maps and complete isometries to complete
  isometries.  Thus $\Unze{(\OA/\ide)}\cong \Unze{\OA}/\ide$.
\end{corollary}

Let~$\NS$ be an abstract operator space, furnish it with the zero
multiplication.  This yields an abstract operator algebra: Let $\rho\colon
\NS\to\Bound(\HilS)$ be any completely isometric linear representation, define
$\rho_2\colon \NS\to \Bound(\HilS\oplus\HilS)$ by putting
$$
\rho_2(x) = \begin{pmatrix} 0 & \rho(x) \\ 0 & 0 \end{pmatrix}.
$$
Then $\rho_2(x)\rho_2(y)=0$ for all $x,y\in \NS$, i.e.~$\rho_2$ is
multiplicative.  Now let~$\Unze{\NS}$ be the unitization of~$\NS$, viewed as an
abstract operator algebra with zero multiplication.  $\Unze{\NS}$ is called the
\emph{trivial unitization of~$\NS$}.  Here trivial does not mean that the
resulting objects cannot be very complicated but that~$\NS$ is endowed with the
trivial multiplication.

It is easy to see that every unital homomorphism
$\sigma\colon\Unze{\NS_1}\to\Unze{\NS_2}$ is of the form $\sigma=\Unze\rho$ for
some $\rho\colon \NS_1\to \NS_2$, namely $\rho=\sigma|_{\NS_1}$.  $\NS$~is the
unique maximal ideal of~$\Unze{\NS}$, and the multiplication on~$\NS$ is
trivial.  Conversely, if~$\OA$ is an abstract operator algebra with a
\PrMn{1}codimensional ideal~$\NS$ such that the multiplication on~$\NS$ is
trivial, then $\OA=\Unze{\NS}$.  Especially, the image of a unital homomorphism
$\OA\to\Bound(\HilS)$ again has these properties.

\begin{theorem}  \label{the:unitizationCB}
  Let $\NS_1$ and~$\NS_2$ be abstract operator spaces, let $\rho\colon
  \NS_1\to\NS_2$ be a linear map, and let $\Unze\rho\colon
  \Unze{\NS_1}\to\Unze{\NS_1}$ be its unitization.  Then, for all $n\in\N$,
  $$
        \|\Unze\rho_{(n)}\|
  =     \max\{1,\|\rho\|_{(n)}\}.
  $$
\end{theorem}

\begin{proof}
  The inequality ``$\ge$'' is trivial.  To prove ``$\le$'', assume
  $M=\|\rho_{(n)}\|<\infty$.  If $M\le1$, the assertion follows from
  Theorem~\ref{the:unitizationWellDefined}.  Thus assume $M>1$ and let $m\colon
  \NS_1\to\NS_1$ be the map $T\mapsto M T$.  Then $\rho=\rho\circ m^{-1}\circ
  m$, and $\|(\rho\circ m^{-1})_{(n)}\|=1$.  Hence $\|(\Unze\rho\circ
  \Unze{(m^{-1})})_{(n)}\|\le 1$, so that it remains to prove
  $\cbn{\Unze{m}}\le M$.
  
  Therefore, choose a representation of~$\NS_1$ on some Hilbert space~$\HilS$
  and view $\Unze{\NS_1}\subset\Bound(\HilS\oplus\HilS)$ as above.  Let
  $S=\diag(M^{1/2},M^{-1/2})$ be a diagonal matrix with respect to this
  decomposition, i.e.\ $S=M^{1/2}\ID[\HilS]$ on the first copy of~$\HilS$ and
  $S=M^{-1/2}\ID[\HilS]$ on the second copy.  Then $\Unze{m}(T)=STS^{-1}$ for
  all $T\in\Unze{\NS_1}$ because both sides of this equation are unital maps
  that coincide on~$\NS_1$.  Thus $\cbn{\Unze{m}}\le \|S\|\cdot \|S^{-1}\| = M$
  as desired.
\end{proof}

Any linear representation $\rho\colon \NS\to\Bound(\HilS)$ yields a unital,
multiplicative representation $\Unze{\rho}\colon \Unze{\NS} \to
\Bound(\HilS\oplus\HilS)$ with the same boundedness properties.  Thus the
representation theory of~$\Unze{\NS}$ is precisely as well\PrM{}behaved (or
pathological) as the linear representation theory of~$\NS$.

There seems to be no analogue of Theorem~\ref{the:unitizationCB}
for unitizations of homomorphisms between operator algebras with non-trivial
multiplication.  Already for the two\PrM{}dimensional examples, the estimate
$\|\Unze\rho\|\le \max\{1,\cbn{\rho}\}$ does not hold.  Instead,
\eqref{equ:iotaccNorm} implies
$$
\|\Unze\rho\| = \cbn{\Unze\rho} \le \max\{1,2\|\rho\|\}.
$$
But this probably is a special feature of the two\PrM{}dimensional case.

\begin{theorem}  \label{the:genSmith}
  Let~$\OA$ be a unital, commutative operator algebra, $d\in\N$.  Then any
  \PrM{d-1}contractive unital homomorphism $\rho\colon\OA\to\Mat[d]$ is
  completely contractive.
\end{theorem}

\begin{proof}
  Let $\OA[B]=\rho(\OA)$, let~$\OA[J]$ be a maximal ideal in~$\OA[B]$, and
  $\OA[I]=\rho^{-1}(\OA[J])$.  Then $\OA[A]=\Unze{\OA[I]}$,
  $\OA[B]=\Unze{\OA[J]}$, and $\rho=\Unze{(\rho|_{\OA[I]})}$.  By
  Theorem~\ref{the:unitizationWellDefined}, it suffices to show
  that~$\rho|_{\OA[I]}$ is completely contractive.

  By abstract algebra, there exists a vector $x\in\C^d$ that is annihilated by
  all elements of~$\OA[J]$, because~$\OA[J]$ is a non-unital, commutative
  operator algebra.  Apply the same reasoning to the algebra~$\OA[J]^\ast$ of
  all adjoints of elements of~$\OA[J]$.  This yields $y\in\C^d$ that is
  annihilated by all adjoints~$T^\ast$ of elements $T\in\OA[J]$.

  Thus elements of~$\OA[J]$ can be viewed as operators from $\C^d\ominus x$ to
  $\C^d\ominus y$.  This yields a completely isometric \emph{linear}
  representation of~$\OA[J]$ by \PrM{(d-1)\times (d-1)}matrices.  In general,
  this representation fails to be a homomorphism, but this does not matter.  By
  a result of Smith (see~\cite{Paulsen:86}), every linear representation
  $\phi\colon \NS\to\Mat[n]$ of an operator space~$\NS$ satisfies
  $\cbn{\phi}=\|\phi_{(n)}\|$.  Application of Smith's theorem to~$\rho|_\ide$,
  with $n=d-1$, yields the assertion.
\end{proof}

The following obvious corollary will be strengthened in the next section.

\begin{corollary}  \label{cor:genSmithTwo}
  Let~$\OA$ be a unital, commutative operator algebra.  Then any contractive
  unital homomorphism $\rho\colon\OA\to\Mat[2]$ is completely contractive.
\end{corollary}

\section{Two-dimensional unital operator algebras}
\label{sec:twodimOpalg}

In this section, two\PrM{}dimensional unital operator algebras are classified
up to complete equivalence.  This is the smallest non-trivial dimension by
Lemma~\ref{lem:onedimOpalg}.  Most results will be false for
two\PrM{}dimensional non-unital operator algebras and for operator algebras of
dimension at least~$3$.  This makes the study of two\PrM{}dimensional unital
operator algebras even more important because it is the only case with a
satisfactory theory.  Theorem~\ref{the:unitizationWellDefined} reduces the
classification of two\PrM{}dimensional unital operator algebras~$\OA$ to that
of \PrMn{1}dimensional operator algebras, which is rather trivial.

\begin{lemma}  \label{lem:onedimOpspace}
  Every \PrMn{1}dimensional abstract operator space is completely equivalent
  to~$\C$ with its usual matrix normed structure.
\end{lemma}

\begin{proof}
  Let~$\NS$ be a \PrMn{1}dimensional abstract operator space, choose a
  completely isometric linear representation $\NS\subset\Bound(\HilS)$ on some
  Hilbert space~$\HilS$.  Choose $\ns\in\NS$ with $\|\ns\|=1$.  Then every
  element of $\NS_{(n)}$ can be written as $\ns\otimes T$ for some
  $T\in\Mat$.  Since~$\rho$ is completely isometric,
  $$
  \|\ns\otimes T\|=\|\rho(\ns)\otimes T\|= \|\rho(\ns)\|\cdot\|T\|=\|T\|.
  $$
\end{proof}

For $c\in[0,1]$, let
$$
T_c = \begin{pmatrix} 0 & \sqrt{1-c^2} \\ 0 & c \end{pmatrix}.
$$
Clearly, $\|T_c\|=1$.  Since $T_c^2=cT_c$, $\C\cdot T_c$ is an operator
algebra.  Let $\QA_c=\lin\{1,T_c\}$ be the \emph{unital} operator algebra
generated by~$T_c$.

\begin{theorem}  \label{the:onedimOpalgNU}
  Let~$\OA$ be a \PrMn{1}dimensional abstract operator algebra.  Then there
  exist unique $c\in[0,1]$ and $T\in\OA$ such that $T\mapsto T_c$ defines a
  completely isometric representation of~$\OA$.
\end{theorem}

\begin{proof}
  There exist unique $T\in\OA$ and $c\in[0,\infty]$ with $\|T\|=1$ and
  $T^2=cT$: The first condition determines~$T$ uniquely up to scalar
  multiplication by some $\lambda\in\bdD$.  This can be used to make the
  constant~$c$ real and nonnegative.  In addition, $c=\|cT\|= \|T^2\|\le
  \|T\|^2=1$, so that $c\in[0,1]$.  Using Lemma~\ref{lem:onedimOpspace}, it is
  easy to see that $T\mapsto T_c$ defines a completely isometric
  representation of~$\OA$.
\end{proof}

\begin{theorem}  \label{the:twodimOpalg}
  Let~$\OA$ be a two\PrM{}dimensional unital operator algebra.  Then~$\OA$ is
  completely equivalent to~$\QA_c$ for a unique $c\in[0,1]$ and thus has a
  completely isometric, unital representation by \PrMn{2\times2}matrices.
  Indeed, the algebras~$\QA_c$ are not even isometrically isomorphic for
  different values of~$c$.
\end{theorem}

\begin{proof}
  Let~$\OA$ be a \PrMn{2}dimensional, unital abstract operator algebra.
  Then~$\OA$ is necessarily commutative and thus contains a maximal
  ideal~$\ide$.  Thus $\OA=\Unze{\ide}$.  By Theorem~\ref{the:onedimOpalgNU},
  there is a completely isometric representation $\ide\to\Mat[2]$ whose image
  does not contain the identity.  The unitization of this representation gives
  a unital, completely isometric representation of~$\OA$ by
  \PrM{2\times2}matrices.

  It will follow immediately from Theorem~\ref{the:twodimAut} that the
  algebras~$\QA_c$ for different values of~$c$ are not isometrically
  isomorphic.
\end{proof}

The following theorem lists all (algebraic) isomorphisms, i.e.\ bijective
algebra homomorphisms, between the algebras~$\QA_c$, $c\in[0,1]$.  Such an
isomorphism is necessarily unital and hence determined by the image of the
generator~$T_c$.

\begin{theorem}  \label{the:twodimAut}
  \begin{enumerate}[(i)]
  \item The automorphisms of~$\QA_0$ are the unital maps
    $m_\lambda\colon\QA_0\to\QA_0$ defined by $T_0\mapsto \lambda T_0$ for
    $\lambda\in\C^\ast$.  There are no isomorphisms $\QA_0\to\QA_c$ (or in the
    opposite direction) for $c\in(0,1]$.  The only non-identical automorphism
    of~$\QA_c$ for $c\neq0$ is $\theta_c\colon\QA_c\to\QA_c$ given by
    $$
    T_c \mapsto c - T_c =
    \begin{pmatrix}
      c & -\sqrt{1-c^2} \\ 0 & 0
    \end{pmatrix}.
    $$
    For $c,c'\in(0,1]$, there are two isomorphisms $\QA_c\to\QA_{c'}$,
    namely~$\iota_{c,c'}$ defined  by $\iota_{c,c'}(T_c)= (c/c')T_{c'}$ and
    $\theta_{c'}\circ\iota_{c,c'}=\iota_{c,c'}\circ\theta_c$.
    
  \item These isomorphisms are contractive if and only if they are completely
    contractive and isometric if and only if they are completely isometric.
    
  \item The automorphism~$m_\lambda$ is (completely) contractive iff
    $|\lambda|\le1$ and (completely) isometric iff $|\lambda|=1$.  More
    generally,
    \begin{equation}  \label{equ:mlambdaNorm}
          \|m_\lambda\|
    =     \cbn{m_\lambda}
    =     \max \{1,|\lambda|\}.
    \end{equation}
  
  \item The automorphism~$\theta_c$ is completely isometric for all
    $c\in(0,1]$.  The isomorphism~$\iota_{c,c'}$ is (completely) contractive
    iff $c\le c'$ and isometric only for $c=c'$.  More generally,
    \begin{equation} \label{equ:iotaccNorm}
        \|\iota_{c,c'}\|
    =   \cbn{\iota_{c,c'}}
    =   \max \{1, h(c')/h(c)\},
    \end{equation}
    where
    $$h(c) = c^{-1}(1+\sqrt{1-c^2}).$$
  \end{enumerate}
\end{theorem}

\begin{proof}
  \begin{enumerate}[(i)]
  \item A unital map $\ell\colon \QA_c\to\QA_{c'}$ is an isomorphism iff
    $\ell(T_c)^2 = c\ell(T_c)$.  This easily implies that the maps
    $m_\lambda$, $\theta_c$, and $\iota_{c,c'}$ are isomorphism, and that
    there are no other possibilities.
    
  \item Theorem~\ref{the:unitizationWellDefined} and
    Lemma~\ref{lem:onedimOpspace} imply that an isomorphism $\ell\colon
    \QA_c\to\QA_{c'}$ is (completely) contractive iff its restriction to the
    linear span of~$T_c$ is (completely) contractive iff $\|\ell(T_c)\|\le1$,
    and that it is (completely) isometric iff $\|\ell(T_c)\|=1$.
    
  \item The restriction of~$m_\lambda$ to the linear span of~$T_0$ is just
    scalar multiplication by~$\lambda$ and thus has norm and complete
    norm~$|\lambda|$.  Hence the claim is a special case of
    Theorem~\ref{the:unitizationCB}.

  \item Since $\|\theta_c{T_c}\|=1$ and $\|\iota_{c,c'}(T_c)\|= c/c'$, the
    first part of the assertion follows from the proof of~(ii).  It remains to
    prove~\eqref{equ:iotaccNorm}.
    
    For $c\in(0,1]$, let $\tilde{T}_c=-1+ 2c^{-1}T_c$.  The
    matrix~$\tilde{T}_c$ has the eigenvalues $\pm1$ and thus is a more
    ``symmetric'' generator for~$\QA_c$.  Moreover,
    $\iota_{c,c'}(\tilde{T}_c)=\tilde{T}_{c'}$.  An elementary calculation
    shows $\|\tilde{T}_c\|=h(c)$.  This is, of course, where the function~$h$
    comes from.  Therefore,
    $$
        \cbn{\iota_{c,c'}}
    \ge \|\iota_{c,c'}\|
    \ge \max \{1, h(c')/h(c)\},
    $$
    so that it remains to prove $\cbn{\iota_{c,c'}} \le\max \{1, h(c')/h(c)\}$.
    This estimate is true for $c\le c'$ because then~$\iota_{c,c'}$ is
    completely contractive.
    
    Thus it only remains to show $\cbn{\iota_{c,c'}} \le h(c')/h(c)$ if
    $0<c'<c\le1$.  This follows if there is an invertible $S\in\Mat[2]$ with
    $\iota_{c,c'}(T)= STS^{-1}$ for all $T\in\QA_c$ and $\|S\|\cdot \|S^{-1}\|=
    h(c')/h(c)$.  Let
    $$
    S_{\gamma}=\begin{pmatrix} 1 & \sqrt{1-\gamma^2} \\ 0 & \gamma \end{pmatrix}
    $$
    for $\gamma\in(0,1]$ and $S=S_{c'}S_c^{-1}$.  It is easy to check
    $S_\gamma^{-1}T_\gamma S_\gamma =\gamma T_1$, so that
    $ST_cS^{-1}=(c/c')T_{c'}$.  Thus $STS^{-1}=\iota_{c,c'}(T)$ for all
    $T\in\QA_c$.
    
    The computation of $\|S\|\cdot\|S^{-1}\|$ is quite elementary but tedious,
    so that the details are left to the reader.  A main step is to compute
    $$
      \left\| \begin{pmatrix} 1 & y \\ 0 & x \end{pmatrix} \right\|
      \cdot
      \left\| \begin{pmatrix} 1 & y \\ 0 & x \end{pmatrix}^{-1} \right\|
    = \frac{1}{2|x|}\bigl( 1+|x|^2+|y|^2 + \sqrt{(1+|x|^2+|y|^2)^2
       - 4|x|^2}\bigr)
    $$
    for $x\in\C^\ast$, $y\in\C$.  This can be applied to the matrix~$S$.
    Using the rules
    $$
        \sqrt{2-c^2-{c'}^2-2\sqrt{1-c^2}\sqrt{1-{c'}^2}}
    =   \sqrt{1-{c'}^2}-\sqrt{1-c^2} \\
    $$
    and
    $$
        h(c)^{-1}= c^{-1}(1-\sqrt{1-c^2}),
    $$
    the result can be transformed into $h(c')/h(c)$.
  \end{enumerate}
\end{proof}

See also \cite{Fu-Russo:95} and~\cite{Meyer:96} for different looking but
equivalent versions of \eqref{equ:mlambdaNorm} and~\eqref{equ:iotaccNorm}.  The
discovery that $\|\rho\|=\cbn{\rho}$ for homomorphisms between the
algebras~$\QA_c$ goes back to Holbrook~\cite{Holbrook:77}.

\begin{corollary}  \label{cor:twodimOpalgDisk}
  Every two\PrM{}dimensional unital operator algebra is completely equivalent
  to a quotient of the disk algebra $\Pol(\clD)$.
\end{corollary}

\begin{proof}
  By Theorem~\ref{the:twodimOpalg}, it suffices to show that each of the
  algebras~$\QA_c$ is completely equivalent to a quotient of the disk algebra.
  The quotient algebra $\Pol(\clD)/\ide(0)^2$ is two\PrM{}dimensional and
  algebraically not isomorphic to~$\QA_c$ for $c\neq0$.  Hence $\QA_0\cong
  \Pol(\clD)/\ide(0)^2$ by Theorem~\ref{the:twodimOpalg}.
  
  Let $c\in(0,1]$.  The spectrum of~$\QA_c$ consists of precisely two points,
  called $0$ and~$c$, respectively, such that $T_c(0)=0$, $T_c(c)=c$.
  Clearly,
  \begin{equation}  \label{equ:c}
        \sup\{|T(c)| \mid T\in \QA_c,\|T\|\le1,T(0)=0\} = c.
  \end{equation}
  It is elementary that
  $$
  \sup \{ |f(c)| \mid f\in\Pol(\clD), f(0)=0 \} = c.
  $$
  By definition of the quotient norm, $\Pol(\clD)$ can be replaced by
  $\Pol(\clD)/\ide(0,c)$ here.  Comparison with~\eqref{equ:c} shows
  $\Pol(\clD)/\ide(0,c)\not\cong \QA_{c'}$ for $c'\neq c$.  Thus, by
  Theorem~\ref{the:twodimOpalg}, $\Pol(\clD)/\ide(0,c)\cong \QA_c$.
\end{proof}

The easier parts of Theorem~\ref{the:twodimAut} can be obtained from this
realization of the algebras~$\QA_c$ as Q\PrM{}algebras.
$\theta_c[f]=[f\circ\Shf{c}]$, where $\Shf{c}\in\Aut(\D)$ is characterized by
$\Shf{c}(c)=0$ and $\Shf{c}(0)=c$.  The map $f\mapsto f\circ\Shf{c}$ is an
automorphism of $\clPol(\clD)$ and obviously completely isometric.  By
definition of the quotient norm, this passes down to~$\QA_c$.  Similarly, the
complete contractiveness of $m_\lambda$ and~$\iota_{c,c'}$, for appropriate
$\lambda,c,c'$, can be deduced by lifting these maps to endomorphisms $f\mapsto
f\circ g$ of $\Pol(\clD)$, where $g(z)=\lambda z$.  Moreover, the fact that
every contractive homomorphism $\QA_c\to\blank$ is completely contractive
follows from the \SzNagy dilation theorem \cite{Nagy-Foias:70}.

\begin{corollary}  \label{cor:contr2}
  Let $\OA_1$ and~$\OA_2$ be unital operator algebras and let $\rho\colon
  \OA_1\to \OA_2$ be a unital homomorphism, which has rank at most~$2$ as a
  linear map (for example, $\dim \OA_1=2$ or $\dim \OA_2=2$).  If~$\rho$ is
  contractive, then it is completely contractive.  More generally,
  $\|\rho\|=\cbn{\rho}$.
\end{corollary}

\begin{proof}
  Replace~$\OA_1$ by the unital operator algebra $\OA_1/\Ker\rho$ and
  replace~$\OA_2$ by the unital operator algebra $\Ran\rho$.  This does not
  change the norm or complete norm of~$\rho$.  Thus both $\OA_1$ and~$\OA_2$
  may be assumed to have dimension at most~$2$.

  The one\PrM{}dimensional case follows immediately from
  Lemma~\ref{lem:onedimOpalg}.  The two\PrM{}dimensional case follows from
  Theorem~\ref{the:twodimOpalg} and Theorem~\ref{the:twodimAut}.
\end{proof}

A main tool in Agler's proof of \Lempert's theorem \cite{Agler:90} is that for
$\Md\subset\C^n$, contractive representations of $\BH(\Md)$ by
\PrM{2\times2}matrices are always completely contractive.  This result about
representations by \PrM{2\times2}matrices has since then been more and more
generalized \cite{Salinas:91}, \cite{Chu:92}, \cite{Meyer:96}.  The following
should be the most general form possible.

\begin{corollary}
  Let~$\OA$ be a commutative, unital operator algebra.  Then every unital,
  contractive representation of~$\OA$ by \PrM{2\times2}matrices is completely
  contractive.  More generally, if $\rho\colon \OA\to\Mat[2]$ is a unital
  homomorphism, then $\cbn{\rho}=\|\rho\|$.
\end{corollary}

\begin{proof}
  It is not difficult to see that every commutative subalgebra~$\OA_2$
  of~$\Mat[2]$ has dimension at most~$2$: Otherwise $\OA_2\cap \OA_2^\ast$ has
  codimension at most~$2$ and thus contains an element that is not a multiple
  of~$1$.  Hence~$\OA_2$ contains all diagonal matrices in a particular basis.
  But only diagonal matrices commute with all diagonal matrices in a
  particular basis, so that $\dim\OA_2=2$ contrary to assumption.  Since every
  representation of a commutative operator algebra has commutative range,
  Corollary~\ref{cor:contr2} can be applied.
\end{proof}

\section[Quotient distance]{The quotient distance and metric for unital
  operator algebras}
\label{sec:CaraDist}

\begin{example}  \label{exa:twodimQalg}
  Let~$\Md$ be a complex manifold and either $\pM_1,\pM_2\in\Md$ or
  $(\pM;X)\in\Tg\Md$.  Let $\UA=\BH(\Md)$ and either $\ide=\ide(\pM_1,\pM_2)$
  or
  $$
  \ide=\ide(\pM;X)=\{f\in\UA\mid f(\pM)=0, Df(\pM;X)=0\}.
  $$
  Assume that the codimension of~$\ide$ is two.

  In the first case, the quotient algebra~$\UA/\ide$ is completely equivalent
  to~$\QA_c$ with $c=c_\UA^\ast(\pM_1,\pM_2)\neq0$ the \Caratheodory{}\Ast
  pseudodistance of $\pM_1$ and~$\pM_2$.  In the second case, the quotient
  algebra~$\UA/\ide$ is completely equivalent to~$\QA_0$.
\end{example}

\begin{proof}
  $\UA/\ide$ must be completely isometrically isomorphic to some~$\QA_c$,
  $c\in[0,1]$.  In the first case, there are two distinct characters on
  $\UA/\ide$ because $\dim \UA/\ide=2$.  This excludes the possibility $c=0$.
  Moreover, comparing~\eqref{equ:c} with the definition of the
  \Caratheodory{}\Ast pseudodistance shows that $c=c_\UA^\ast(\pM_1,\pM_2)$.

  In the second case, there is one character and a non-trivial derivation
  on~$\UA/\ide$.  Thus $\UA/\ide\not\cong\QA_c$ for any $c\in(0,1]$, forcing
  $\UA/\ide\cong\QA_0$.
\end{proof}

For a unital operator algebra~$\OA$, let $\Spec(\OA)$ be the set of all
nonzero, continuous homomorphisms $\OA\to\C$.  Endowed with the weak topology
from $\Spec(\OA)\subset \clBall (\OA')$, this becomes a compact \Hausdorff space.
If $f\in\OA$, $\omega\in\Spec(\OA)$, write $f(\omega)$ for $\omega(f)$.  If
$\omega\in\Spec(\OA)$, a linear functional $d\colon A\to \C$ with
$d(fg)=f(\omega)d(g) + d(f)g(\omega)$ for all $f,g\in\OA$ is called a
\emph{derivation at~$\omega$}.  Write $\Tg_\omega\OA$ for the set of all
derivations at~$\omega$ and $\Tg\OA=\coprod_{\omega\in\Spec(\OA)}
\Tg_\omega\OA$.  Obviously, $\Tg_\omega\OA$ is a complex vector space, called
the \emph{tangent space of~$\OA$ at~$\omega$}.  $\Tg\OA$ is called the
\emph{tangent space} of~$\OA$.

These definitions can be made for not necessarily commutative operator
algebras.  However, it is easy to check that $[f,g](\omega)=0$ and $d[f,g]=0$
for all $f,g\in\OA$, $\omega\in\Spec(\OA)$, $d\in\Tg_\omega\OA$, and
$[f,g]=fg-gf$.  Let $[\OA,\OA]$ be the ideal generated by all commutators.
Then all elements of the spectrum and all derivations at some point of the
spectrum annihilate $[\OA,\OA]$.  Hence they factor through the commutative
operator algebra $\OA/[\OA,\OA]$.  All the constructions in this section will
ignore any noncommutativity of~$\OA$ in this way.

\begin{definition}
  Let~$\OA$ be a \emph{unital} operator algebra and
  $\omega_1,\omega_2\in\Spec(\OA)$.  Define
  $$
  c_{\OA}^\ast(\omega_1,\omega_2) = \sup\{ |f(\omega_2)| \mid
	\text{$f\in\Ball(\OA)$ and $f(\omega_1)=0$}\}
  $$
  and let $c_{\OA}={\Artanh}\circ c_{\OA}^\ast$ with the convention
  $\Artanh(1)=\infty$.  The functions $c_{\OA}^\ast$ and~$c_{\OA}$ are called
  the \emph{quotient\Ast distance} and \emph{quotient distance} respectively.

  If $(\omega,d)\in\Tg\OA$, let
  $$
  \gamma_{\OA}(\omega,d) = \sup \{ |d(f)| \mid
	\text{$f\in\Ball(\OA)$ and $f(\omega)=0$}\}.
  $$
  $\gamma_{\OA}$ is called the \emph{quotient metric} for~$\OA$.
\end{definition}

Clearly, $c_{\OA}^\ast(\omega_1,\omega_2)\le1$ for all
$\omega_1,\omega_2\in\Spec(\OA)$, so that~$c_{\OA}$ is a well\PrM{}defined function
from $\Spec(\OA)^2$ to $[0,\infty]$.  We write~$c_{\OA}\bast$ if an assertion
holds both for~$c_{\OA}$ and for~$c_{\OA}^\ast$.

A \emph{distance} on a set~$X$ is a symmetric function $d\colon X\times
X\to[0,\infty]$ satisfying the triangle inequality and $d(x,y)=0$ iff $x=y$.
Thus infinite distances are allowed.  This is necessary because it can easily
happen that $c_{\OA}^\ast(\omega_1,\omega_2)=1$ and thus
$c_{\OA}(\omega_1,\omega_2)=\infty$.  However, it will be shown below that
$c_{\OA}^\ast$ and~$c_{\OA}$ are distances on $\Spec(\OA)$ in the above sense.

A first justification for the names ``quotient distance'' and ``quotient
metric'' is that they behave well with respect to taking quotients:

\begin{lemma}  \label{lem:CaraQuotient}
  Let $\OA$ be a unital operator algebra, $\ide\subset\OA$ an ideal, and
  $\omega_1,\omega_2\in\Spec(\OA/\ide)$.  View $\omega_1,\omega_2\in\Spec(\OA)$
  by putting $f(\omega_j)=[f](\omega_j)$.  Then
  $$
        c_{\OA/\ide}^\ast(\omega_1,\omega_2)
  =     c_{\OA}^\ast(\omega_1,\omega_2).
  $$
  Similarly, if $(\omega,d)\in\Tg_\omega(\OA/\ide)$, then
  $$
        \gamma_{\OA/\ide}(\omega,d)
  =     \gamma_{\OA}(\omega,d).
  $$
\end{lemma}

\begin{proof}
  Trivial.
\end{proof}

\begin{theorem}  \label{the:CaraDistance}
  Let~$\OA$ be a unital operator algebra, $\omega_1,\omega_2\in\Spec(\OA)$,
  $\omega_1\neq\omega_2$.  Define
  $\OA(\omega_1,\omega_2)=\OA/\ide(\omega_1,\omega_2)$.  This is again a unital
  operator algebra and completely equivalent to~$\QA_c$ for
  $c=c_{\OA}^\ast(\omega_1,\omega_2)$.

  There exists $f\in\Ball(\OA)$ with $f(\omega_j)=\lambda_j$, $j=1,2$,
  if and only if $\lambda_1,\lambda_2\in\D$ and
  $m(\lambda_1,\lambda_2)<c_{\OA}^\ast(\omega_1,\omega_2)$, where~$m$ denotes
  the \Moebius distance.  Thus
  $$
        c_{\OA}^\ast(\omega_1,\omega_2)
  =     \sup\bigl\{ m\bigl(f(\omega_1),f(\omega_2)\bigr) \bigm|
	 f\in\Ball(\OA) \}.
  $$
  In particular, $c_{\OA}^\ast$ and~$c_{\OA}$ are distances on $\Spec(\OA)$.
\end{theorem}

\begin{proof}
  By Theorem~\ref{the:twodimOpalg}, $\OA(\omega_1,\omega_2)\cong \QA_c$ for
  some $c\in[0,1]$.  Clearly, $c=0$ is impossible and
  $c=c_{\OA}^\ast(\omega_1,\omega_2)$ follows immediately from~\eqref{equ:c}.
  For the special case $\OA=\Pol(\clD)$, the second assertion is just the
  classical Schwarz-Pick lemma, and it is well\PrM{}known that
  $m=c_{\Pol(\clD)}^\ast$ and $p=c_{\Pol(\clD)}$ are distances on~$\clD$.
  Lemma~\ref{lem:CaraQuotient} and the definition of the quotient norm yield
  the second assertion for the algebras~$\QA_c$ since they are quotients of
  $\Pol(\clD)$.  The general case follows from this in the same way.  (It is
  not difficult to give a direct proof using von Neumann's inequality,
  paralleling the argument in~\cite{Meyer:97} for uniform algebras.)

  If $f\in\Ball(\OA)$ then $f(\omega)\in\D$ for each
  $\omega\in\Spec(\OA)$.  Hence
  $$
  f^\ast m= m\circ(f,f)\colon \Spec(\OA)\times\Spec(\OA)\to\R_+
  $$
  is a well\PrM{}defined pseudodistance on $\Spec(\OA)$ since~$m$ is a distance
  on~$\D$.  Clearly, the supremum of a family of pseudodistances is again a
  pseudodistance.  It is trivial that $c_{\OA}^\ast(\omega_1,\omega_2)=0$
  implies $\omega_1=\omega_2$.  Thus~$c_{\OA}^\ast$ is a distance.
  Replacing~$m$ by~$p$, the same argument yields that~$c_{\OA}$ is a distance.
\end{proof}

If $(\omega,d)\in\Tg_\omega\OA$ with $d\neq0$, define
$\ide(\omega,d)=\{f\in\OA\mid f(\omega)=d(f)=0\}$ and
$\OA(\omega,d)=\OA/\ide(\omega,d)$.  Clearly, $\OA(\omega,d)=\OA(\omega,\lambda
d)$ for all $\lambda\in\C^\ast$.

\begin{theorem}  \label{the:CaraMetric}
  Let $(\omega,d)\in\Tg_\omega\OA$, $d\neq0$.  Then there is a complete
  equivalence $\phi\colon \OA(\omega,d)\to\QA_0$.  There exists
  $f\in\Ball(\OA)$ with $f(\omega)=\lambda$, $d(f)=l$ iff $\lambda\in\D$,
  $l\in\C$, and $\gamma(\lambda,l)< \gamma_{\OA}(\omega,d)$.  Thus
  $\gamma_{\OA}(\omega,d)=\|d\|$, and this is a norm on~$\Tg_\omega\OA$.
\end{theorem}

\begin{proof}
  By Theorem~\ref{the:twodimOpalg}, there exists a complete equivalence
  $\phi\colon\OA(\omega,d)\to\QA_0$.  Thus~$d$ induces some derivation $d\circ
  \phi^{-1}$ on~$\QA_0$.  By Lemma~\ref{lem:CaraQuotient}, it suffices to prove
  the remaining claims for the special case $\OA=\QA_0$.  This case can further
  be translated to~$\Pol(\clD)$, where everything follows from the Schwarz-Pick
  lemma and $\gamma(\lambda,l)\le \gamma(0,l)=|l|$ for all $\lambda\in\D$,
  $l\in\C$.
\end{proof}

\begin{corollary}  \label{cor:boundaryTrivial}
  Let $f\in\OA$, $\omega\in\Spec(\OA)$, with $f(\omega)=\|f\|$.  Then $d(f)=0$
  for all $d\in\Tg_\omega\OA$ and $c_{\OA}^\ast(\omega,\omega_2)=1$ for all
  $\omega_2\in\Spec(\OA)$ with $f(\omega_2)\neq f(\omega)$.
\end{corollary}

\begin{proof}
  Apply Theorem~\ref{the:CaraDistance} and~\ref{the:CaraMetric} to
  $(\|f\|+\epsilon)^{-1} f$ for $\epsilon>0$.
\end{proof}

\begin{example}
  Let $\OA=\NBC(\Omega)$, $\Omega$ a compact \Hausdorff space.  Then
  $\Spec(\OA)$ with the weak topology is homeomorphic to~$\Omega$.  However, it
  is easy to see that $c_{\OA}^\ast(\omega_1,\omega_2)=1$ for all
  $\omega_1,\omega_2\in\Omega$ with $\omega_1\neq\omega_2$: there exists
  $f\in\OA$ with $\|f\|=1=f(\omega_1)$, $f(\omega_2)=0$.  Moreover,
  Corollary~\ref{cor:boundaryTrivial} implies that $\Tg_\omega\OA=\{0\}$ for
  all $\omega\in\Omega$.  Thus the topology on~$\Omega$ defined
  by~$c_{\OA}^\ast$ is always the discrete topology, which is different from
  the usual topology unless~$\Omega$ is finite.
\end{example}

Thus the quotient distance and metric do not say anything interesting about
commutative \CstarAlgebra{}s.  In a sense, they measure how much a commutative
operator algebra deviates from being self-adjoint.  For function algebras they
measure whether there is a relation (in form of inequalities) between the
function values at different points.

The following theorem is the analogue of the holomorphic contractiveness of the
classical \Caratheodory{}\bAst pseudodistance and the \Caratheodory-\Reiffen
pseudometric.

\begin{theorem}  \label{the:contractivity}
  Let $\rho\colon\OA_1\to\OA_2$ be a contractive, unital homomorphism.
  Then~$\rho$ induces natural maps $\rho^\ast\colon
  \Spec(\OA_2)\to\Spec(\OA_1)$ and $D\rho^\ast(\omega)\colon
  \Tg_\omega\OA_2\to\Tg_{\rho^\ast\omega}\OA_1$ for all
  $\omega\in\Spec(\OA_2)$.  This yields a map $D\rho^\ast\colon
  \Tg\OA_2\to\Tg\OA_1$ over~$\rho^\ast$.  These maps are contractions for the
  quotient\bAst distance and the quotient metric, respectively, i.e.
  \begin{align*}
	c_{\OA_1}\bast(\rho^\ast\omega_1,\rho^\ast\omega_2)
  &\le	c_{\OA_2}\bast(\omega_1,\omega_2),
\\	\gamma_{\OA_1}\bigl(\rho^\ast\omega;D\rho^\ast(\omega;d)\bigr)
  &\le	\gamma_{\OA_2}(\omega;d).
  \end{align*}
\end{theorem}

\begin{proof}
  Define $\rho^\ast(\omega)=\omega\circ\rho$ and
  $D\rho^\ast(\omega;d)=(\omega\circ\rho;d\circ\rho)$.  All the algebraic
  properties are easy to check.  To get the contractiveness of the maps, note
  that, for $\omega\in\Spec(\OA_2)$, $\rho^{-1}\bigl(\ide(\omega)\bigr)=
  \ide\bigl(\rho^\ast(\omega)\bigr)$ and that the restriction $\rho\colon
  \ide(\rho^\ast\omega)\to\ide(\omega)$ is still contractive.  Hence, for any
  linear functional~$l$ on~$\OA_2$,
  $$
        \sup\bigl\{ |l\circ\rho(f)| \bigm|
	 f\in\Ball\bigl(\ide(\rho^\ast\omega)\bigr) \bigr\}
  \le   \sup\bigl\{ |l(f)| \bigm| f\in\Ball\bigl(\ide(\omega)\bigr) \bigr\}.
  $$
  Now the result follows by specializing to derivations and the linear
  functionals $f\mapsto f(\omega_2)$.
\end{proof}

\begin{proposition}  \label{pro:CaraTopFiner}
  The topology on~$\Spec(\OA)$ generated by the distance~$c_{\OA}^\ast$ is
  finer than the weak topology.
\end{proposition}

\begin{proof}
  Assume that the net~$(\omega_j)$ in~$\Spec(\OA)$ converges towards some
  $\omega_\infty\in\Spec(\OA)$ in the distance topology.  This implies
  $m\bigl(f(\omega_j),f(\omega_\infty)\bigr)\to0$ for all $f\in\Ball(\OA)$.
  Since~$m$ is a distance on~$\D$, the net $\bigl(f(\omega_j)\bigr)$ converges
  to $f(\omega_\infty)$, even for $f\in\OA$ with~$\|f\|$ arbitrary.  Hence
  $\omega_j\to\omega_\infty$ in the weak topology.
\end{proof}

\begin{example}
  Let~$\Md$ be a complex manifold, $\OA=\BH(\Md)$.  There is a canonical map
  $\Md\to\Spec(\OA)$.  The quotient\bAst distance on~$\Spec(\OA)$ pulls back
  to a pseudodistance on~$\Md$ (which may fail to separate points if $\BH(\Md)$
  does not separate the points of~$\Md$).  This is nothing but the classical
  \Caratheodory{}\bAst pseudodistance.  Similarly, there is a canonical map
  $\Tg\Md\to \Tg\OA$ by restricting a derivation on $\CINF(\Md)$ to $\BH(\Md)$.
  Under this map, the quotient metric on $\Spec(\OA)$ pulls back to the
  classical \Caratheodory-\Reiffen metric on~$\Md$.
\end{example}

Complex analysts may be surprised by Proposition~\ref{pro:CaraTopFiner} because
the \PrMn{c}topology, which is the topology on~$\Md$ defined by the
distance~$c_\Md\bast$ is always \emph{weaker} than the usual topology on~$\Md$
as a complex manifold, and may be strictly weaker (even if $\Md\to\Spec(\OA)$
is injective) \cite{Jarnicki-Pflug:93}.

On $\BH(\Md)$, there is also the topology of locally uniform convergence, which
is weaker than the norm topology and has the nice property that the closed unit
ball of $\BH(\Md)$ is compact in this topology.  This implies that the
function~$c_\Md^\ast$ is continuous with respect to the manifold topology on
$\Md\times\Md$.  Thus the \PrMn{c}topology is weaker than the manifold
topology.  However, the mapping $\Md\to\Spec(\OA)$ may fail to be a
homeomorphism onto its image.  In fact, by Proposition~\ref{pro:CaraTopFiner}
this must happen if the \PrMn{c}topology on~$\Md$ is strictly weaker than the
manifold topology.

The cases $c_{\OA}(\omega_1,\omega_2)=\infty$ and
$c_{\OA}(\omega_1,\omega_2)<\infty$ are qualitatively different.  Call
$\omega_1,\omega_2\in\Spec(\OA)$ \emph{\PrMn{\OA}related} if
$c_{\OA}(\omega_1,\omega_2)<\infty$ and write $\omega_1\sim\omega_2$.
Since~$c_{\OA}$ satisfies the triangle inequality, this is an equivalence
relation.  Call $\Spec(\OA)$ \emph{\PrMn{\OA}connected} if all elements
of~$\Spec(\OA)$ are \PrMn{\OA}related.  Conversely, if $\omega_1\sim\omega_2$
implies $\omega_1=\omega_2$, then $\Spec(\OA)$ is called \emph{totally
\PrMn{\OA}disconnected}.

By definition, $\Spec(\OA)$ is \PrMn{\OA}connected if and only if there are no
two\PrM{}dimensional quotients that are equivalent to~$\QA_1$.  $\QA_1$~is the
only two\PrM{}dimensional unital operator algebra that can be written as an
orthogonal direct sum of two one\PrM{}dimensional operator algebras.  Call an
operator algebra \emph{decomposable} if it is an orthogonal direct sum of two
nonzero operator algebras, and \emph{indecomposable} otherwise.
Indecomposability is closely linked with \PrMn{\OA}connectedness of the
spectrum.

Let~$\OA$ be a finite dimensional commutative operator algebra.  Then
$\Spec(\OA)$ is a finite set, and both the weak and the \PrMn{c}topology are
discrete.  First, recall some general algebraic facts about finite dimensional
algebras.

Let $\omega\in\Spec(\OA)$, and let $\ide=\ide(\omega)\subset\OA$ be the
corresponding maximal ideal.  The sequence of ideals $\ide\supset \ide^2\supset
\cdots$ becomes constant, i.e.\ $\ide^k=\ide^{k+1}$ for some $k\in\N$.  Write
$\ide^\infty=\ide^k$.  These ideals are still coprime, i.e.\
$\ide(\omega_1)^\infty + \ide(\omega_2)^\infty = \OA$ for
$\omega_1,\omega_2\in\Spec(\OA)$, $\omega_1\neq\omega_2$.  Moreover,
$\bigcap_{\omega\in\Spec(\OA)} \ide(\omega)^\infty = \{0\}$.  The Chinese
remainder theorem yields a direct sum decomposition
$$
\OA \cong \bigoplus_{\omega\in\Spec(\OA)} \OA/\ide(\omega)^\infty \quad
\mbox{(algebraically)}.
$$
Given a completely isometric representation of~$\OA$ on~$\HilS$, this
corresponds to a decomposition of~$\HilS$ into generalized eigenspaces
\begin{equation}  \label{equ:HilsOm}
\begin{split}
        \HilS(\omega)
  &=    \bigl\{ \xi\in\HilS \bigm|
          \mbox{$\bigl(f-f(\omega)\bigr)^k\xi=0$ for some
          $k\in\N$, all $f\in\OA$} \bigr\}
\\&=    \bigcup_{k\in\N} \bigcap_{f\in\OA} \Ker\bigl(f-f(\omega)\bigr)^k
   =    \bigcap_{f\in\OA} \bigcup_{k\in\N} \Ker\bigl(f-f(\omega)\bigr)^k.
\end{split}
\end{equation}

The subspaces~$\HilS(\omega)$ are \PrMn{\OA}invariant,
and~$\ide(\omega)^\infty$ acts on~$\HilS(\omega)$ by zero.  Moreover, in a
suitable orthonormal basis, the action of all $f\in\OA$ on~$\HilS(\omega)$ is
(jointly) upper triangular with entries $f(\omega)$ in the diagonal.  In most
cases, the spaces~$\HilS(\omega)$ will not be orthogonal.  But always $\sum
\HilS(\omega) = \HilS$.

Now decompose $\Spec(\OA)$ into classes of related elements $C_1,\dots,C_m$,
and let $\ide(C_j)=\bigcap_{\omega\in C_j} \ide(\omega)^\infty$.  Since
$\bigcup C_k =\Spec(\OA)$, the Chinese remainder theorem yields that,
algebraically, $\OA\cong\bigoplus_{k=1}^m \OA/\ide(C_k)$.  However, in this
case the decomposition is orthogonal:

\begin{theorem} \label{the:OrthDirectSum}
  Let~$\OA$ be a finite dimensional, unital, commutative operator algebra
  and let $C_1,\dots,C_m$ be the \PrMn{\sim}equivalence classes of
  $\Spec(\OA)$.  Then the canonical isomorphism $\phi\colon \OA\to
  \bigoplus_{k=1}^m \OA/\ide(C_k)$ is completely isometric, where the right
  side has the orthogonal direct sum matrix normed structure.  Furthermore, the
  operator algebras $\OA/\ide(C_k)$ are indecomposable.

  In fact, if the commutative operator algebra~$\OA[B]$ has a
  \PrMn{\OA[B]}connected spectrum, then~$\OA[B]$ is not even isometrically
  isomorphic to a non-trivial direct sum, where non-trivial means that both
  summands are nonzero.
\end{theorem}

\begin{proof}
The quotient maps $\OA\to \OA/\ide(C_k)$ for $k=1,\dots,m$ induce a completely
contractive map $\phi\colon \OA\to \bigoplus_{k=1}^m \OA/\ide(C_k)$ by
definition of the orthogonal direct sum of operator algebras.  Moreover,
$\phi$~is an isomorphism by abstract algebra.  The point is to show
that~$\phi^{-1}$ is completely contractive.  Choose a completely isometric
representation $\rho\colon \OA\to\Bound(\HilS)$ on some Hilbert space~$\HilS$.
Let $\omega\in\Spec(\OA)$, and let $\HilS(\omega)$ be as in~\eqref{equ:HilsOm}.
Choose an orthonormal basis of~$\HilS(\omega)$ making the \PrMn{\OA}action
upper triangular.

Let~$\omega_j$, $j=1,2$, lie in different \PrMn{\sim}equivalence classes.
Since~$\OA$ is finite dimensional, its unit ball is compact.  Thus by
Theorem~\ref{the:CaraDistance} there is $f\in\OA$ with $\|f\|\le1$,
$f(\omega_1)=1$, and $f(\omega_2)=-1$.  For $j=1,2$, the restriction
of~$\rho(f)$ to~$\HilS(\omega_j)$ must be of the form
$$
\begin{pmatrix}
f(\omega_j) & \ast & \ast & \ast & \cdots \\
0 & f(\omega_j) & \ast & \ast & \cdots \\
0 & 0 & f(\omega_j) & \ast & \cdots \\
0 & 0 & 0 & f(\omega_j) & \cdots \\
0 & 0 & 0 & 0 & f(\omega_j) \\
\end{pmatrix}.
$$
Since $\|\rho(f)\|=\|f\|=1$ and $|f(\omega_j)|= 1$, this is only possible
if~$\rho(f)|_{\HilS(\omega_j)}$ is diagonal for $j=1,2$.  Pick unit vectors
$x_j\in\HilS(\omega_j)$ for $j=1,2$ and let $\lambda=\5{x_1}{x_2}$.  The goal
is to show $\lambda=0$.  Since $\rho(f)x_1=x_1$ and $\rho(f)x_2=-x_2$,
\begin{align*}
	|a_1|^2 + |a_2|^2 - 2\RE(a_1a_2\lambda)
  &=	\|a_1x_1-a_2x_2\|^2
   =	\|\rho(f)(a_1x_1+a_2x_2)\|^2
\\&\le	\|f\|^2\|a_1x_1+a_2x_2\|^2
   =	|a_1|^2 + |a_2|^2 + 2\RE(a_1a_2\lambda)
\end{align*}
for all $a_1,a_2\in\C$.  Thus $\RE(a_1a_2\lambda)\ge 0$ for all $a_1,a_2\in\C$.
This is only possible if $\lambda=0$, so that
$\HilS(\omega_1)\bot\HilS(\omega_2)$.  Thus the subspaces $\HilS_j=
\sum_{\omega\in C_j} \HilS(\omega)$ are orthogonal to each other.  Clearly,
$\bigoplus \HilS_j=\HilS$.

Define $\rho_j\colon \OA\to\Bound(\HilS_j)$ by
$\rho_j(f)=P_{\HilS_j}\rho(f)|_{\HilS_j}$.  Then $\rho_j|_{\ide(C_j)}=0$, so
that~$\rho_j$ induces a completely contractive representation $\bar\rho_j\colon
\OA/\ide(C_j)\to \Bound(\HilS_j)$.  Together, these mappings induce a
completely contractive representation
$$
\bigoplus \bar\rho_j\colon
        \bigoplus \OA/\ide(C_j)\to \Bound\left(\bigoplus \HilS_j\right)
	=\Bound(\HilS).
$$
Clearly, $\left(\bigoplus \bar\rho_j\right)\circ\phi=\rho$.  Since~$\rho$ is
completely isometric and both $\bigoplus \bar\rho_j$ and~$\phi$ are completely
contractive, $\phi$~must be completely isometric.  This proves the first claim.

For the second claim, it suffices to show that if~$\OA[B]$ is isometrically
isomorphic to the direct sum $\OA[B]_1\oplus \OA[B]_2$ for some nonzero
operator algebras $\OA[B]_1$ and~$\OA[B]_2$, then $\Spec(\OA[B])$ is not
\PrMn{\OA[B]}connected.  $\OA[B]_1$ and~$\OA[B]_2$ must be ideals of~$\OA[B]$
and $\OA[B]/\OA[B]_1\cong \OA[B]_2$, $\OA[B]/\OA[B]_2\cong \OA[B]_1$
isometrically.  Choose $\omega_j\in\Spec(\OA[B])$ such that
$\ide(\omega_j)\supset\OA[B]_j$.  Write $1\in \OA[B]$ as a sum $e_1+e_2$,
$e_j\in \OA[B]_j$, $j=1,2$.  Then
$$
1=\|1\|=\|e_1+e_2\|=\max\{\|e_1\|,\|e_2\|\}=\|e_1-e_2\|.
$$
Moreover,
\begin{align*}
	\omega_1(e_1-e_2)
  &=	-\omega_1(1)+2\omega_1(e_1)
   =	-1,
\\	\omega_2(e_1-e_2)
  &=	\omega_2(1)-2\omega_2(e_2)
   =	1.
\end{align*}
Together with $\|e_1-e_2\|=1$ this implies
$c_{\OA[B]}^\ast(\omega_1,\omega_2)=1$, so that $\Spec(\OA[B])$ is not
\PrMn{\OA[B]}connected.
\end{proof}

It is easy to find examples of indecomposable operator algebras whose spectrum
is not \PrMn{\OA}connected, for example $\NBC(\Omega)$ for a connected
space~$\Omega$.  However, this only occurs for infinite dimensional operator
algebras, where there are also topological obstructions to decomposability.

If $\Spec(\OA)$ is totally \PrMn{\OA}disconnected it cannot be concluded that
$\OA\cong\NBC\bigl(\Spec(\OA)\bigr)$ because of examples like~$\QA_0$.  Such
operator algebras can be distinguished from the self-adjoint case by their
tangent space.

Let~$\OA$ be a unital operator algebra, $\omega\in\Spec(\OA)$, and
$d\in\Tg_\omega\OA$.  Then $d|_{\ide(\omega)^2}=0$, so that $d|_{\ide(\omega)}$
determines a continuous linear functional~$\delta$ on
$\ide(\omega)/\ide(\omega)^2$.  Moreover, $d$~is uniquely determined
by~$\delta$ because $d(1)=0$.  Conversely, if $\delta\colon
\ide(\omega)/\ide(\omega)^2\to\C$ is a linear functional, then
$d(f)=\delta\bigl([f-f(\omega)]\bigr)$ defines a derivation at~$\omega$.
Hence~$\Tg_\omega\OA$ is the dual space of $\ide(\omega)/\ide(\omega)^2$
(compare this with Exercise~2.12 of~\cite{Spivak:79}).

\begin{definition}
  Let~$\OA$ be a unital operator algebra, $\omega\in\Spec(\OA)$.  The
  \emph{cotangent space of~$\OA$ at~$\omega$} is the abstract operator algebra
  $$
  \CTg_\omega\OA=\ide(\omega)/\ide(\omega)^2.
  $$
  $\OA(\omega)=\OA/\ide(\omega)^2$ is the trivial unitization
  of~$\CTg_\omega\OA$.
\end{definition}

Let $\natural\colon \Tg_\omega\OA\to (\CTg_\omega\OA)'$ be the bijection
constructed above.  One of the consequences of Theorem~\ref{the:CaraMetric} is
that, for any $d\in\Tg_\omega\OA$,
$$
        \|d\|
  =     \|d|_{\ide(\omega)}\|
  =     \gamma_{\OA}(\omega,d)
  =     \|\natural(d)\|.
$$
Thus~$\Tg_\omega\OA$ with the quotient metric is the normed dual
of~$\CTg_\omega\OA$.

For the special case of the function algebras $\Rat(K)$, $K\subset\C^n$
compact, the tangent and cotangent spaces at points $k\in K$ were already
introduced by Paulsen in~\cite{Paulsen:92}.  He also endowed~$\Tg_\omega\OA$
with a \PrMn{\LP[1]}matricially normed structure making it the standard
operator space dual of~$\CTg_k\Rat(K)$.

A commutative, unital operator algebra~$\OA$ is said to have \emph{zero tangent
space} if $\Tg_\omega\OA=\{0\}$ for all $\omega\in\OA$.  Equivalently,
$\CTg_\omega\OA=\{0\}$ for all $\omega\in\Spec\OA$, i.e.\
$\ide(\omega)=\ide(\omega)^2$ for all $\omega\in\Spec\OA$.  For example,
if~$\Omega$ is a compact \Hausdorff space, then $\NBC(\Omega)$ has zero
tangent space.

\begin{theorem}  \label{the:finiteDimSelfAdjoint}
  Let~$\OA$ be a unital, commutative operator algebra of dimension $n\in\N$.
  Then~$\OA$ is completely equivalent to $\NBC(\{1,\dots,n\})$ if and only
  if~$\OA$ has zero tangent space and $\Spec(\OA)$ is totally
  \PrMn{\OA}disconnected.
\end{theorem}

\begin{proof}
  Of course, $\NBC(\{1,\dots,n\})$ has zero tangent space and totally
  \PrMn{\OA}disconnected spectrum.  Assume conversely that~$\OA$ has zero
  tangent space and totally \PrMn{\OA}disconnected spectrum.  By
  Theorem~\ref{the:OrthDirectSum}, $\OA\cong \bigoplus
  \OA/\ide(\omega)^\infty$, where the sum runs over all $\omega\in\Spec(\OA)$.
  But $\CTg_\omega\OA$ means $\ide(\omega)^2=\ide(\omega)$ and thus
  $\ide(\omega)^\infty=\ide(\omega)$, so that $\OA/\ide(\omega)^\infty\cong\C$.
  Thus~$\OA$ is an orthogonal direct sum of several copies of~$\C$, which is
  completely equivalent to $\NBC(\{1,\dots,n\})$.
\end{proof}

\begin{remark}
  There are commutative, unital operator algebras $\OA_1$ and~$\OA_2$ that have
  zero tangent space and isometric spectra $(\Spec(\OA),c_{\OA}^\ast)$, but
  which are not isometric.  Examples are appropriate quotients of
  $\Pol(\cl{\D_d})$ and $\Shift_d$.  Thus a lot of information is lost by
  looking only at the tangent space and the spectrum with the quotient\bAst
  distance.
\end{remark}

\subsection{The quotient distance and metric for tensor products}
\label{sec:quotTensor}

Let $\OA_1$ and~$\OA_2$ be unital operator algebras.  There is no unique way to
turn their algebraic tensor product into an operator algebra.  The most natural
choices are the \emph{spatial} and the \emph{maximal} tensor
product~\cite{Paulsen-Power:92}.  If $\OA_j\subset\Bound(\HilS_j)$, then the
spatial tensor product structure comes from the natural representation
$\OA_1\otimes\OA_2\to\Bound(\HilS_1\otimes\HilS_2)$.  The maximal tensor
product structure has the maximal matrix norms for which the embedding
$\OA_1\to\OA_1\otimes\OA_2$ and $\OA_2\to\OA_1\otimes\OA_2$ are still
completely contractive.  Thus a representation of the maximal tensor product is
completely contractive iff its restrictions to $\OA_1\otimes\{\ID\}$ and
$\{\ID\}\otimes \OA_2$ are completely contractive.  Their is no such criterion
for a representation of the spatial tensor product to be completely
contractive.  However, the maximal tensor product, like most universal objects,
does not come with an interesting completely isometric representation.

Both tensor product structures are functorial for completely contractive maps
in the sense that if $\rho_j\colon \OA_j\to\OA[B]_j$, $j=1,2$, are completely
contractive maps, then $\rho_1\otimes\rho_2$ is also completely contractive.
This is trivial for the maximal tensor product, but not for the spatial tensor
product (the corresponding statement for contractive representations is false
because otherwise every contractive map would be completely contractive).
See~\cite{Paulsen:86} for a proof.

If $\omega_j\in\Spec(\OA_j)$, $j=1,2$, then $(\omega_1,\omega_2)\colon f\otimes
g\mapsto f(\omega_1)\cdot g(\omega_2)$ defines a character of
$\OA_1\otimes\OA_2$, which is continuous both for the spatial and the maximal
tensor product structure.  Moreover, every character of $\OA_1\otimes\OA_2$ is
of this form.  If $d_j\in\Tg_{\omega_j}\OA_j$, $j=1,2$, then
$$
(d_1,d_2)\colon f\otimes g\mapsto f(\omega_1) d_2(g) + d_1(f) g(\omega_2)
$$
is a continuous derivation at $(\omega_1,\omega_2)$.  Moreover, every
derivation at $(\omega_1,\omega_2)$ is of this form.

The spatial tensor product of uniform algebras $\UA_j\subset\NBC(\Omega_j)$,
$j=1,2$, is a unital function algebra on $\Omega_1\times\Omega_2$.  In this
special case, there is the following formula for the quotient distance and
metric of the spatial tensor product:

\begin{theorem}[\cite{Meyer:97}]  \label{the:ppua}
  Let $\UA_1$ and~$\UA_2$ be uniform algebras or, more generally, unital
  function algebras.  If $\omega_{j,k}\in\Spec(\UA_k)$, $j,k\in\{1,2\}$, then
  \begin{equation}  \label{equ:ppcd}
        c_{\UA_1\otimes\UA_2}\bast\bigl(
         (\omega_{1,1},\omega_{1,2}),(\omega_{2,1},\omega_{2,2})\bigr)
  =     \max\{  c_{\UA_1}\bast(\omega_{1,1},\omega_{2,1}),
                c_{\UA_2}\bast(\omega_{1,2},\omega_{2,2}) \}.
  \end{equation}
  If $d_j\in\Tg_{\omega_j}\UA_j$, then
  \begin{equation}  \label{equ:ppcm}
	\gamma_{\UA_1\otimes\UA_2}(d_1,d_2)
  =	\max\{ \gamma_{\UA_1}(d_1),\gamma_{\UA_2}(d_2) \}.
  \end{equation}
\end{theorem}

The inclusion $\OA_1\to\OA_1\otimes\OA_2$ is completely contractive for any
reasonable tensor product structure.  Hence Theorem~\ref{the:contractivity}
shows easily that the inequality~``$\ge$'' holds in \eqref{equ:ppcd}
and~\eqref{equ:ppcm} also for more general operator algebras and any reasonable
tensor product structure.  Moreover, the estimate~``$\le$'' is trivial for the
maximal tensor product structure.  Hence the analogues of \eqref{equ:ppcd}
and~\eqref{equ:ppcm} hold for the maximal tensor product of unital operator
algebras.  However, this is not a generalization of Theorem~\ref{the:ppua}
because the spatial tensor product of uniform algebras considered there is
usually different from the maximal tensor product.  E.g., the maximal tensor
product of $\Pol(\clD^2)$ and $\Pol(\clD)$ is the universal operator algebra
for three commuting contractions and thus different from $\Pol(\clD^3)$.

However, the situation is more complicated for the spatial tensor product.

\begin{theorem}  \label{the:ppWorstCase}
  For $c,d\in(0,1]$, let
  $$
	\phi(c,d)
  =	c_{\QA_c\otimes\QA_d}^\ast\bigl((0,0),(c,d)\bigr).
  $$
  If $\OA_1$ and~$\OA_2$ are any unital operator algebras and
  $\omega_{j,k}\in\Spec(\OA_k)$, then
  \begin{equation}
  \begin{split}  \label{equ:ppWorstCase}
  &\phantom{\le}
	\max\{	c_{\OA_1}^\ast(\omega_{1,1},\omega_{2,1}),
		c_{\OA_2}^\ast(\omega_{1,2},\omega_{2,2}) \} \\
  &\le	c_{\OA_1\otimes\OA_2}^\ast\bigl(
	 (\omega_{1,1},\omega_{1,2}),(\omega_{2,1},\omega_{2,2})\bigr) \\
  &\le	\phi\bigl(c_{\OA_1}^\ast(\omega_{1,1},\omega_{2,1}),
		c_{\OA_2}^\ast(\omega_{1,2},\omega_{2,2})\bigr).
  \end{split}
  \end{equation}
\end{theorem}

\begin{proof}
  The lower bound in~\eqref{equ:ppWorstCase} follows from the above discussion.
  It remains to prove the upper bound.  Let
  $c=c_{\OA_1}^\ast(\omega_{1,1},\omega_{2,1})$ and
  $d=c_{\OA_2}^\ast(\omega_{1,2},\omega_{2,2})$.  There are completely
  isometric isomorphisms $\OA_1/\ide(\omega_{1,1},\omega_{2,1})\to\QA_c$ and
  $\OA_2/\ide(\omega_{1,2},\omega_{2,2})\to\QA_d$.  The induced map
  $\OA_1\otimes\OA_2\to \QA_c\otimes\QA_d$ is completely contractive.  Hence
  the induced map $\Spec(\QA_c\otimes\QA_d)\to\Spec(\OA_1\otimes\OA_2)$ is a
  contraction for the quotient distance.  This is the assertion.
\end{proof}

Hence, with respect to estimating the quotient distance for spatial tensor
products, the operator algebras~$\QA_c$ are the worst case.  Unfortunately, the
computation of~$\phi(c,d)$ is more complicated than one might expect.

\begin{problem}
  Compute the function~$\phi$ of Theorem~\ref{the:ppWorstCase}.
\end{problem}

Compressing the standard representation of
$\QA_c\otimes\QA_d\subset\Bound(\C^4)$ to the subspace spanned by the joint
eigenvectors corresponding to $(0,0)$ and $(c,d)$ yields a completely
contractive representation of $\QA_c\otimes\QA_d$ by \PrM{2\times2}matrices.
This yields an upper bound for~$\phi(c,d)$, which can be computed from the
angle between the two joint eigenvectors.  The result is
$$
	\phi(c,d)\le \sqrt{c^2+d^2-c^2d^2}.
$$
This is better than the estimate $\phi(c,d)\le (c+d)/(1-cd)$ which follows
because the quotient distance is a distance and from the addition theorem
for~$\tanh$.

The computations are much cleaner for the quotient metric:

\begin{theorem}  \label{the:QzeroQzero}
  Let $\OA=\QA_0\otimes\QA_0\subset\Mat[4]$ be the spatial tensor product.
  Then~$\OA$ has a unique \PrMn{1}dimensional ideal~$\ide$ spanned by
  $T_0\otimes T_0$.  The quotient $\OA/\ide$ is the trivial unitization of a
  two\PrM{}dimensional operator space~$\OA[B]$ with the matrix normed structure
  \begin{equation}  \label{equ:QzeroQzero}
	\| [T_0\otimes 1]\otimes A + [1\otimes T_0]\otimes B \|_{(n)}
  =	\max \{ \|AA^\ast+BB^\ast\|^{1/2}, \|A^\ast A+B^\ast B\|^{1/2} \}
  \end{equation}
  for $A,B\in\Mat$.
\end{theorem}

\begin{proof}
  The algebraic assertions are all easy.  If $A,B,C\in\Mat$, then
  $$
	M
  =	T_0\otimes 1\otimes A + 1\otimes T_0\otimes B + T_0\otimes T_0\otimes C
  \in	\QA_0\otimes\QA_0\otimes\Mat
  $$
  is represented by the block matrix
  $$
  \begin{pmatrix}  0 & A & B & C \\ 0 & 0 & 0 & B \\ 0 & 0 & 0 & A
  \end{pmatrix}.
  $$
  There exist unitary matrices $U,V\in\Mat[2n]$ with
  $$
  U \begin{pmatrix} B \\ A \end{pmatrix} = 
  \begin{pmatrix} (A^\ast A + B^\ast B)^{1/2} \\ 0 \end{pmatrix},\qquad
  \begin{pmatrix} A & B \end{pmatrix} V = 
  \begin{pmatrix} 0 & (A A^\ast + B B^\ast)^{1/2} \end{pmatrix}.
  $$
  Hence the matrix~$M$ is unitarily equivalent to the matrix
  $$
  \begin{pmatrix}
	0 & 0 & (AA^\ast + BB^\ast)^{1/2} & C \\
	0 & 0 & 0 & (A^\ast A + B^\ast B)^{1/2} \\
	0 & 0 & 0 & 0 \\
	0 & 0 & 0 & 0
  \end{pmatrix}.
  $$
  Clearly, the norm of this is minimal for $C=0$ and has the value that occurs
  in~\eqref{equ:QzeroQzero}.
\end{proof}

\begin{corollary}
  For $j=1,2$, let $\OA_j$ be a unital operator algebra, $\omega_j\in\OA_j$,
  and $d_j\in\Tg_{\omega_j}\OA_j$.  Then, if $\OA_1\otimes\OA_2$ is the spatial
  tensor product,
  $$
	\max\{ \gamma_{\OA_1}(d_1), \gamma_{\OA_2}(d_2) \}
  \le	\gamma_{\OA_1\otimes\OA_2}(d_1,d_2)
  \le	\bigl(\gamma_{\OA_1}(d_1)^2 + \gamma_{\OA_2}(d_2)^2\bigr)^{1/2}.
  $$
\end{corollary}

\begin{proof}
  If $\OA_1=\OA_2=\QA_0$, this follows at once from Theorem~\ref{the:QzeroQzero}
  because the norm on $\QA_0\otimes\QA_0/(T_0\otimes T_0)$ is of the form
$$
	\| \lambda\cdot T_0\otimes 1 + \mu\cdot 1\otimes T_0\|
  =	\bigl(|\lambda|^2 +|\mu|^2\bigr)^{1/2}.
$$
  For general $\OA_1,\OA_2$, the estimate follows as in the proof of
  Theorem~\ref{the:ppWorstCase}.
\end{proof}

\section{Transposition}
\label{sec:Transposition}

Let~$\HilS$ be a Hilbert space and let $U\colon \HilS\to\HilS$ be an
\emph{anti}-unitary operator.  Then $T\mapsto UT^\ast U^{-1}$ defines a
\emph{linear} isometry $\blank^t\colon \Bound(\HilS)\to\Bound(\HilS)$ called a
\emph{transposition}.  Let $\NS\subset\Bound(\HilS)$ be an operator space.
Then $\NS^t= \blank^t(\NS)$ is another operator space and
$\blank^t\colon\NS\to\NS^t$ is an isometric representation of~$\NS$, which
usually fails to be completely contractive.

View~$\NS^t$ as an abstract operator algebra.  A priori, $\NS^t$ depends on
the choice of a completely isometric representation of~$\NS$ and on the
anti-unitary~$U$.  But the matrix normed structure of~$\NS^t$ turns out to be
independent of these choices and can be defined intrinsically as follows:

\begin{definition}
  Let~$\OA$ be an abstract operator algebra.  Define the \emph{transposed
  algebra}~$\OA^t$ as follows: Algebraically, $\OA^t$~is the opposite algebra,
  i.e.\ has the same vector space structure, but the order of multiplication is
  reversed to $\hat{x}\circ \hat{y}=\widehat{yx}$.  Here~$\hat{\ }$ is used to
  signify that $x,y$ are viewed as elements of~$\OA^t$.  However, the
  \PrMn{\Mat}bimodule structure of $\OA_{(n)}$ is as usual: $S\cdot
  \hat{x}\cdot T= \widehat{SxT}$ for $x\in \OA_{(n)}$, $S,T\in\Mat$.

  If $x=\sum_{j=1}^N a_j\otimes T_j$, $a_j\in \OA$, $T_j\in\Mat$, define
  $x^t=\sum_{j=1}^N a_j\otimes T_j^t$, where~$T_j^t$ is the transpose of~$T_j$.
  Writing $\blank^t\colon \Mat\to\Mat$ for the transpose operation, this is
  just $\ID[\OA]\otimes \blank^t\colon \OA_{(n)}\to \OA_{(n)}$ and
  thus well\PrM{}defined.  The norm~$\|\blank\|_{(n)}^t$ on $\OA^t_{(n)}$ is
  now defined by $\|\hat{x}\|_{(n)}^t=\|x^t\|_{(n)}$.
\end{definition}

It is clear how to define~$\NS^t$ for an abstract operator space.

\begin{theorem}  \label{the:transposition}
  Let~$\NS$ be an abstract operator space.  Choose a completely isometric
  representation $\rho\colon \NS\to\Bound(\HilS)$ and an anti-unitary
  $U\in\Bound(\HilS)$.  Then $\hat{T}\mapsto U\rho(T)^\ast U^{-1}$, viewed as a
  representation of~$\NS^t$, is a completely isometric linear representation.
  Moreover, $(\NS^t)^t$ is completely equivalent to~$\NS$.  The ``identity
  map'' $\iota\colon \NS\ni \ns\to \hat{\ns}\in \NS^t$ is isometric.  If it is
  \PrMn{n}contractive, it is necessarily \PrMn{n}isometric.
  
  If~$\NS$ is a (unital) operator algebra to start with, then this
  representation is also multiplicative (and unital).  Moreover, $(\NS^t)^t$
  is completely equivalent to~$\NS$, and~$\iota$ is an anti-isomorphism.
\end{theorem}

\begin{proof}
  Let $\ns_j\in\NS$, $T_j\in\Mat$ for $j=1,\dots,N$.  Let $V\colon
  \C^n\to\C^n$ be the standard anti-unitary operator given by
  $(x_1,\dots,x_n)\mapsto (\conj{x_1},\dots,\conj{x_n})$.  Then the usual
  transposition operation for matrices is given by $T^t= VT^\ast V^{-1}$.
  Moreover, $U\otimes V^{-1}$ is an anti-unitary operator on $\HilS\otimes\C^n$.
  \begin{align*}
     \left\| \sum_{j=1}^N U\rho(\ns_j)^\ast U^{-1}\otimes T_j
     \right\|
  &= \left\|
	(U\otimes V^{-1})
	\sum_{j=1}^N \rho(\ns_j)^\ast \otimes (T_j^\ast)^t
	(U\otimes V^{-1})^{-1}
     \right\|
\\&= \left\| \sum_{j=1}^N \rho(\ns_j)\otimes T_j^t \right\|
   = \left\| \sum_{j=1}^N \ns_j\otimes T_j \right\|^t_{(n)}.
  \end{align*}
  Thus the representation is completely isometric.
  
  The algebraic assertions are all trivial, it only remains to show that
  if~$\iota$ is \PrMn{n}contractive, it is also \PrMn{n}isometric.  If this
  were false there would exist $X\in\NS_{(n)}$ with
  $\|X\|_{(n)}>\|\hat{X}\|_{(n)}^t$ and~$\iota_{(n)}$ would be contractive.
  But then
  $$
  \|X^t\|_{(n)}=\|\hat{X}\|_{(n)}^t<\|X\|_{(n)}=
  \|\widehat{X^t}\|_{(n)}^t,
  $$
  contrary to the assumption that~$\iota_{(n)}$  is contractive.
\end{proof}

Nevertheless, transposition is functorial for completely contractive maps as
follows:

\begin{lemma}  \label{lem:transposeFunctorial}
  Let $\rho\colon \NS_1\to\NS_2$ be a linear map between abstract operator
  spaces.  Then~$\rho$ gives rise to a transposed linear map $\rho^t\colon
  \NS_1^t\to\NS_2^t$ mapping~$\hat{\ns}$ to~$\widehat{\rho(\ns)}$.  This map
  satisfies $\|\rho\|_{(n)}=\|\rho^t\|_{(n)}$ for all $n\in\N$ and hence
  $\cbn{\rho}=\cbn{\rho^t}$.  Moreover, if~$\rho_{(n)}$ is a quotient map, then
  so is~$\rho^t_{(n)}$ and if~$\rho_{(n)}$ is isometric, then so
  is~$\rho^t_{(n)}$.  The same holds for (unital) homomorphisms between
  (unital) abstract operator algebras.
\end{lemma}

\begin{proof}
  This is immediate from the abstract definition of the transposed operator
  algebra.
\end{proof}

\begin{corollary}  \label{cor:Qtranspose}
  Let~$\QA$ be a Q\PrM{}algebra.  Then $\QA\to \QA^t$ is completely isometric.
\end{corollary}

\begin{proof}
  It is easy to prove this directly, but here we give an abstract nonsense
  proof using Lemma~\ref{lem:transposeFunctorial}.  Let
  $\FA\subset\NBC(\Omega)$ be a function algebra and $\ide\subset\FA$ a closed
  ideal such that $\QA=\FA/\ide$.  Then $\QA^t\cong \FA^t/\ide^t$ naturally.
  The transposition map $\QA\to\QA^t$ lifts to the transposition map
  $\FA\to\FA^t$.  By definition of the quotient structure, it suffices to show
  that this lifted map is completely isometric.  By
  Theorem~\ref{the:transposition}, this map is isometric and it suffices to
  show that its inverse $\FA^t\to\FA\subset\NBC(\Omega)$ is completely
  contractive.  This is clear because any contractive map into $\NBC(\Omega)$
  is completely contractive.
\end{proof}

\section{Cotangent spaces and counterexamples}
\label{sec:Cotangent}

It is easy to classify those unital operator algebras that have a completely
isometric representation by \PrM{3\times3}matrices and are of the
form~$\Unze{\NS}$ for some two\PrM{}dimensional operator space~$\NS$.

Fix a basis $X,Y$ for~$\NS$ and let $N= \Ran X + \Ran Y$.  Then $X$ and~$Y$
must vanish on~$N$ because the multiplication on~$\NS$ is trivial.  $X,Y\neq0$,
implies $N\neq0$, $N\neq\C^3$.  Choose an orthonormal basis of~$N$ and extend
it to an orthonormal basis of~$\C^3$.  First consider the case $\dim N=2$.
Then all elements of~$\NS$ must be of the form
$$
\begin{pmatrix} 0 & 0 & \ast \\ 0 & 0 & \ast \\ 0 & 0 & 0
\end{pmatrix}.
$$
Since $\dim \NS=2$, all such matrices occur, yielding the first example
$$
\Shift_2(0) = \left\{
\begin{pmatrix} c & 0 & a \\ 0 & c & b \\ 0 & 0 & c
\end{pmatrix},\quad a,b,c\in\C \right\}.
$$
$\Shift_2$ is defined in Section~\ref{sec:modelD}.  It will turn out that the
operator algebra above is indeed completely equivalent to $\Shift_2/\ide(0)^2$.
For the time being, this name should be viewed just as a symbol.

If $\dim N=1$, all elements of~$\NS$ must be of the form
$$
\begin{pmatrix} 0 & \ast & \ast \\ 0 & 0 & 0 \\ 0 & 0 & 0
\end{pmatrix}.
$$
Again, this determines~$\NS$.  After permuting the basis of~$\C^3$, this is
just the operator algebra $\Shift_2(0)^t$ of all transposes of elements of
$\Shift_2(0)$.

We have just proved the following theorem, which will be essential in
Section~\ref{sec:Complexity}.

\begin{theorem}  \label{the:twodimMCTU}
  Let $\OA\subset\Mat[3]$ be a trivial unitization of some operator space
  $\NS\subset\OA$.  Then~$\OA$ is unitarily equivalent to $\Shift_2(0)$ or
  $\Shift_2(0)^t$.
\end{theorem}

\begin{theorem}  \label{the:transpositionBad}
  The maximal ideals of both $\Shift_2(0)$ and $\Shift_2(0)^t$ are isometric
  to~$\ell_2^2$.  The transposition map $\iota\colon \Shift_2(0)\to
  \Shift_2(0)^t$ is an isometric isomorphism that is not \PrMn{2}contractive,
  and whose inverse $\iota^{-1}\colon \Shift_2(0)^t\to \Shift_2(0)$ is not
  \PrMn{2}contractive either.
\end{theorem}

\begin{proof}
  Transposition always is an isometric anti-homomorphism by
  Theorem~\ref{the:transposition}.  Thus for commutative operator algebras it
  is a homomorphism.  The map~$\iota_{(n)}$ is given by
  $$
  \iota_{(n)}\colon
  \begin{pmatrix} 0 & 0 & A \\ 0 & 0 & B \\ 0 & 0 & 0
  \end{pmatrix} \mapsto
  \begin{pmatrix} 0 & A & B \\ 0 & 0 & 0 \\ 0 & 0 & 0
  \end{pmatrix}.
  $$
  The norm of the matrix on the left is $\|A^\ast A+B^\ast B\|$, whereas the
  matrix on the right has norm $\|AA^\ast+BB^\ast\|$.  Especially, if
  $A,B\in\Mat[1]\cong\C$, then both norms are $|A|^2+|B|^2$, so that the
  maximal ideal of $\Shift_2(0)$ is isometric to~$\ell_2^2$.  For $A=e_{11}$,
  $B=e_{12}$, the left and right side are $1$ and~$2$, respectively, so
  that~$\iota_{(2)}$ is not contractive.  If~$\iota_{(2)}^{-1}$ were
  contractive then so would be~$\iota_{(2)}$ by
  Theorem~\ref{the:transposition}.
\end{proof}

Thus for three\PrM{}dimensional operator algebras, it is no longer true that
every contractive (or even isometric) homomorphism is completely contractive.
The proof shows that this also fails for two\PrM{}dimensional operator algebras
without unit.

\begin{corollary}  \label{cor:notQalgebra}
  Neither $\Shift_2(0)$ nor $\Shift_2(0)^t$ is completely equivalent to a
  Q\PrM{}algebra.
\end{corollary}

\begin{proof}
  The transpose map $\Shift_2(0)\to \Shift_2(0)^t$ is not completely isometric.
  But it is completely isometric for Q\PrM{}algebras by
  Corollary~\ref{cor:Qtranspose}.
\end{proof}

Thus no three\PrM{}dimensional Q\PrM{}algebra with a unique maximal ideal can be
represented completely isometrically by \PrM{3\times3}matrices.  The
three\PrM{}dimensional Q\PrM{}algebra in Example~\ref{exa:ellOne} has no finite
dimensional isometric representation at all.

For any normed algebra~$\OA$, define the corresponding \emph{maximal operator
algebra structure} $\MAX(\OA)$ as follows: Let~$(\rho_j)$ be the class of all
contractive representations of~$\OA$ on Hilbert spaces.  These induce
representations~$(\rho_j)_{(n)}$ of~$\OA_{(n)}$.  For
$X\in\OA_{(n)}$, let $\|X\|_{(n)}= \sup \|(\rho_j)_{(n)}(X)\|$.  It is
easy to check that this yields an abstract operator algebra structure on~$\OA$.
The norm $\|\blank\|_{(1)}$ coincides with the given norm on~$\OA$ iff~$\OA$ is
isometric to an operator algebra.  Evidently, every contractive representation
of~$\OA$ with this maximal operator algebra structure is completely
contractive.

The above construction works equally well for normed vector spaces.  For these
there also is a \emph{minimal operator space structure} $\MIN(\NS)$ given by
the embedding $\NS\to\NBC\bigl(\Ball(\NS')\bigr)$, where~$\NS'$ is the dual
space.  Since every contractive map into~$\NBC(\Omega)$ is completely
contractive \cite{Paulsen:86}, every contractive map into $\MIN(\NS)$ is
completely contractive.

For any normed space~$\NS$, define the map $\iota=\ID[\NS]\colon
\MIN(\NS)\to\MAX(\NS)$ and let $\alpha(\NS)=\cbn{\iota}$.  One of the main
results of~\cite{Paulsen:92} is that
$$
        \sqrt{\dim\NS}/2 \le \alpha(\NS) \le \dim\NS.
$$
Especially, $\alpha(\NS)=\infty$ for $\dim\NS=\infty$, and $\alpha(\NS)>1$ for
$\dim\NS\ge 5$.  Of course, $\alpha(\NS)=1$ if and only if $\MIN(\NS)$ is
completely equivalent to $\MAX(\NS)$.  The only known cases with
$\alpha(\NS)=1$ are $\ell_1^\infty\cong\C$, $\ell_2^\infty$, and~$\ell_2^1$,
i.e.\ this is a very rare behavior.

\begin{example}  \label{exa:MultBd}
  Paulsen shows in~\cite{Paulsen:92} that $\alpha(\ell_n^2)\le n/\sqrt{2}$.
  Especially, $\alpha(\ell_2^2)\le\sqrt{2}$.  This bound is realized by the
  algebras $\Shift_2(0)$ and $\Shift_2(0)^t$.  Both are isometric to~$\ell_2^2$
  and thus yield isometric representations $\rho\colon
  \MIN(\ell_2^2)\to\Mat[3]$.   An easy calculation shows that
  $\|\rho\|_{(2)}\ge\sqrt{2}$.  Hence
  $$
	\|\rho_{(2)}\|=\cbn{\rho}=\alpha(\ell_2^2)=\sqrt{2}.
  $$
\end{example}

Let~$\Md$ be a complex manifold and $\pM\in\Md$.  Let
$$
\ide_1=\{f\in\BH(\Md)\mid f(\pM)=0\},\quad
\ide_2=\{f\in\BH(\Md)\mid \text{$f(\pM)=0$ and $Df(\pM)=0$}\},
$$
and let $\Md(\pM)=\BH(\Md)/\ide_2$, $\CTg_\pM\Md=\ide_1/\ide_2$.
Taking~$\ide_2$ instead of $\ide(\pM)^2=\cl{\ide_1\cdot\ide_1}$ avoids the
problem whether there exist derivations on~$\BH(\Md)$ at~$\pM$ that do not come
from a geometric tangent vector of~$\Md$ at~$\pM$.  If~$\Md$ is sufficiently
nice (``\PrMn{\gamma}hyperbolic at~$\pM$'') then $\CTg_\pM\Md$ as a vector
space is the usual cotangent space.  If $\Md\subset\C^n$ is a bounded domain,
there are also the Q\PrM{}algebras $\HO(\cl\Md)/\ide(\pM)^2$,
$\Rat(\cl\Md)/\ide(\pM)^2$, and $\Pol(\cl\Md)/\ide(\pM)^2$.  They are always
algebraically isomorphic to $\Md(\pM)=\Unze{(\CTg_\pM\Md)}$ but they may have
different matrix normed structures.

It is shown in~\cite{Paulsen:92} that cotangent spaces of absolutely convex
domains\footnote{A domain is absolutely convex iff it is the unit ball of a
normed space.} at the origin are completely equivalent to a certain
$\MIN(\NS)$.  For a proof, it is convenient to use the following formula for
the \Caratheodory-\Reiffen metric on balanced domains\footnote{A domain
$\Bd\subset\C^n$ is balanced iff $\lambda x\in\Bd$ whenever $\lambda\in\clD$
and $x\in\Bd$.}.

\begin{lemma}[\cite{Jarnicki-Pflug:93}]  \label{lem:balCRmetric}
  Let $\Bd\subset\C^n$ be a balanced domain.  Under the usual identification
  $\Tg_0\Bd\cong\C^n$, the open unit ball for the \Caratheodory-\Reiffen
  metric on~$\Tg_0\Bd$ corresponds to the convex hull $\co(\Bd)$ of~$\Bd$.
\end{lemma}

\begin{theorem}  \label{the:ACQClassify}
  Let $\Bd\subset\C^n$ be a bounded, balanced domain, let~$\NS$ be the normed
  space with $\co(\Bd)$ as its open unit ball, and let~$\NS'$ be the
  corresponding dual space.  Then
  $$
        \CTg_0\Bd
  \cong \CTg_0\HO(\cl{\Bd})
  \cong \CTg_0\Rat(\cl{\Bd})
  \cong \CTg_0\Pol(\cl{\Bd})
  \cong \MIN(\NS').
  $$
\end{theorem}

\begin{proof}
  The claim is that the unit ball of $\CTg_0\blank(\Bd)_{(n)}$ coincides with
  the unit ball of $\MIN(\NS')_{(n)}$ for $\blank\in\BH,\HO,\Rat,\Pol$.  Let
  $f\in\Ball\bigl(\MIN(\NS')_{(n)}\bigr)$.  According to the definition, the
  components of~$f$ are to be viewed as functions on
  $\Ball(\NS'')\cong\Ball(\NS)$.  Thus $\|f\|<1$ means that~$f$ maps $\co(\Bd)$
  to $\Ball(\Mat)$.  Especially, $f(\Bd)\subset\Ball(\Mat)$ and of course
  $f(0)=0$.  Thus~$f$, viewed as a (polynomial, rational, or holomorphic) map
  from~$\Bd$ to $\Ball(\Mat)$ defines an element of $\CTg_0\Bd_{(n)}$, etc., of
  norm at most~$1$.  This yields a canonical completely contractive map from
  $\MIN(\NS')$ to $\CTg_0\dots$.

  Conversely, if $f\in\HO\bigl(\Bd,\Ball(\Mat)\bigr)$, $f(0)=0$, then~$Df|_\Bd$
  is a linear map with the same derivative and function value at the origin.
  Since~$Df|_0$ is a contraction for the respective \Caratheodory-\Reiffen
  metrics, Lemma~\ref{lem:balCRmetric} shows that~$Df|_0$ maps $\co(\Bd)$ into
  $\Ball(\Mat)$.  Hence the image of the closed unit ball of $\MIN(\NS')_{(n)}$
  contains the open unit ball of $\CTg_0\Bd$.  This is true a fortiori for
  $\CTg_0\HO(\cl\Bd)$, $\CTg_0\Rat(\cl\Bd)$, and $\CTg_0\Pol(\cl\Bd)$.
\end{proof}

It is important that there is no version of Theorem~\ref{the:ACQClassify} for
arbitrary points of a domain.  The \Caratheodory-\Reiffen metric usually can
only be computed in highly symmetric situations as above, so that there is no
general method to compute the normed structure of $\CTg_\pM\Md$.  Furthermore,
$\CTg_\pM\Md$ can fail to be of the form $\MIN(\NS)$.

\begin{corollary}  \label{cor:contrIntoComplBal}
  Let $\Bd\subset\C^n$ be a bounded, balanced domain, $\OA$ any unital operator
  algebra, and $\QA=\Bd(0)$.  Then every unital, contractive homomorphism
  $\OA\to\QA$ is completely contractive.
\end{corollary}

\begin{proof}
  By Theorem~\ref{the:unitizationWellDefined}, it suffices to check that the
  restriction of a unital homomorphism to the preimage of $\CTg_0\Bd$ is
  completely contractive.  This follows from $\CTg_0\Bd\cong\MIN(\NS')$.
\end{proof}

\begin{theorem}[Paulsen \cite{Paulsen:92}]  \label{the:contrComplBal}
  Let~$\Bd$ be a balanced domain and let~$\NS$ be the normed space with unit
  ball $\co(\Bd)$.  If all contractive representations of $\Rat(\cl{\Bd})$ are
  completely contractive, then $\alpha(\NS)=1$.
\end{theorem}

\begin{proof}
  The homomorphism $\iota\colon \MIN(\NS')\to\MAX(\NS')$ extends to a
  contractive homomorphism $\Unze{\iota}\colon \Unze{\MIN(\NS')}\to
  \Unze{\MAX(\NS')}$.  Since~$\CTg_0\Bd$ is completely isometric to
  $\MIN(\NS')$ by Theorem~\ref{the:ACQClassify}, $\Unze{\MIN(\NS')}$ is
  completely equivalent to~$\Bd(0)$.  If every contractive representation of
  $\Rat(\cl{\Bd})$ is completely contractive, then every contractive
  representation of~$\Bd(0)$ is completely contractive because~$\Bd(0)$ is a
  quotient of $\Rat(\cl{\Bd})$.  Especially, $\Unze{\iota}$ is completely
  contractive, forcing~$\iota$ to be completely contractive, i.e.\
  $\alpha(\NS')=1$.  It is shown in~\cite{Paulsen:92} that
  $\alpha(\NS)=\alpha(\NS')$.
\end{proof}

Especially, for all balanced domains in~$\C^n$ with $n\ge5$, there is a
contractive representation of $\Rat(\cl{\Bd})$ that is not completely
contractive.

Recall the definition $\D_2=\Ball(\ell_2^2)$.

\begin{example}  \label{exa:IsoNotCI}
  There exists a strongly convex, bounded domain $\Md\subset\C^2$ with smooth
  boundary and $\pM\in\Md$ with $\CTg_\pM\Md$ isometric but not
  \PrMn{2}isometric to $\CTg_0\D_2$.  This yields an isometric, completely
  contractive isomorphism $\rho\colon \Md(\pM)\to \D_2(0)$ between
  Q\PrMn{}algebras that is not \PrMn{2}isometric.
\end{example}

\begin{proof}
  In~\cite{Lempert:88}, \Lempert proves the existence of a strongly convex,
  bounded domain $\Md\subset\C^2$ with smooth boundary that is not
  biholomorphic to~$\D_2$, but such that for some $\pM\in\Md$ the unit ball for
  the \Caratheodory-\Reiffen metric on~$\Tg_\pM\Md$ is~$\D_2$.  Dualizing an
  isometry $T'\colon \ell_2^2\to\Tg_\pM\Md$ yields an isometry $T\colon
  \CTg_\pM\Md\to\ell_2^2\cong\CTg_0\D_2$.  Let $\rho=\Unze{T}$.  In order to
  prove the various properties of~$\rho$, it suffices to check them for~$T$ and
  to invoke Theorem~\ref{the:unitizationWellDefined}.

  By construction, $T$ is isometric.  Since $\CTg_0\D_2\cong\MIN(\ell_2^2)$ by
  Theorem~\ref{the:ACQClassify}, $T$~is completely contractive.  There exist
  contractive linear maps $\iota\colon \ell_2^2\to\Mat[2]$ and $\pi\colon
  \Mat[2]\to\ell_2^2$ with $\pi\circ\iota =\ID[\ell_2^2]$, i.e.~$\ell_2^2$ is a
  linear retract of~$\Mat[2]$.  For example, put
  $$
  \iota(x,y) = \begin{pmatrix} x & 0 \\ y & 0 \end{pmatrix},
  \quad
  \pi \begin{pmatrix} x & a \\ y & b \end{pmatrix} = (x,y).
  $$
  Consider $[\iota]\in (\CTg_0\D_2)_{(2)}\cong \MIN(\ell_2^2)_{(2)}$, then
  $\|[\iota]\|_{(2)}\le1$.  View~$\pi$ as a holomorphic map
  $\Ball(\Mat[2])\to\D_2$.

  If~$T^{-1}$ were \PrMn{2}contractive, then
  $\|T^{-1}_{(2)}[\iota]\|_{(2)}\le1$, i.e.\ there would be
  $\hat\iota\in\BH(\Md)_{(2)}$ with $T_{(2)}[\hat\iota]=[\iota]$ and
  $\|\hat\iota\|_{(2)}\le1$ (using that the unit ball of $\BH(\Md)$ is compact
  in the topology of locally uniform convergence).  Thus
  $\hat\iota(\Md)\subset \clBall (\Mat[2])$.  Since $\hat\iota(0)=0$ and
  $\Ball(\Mat[2])$ is convex, even $\hat\iota(\Md)\subset\Ball(\Mat[2])$ (use
  Proposition~1.6 of~\cite{Meyer:97}).  Put $\phi=\pi\circ\hat\iota$.

  The condition $T_{(2)}[\hat\iota]=\iota$ determines $\hat\iota$ and
  hence~$\phi$ up to first order at~$\pM$.  It is easy to check that
  $\phi(\pM)=0$ and $D\phi(\pM)=(T')^{-1}$.  This is an isometry for the
  \Caratheodory-\Reiffen metric at~$\pM$, so that~$\phi$ is biholomorphic by a
  theorem of \Vigue (Proposition~8.7.2 of~\cite{Jarnicki-Pflug:93}).
  Thus~$\Md$ is biholomorphic to~$\D_2$, contrary to assumption.  Hence~$T$ is
  not \PrMn{2}isometric.
\end{proof}

Examples of isometric but not completely isometric homomorphisms between
Q\PrM{}algebras have also been obtained by Paulsen in~\cite{Paulsen:97}.  Thus,
already for Q\PrM{}algebras, the matrix normed structure on the cotangent space
contains additional information besides the quotient metric.  Notice that the
unital Q\PrM{}algebras involved are \PrMn{3}dimensional and that they are realized
as quotients of $\BH(\Md)$ for what are supposed to be very well\PrM{}behaved
domains in~$\C^2$.

The next goal is to show that every strongly pseudoconvex domain
$\Md\subset\C^d$ has a contractive representation by \PrM{3\times3}matrices
that is not \PrMn{2}contractive.  The idea of the proof is that near a point in
the boundary of~$\Md$, the cotangent spaces~$\CTg_\pM\Md$ look more and more
like~$\CTg_0\D_d$.  This reduces the assertion to Example~\ref{exa:MultBd}.

Let $\Md\subset\C^d$ be a bounded, strongly pseudoconvex domain (with
\PrM{C^2}boundary).  Then the \Caratheodory-\Reiffen metric can be computed
approximately near the boundary of~$\Md$ (see~\cite{Graham:75}, \cite{Ma:91a},
or~\cite{Jarnicki-Pflug:93}).  A consequence of this is that the unit ball of
$\Tg_\pM\Md$ looks more and more like~$\D_d$ for $\pM\to\bd\Md$.  To make this
precise, define the distance of two (bicontinuous) Banach spaces $\NS_1$,
$\NS_2$ to be
$$
        \dist(\NS_1,\NS_2)
  =     \log \bigl(\inf \{ \|\rho\|\cdot\|\rho^{-1}\| \mid
         \mbox{$\rho\colon\NS_1\to\NS_2$ invertible} \} \bigr).
$$
There is an obvious version of this distance for matrix normed spaces:
$$
        \dist_\infty(\NS_1,\NS_2)
  =     \log \bigl(\inf \{ \cbn{\rho}\cdot\cbn{\rho^{-1}} \mid
         \mbox{$\rho\colon\NS_1\to\NS_2$ invertible} \} \bigr).
$$

\begin{theorem}  \label{the:approxCTg}
  Let $\Md\subset\C^d$ be a bounded, strongly pseudoconvex domain.  Then
  $\dist_\infty(\CTg_\pM\Md,\CTg_0\D_d)\to0$ for $\pM\to\bd\Md$.
\end{theorem}

\begin{proof}
  The estimates given in~\cite{Jarnicki-Pflug:93} imply immediately that
  $\dist(\Tg_\pM\Md,\ell^2_d)\to0$ for $\pM\to\bd\Md$.  Clearly,
  $\dist(\NS_1,\NS_2)=\dist(\NS_1',\NS_2')$, so that $\Tg_\pM\Md$ can be
  replaced by $\CTg_\pM\Md$.  To get the assertion for $\dist_\infty$, it is
  convenient to use some intermediate results of Ma's proof in~\cite{Ma:91a}.

\begin{theorem}  \label{the:MdDdComp}
  Let $\Md\subset\C^d$ be a bounded, strongly pseudoconvex domain, choose
  $\epsilon>0$.  Then for all $\pM\in\Md$ sufficiently near to~$\bd\Md$, there
  exist $\phi\in\HO(\Md,\D_d)$ and $\psi\in\HO(\D_d,\Md)$ with $\phi(\pM)=0$,
  $\psi(0)=\pM$, and $D(\phi\circ\psi)|_0= (1-\epsilon)\cdot\ID[\Tg_0\D_d]$.
\end{theorem}

This follows from Theorem~3.3 and Lemma~3.4 of~\cite{Ma:91a}.  Ma's Theorem~3.3
gives maps $\phi,\psi$ as above but with $|\det D(\phi\circ\psi)|_0| \ge 1 -
\epsilon$ instead of $D(\phi\circ\psi)|_0=1-\epsilon$.  By Lemma~3.4
of~\cite{Ma:91a}, this implies $\|D(\phi\circ\psi)|_0 X\|_2 \ge
(1-\epsilon)\|X\|_2$, where $\|\blank\|_2$ denotes the \PrMn{\LP[2]}norm.
Hence the linear mapping that sends $D(\phi\circ\psi)|_0(X)$ to $(1-\epsilon)X$
is contractive with respect to the \PrMn{\LP[2]}norms and thus maps $\D_d$
into~$\D_d$.  Composing~$\phi$ with this endomorphism of~$\D_d$ yields
Theorem~\ref{the:MdDdComp}.

The induced maps $\phi^\ast\colon \CTg_0\D_d\to \CTg_\pM\Md$ and
$\psi^\ast\colon \CTg_\pM\Md\to\CTg_0\D_d$ are obviously completely
contractive.  Moreover, the inverse of~$\phi^\ast$ is~$\psi^\ast$ up to a
factor of~$1-\epsilon$.  Thus
$$
	\exp\bigl(\dist_\infty(\CTg_0\D_d,\CTg_\pM\Md)\bigr)
  \le	\cbn{\phi^\ast} \cbn{(1-\epsilon)^{-1}\psi^\ast}
  \le	(1-\epsilon)^{-1}.
$$
Theorem~\ref{the:approxCTg} follows because this holds for all $\epsilon>0$.
\end{proof}

\begin{corollary}  \label{cor:BadPseudoconvex}
  Let $d>1$ and let $\Md\subset\C^d$ be a bounded, strongly pseudoconvex
  domain.  Then $\BH(\Md)$ has a contractive representation by
  \PrM{3\times3}matrices that is not \PrMn{2}contractive.
\end{corollary}

\begin{proof}
  $\MIN(\ell^2_d)$ has a contractive linear representation $\rho\colon
  \ell^2_d\to\ell^2_2\to\Mat[3]$, where $\ell^2_d\to\ell^2_2$ is some
  orthogonal projection and $\ell^2_2\to\Shift_2(0)$ is the map of
  Example~\ref{exa:MultBd}.  Clearly, $\|\rho\|_{(2)}=\sqrt{2}$.  Choose
  $\epsilon>0$ such that $(1-\epsilon)\sqrt{2}>1$.

  According to Theorem~\ref{the:approxCTg}, for~$\pM$ sufficiently near to some
  boundary point of~$\Md$, there exist completely contractive maps
  $\phi^\ast\colon \CTg_0\D_d\to\CTg_\pM\Md$, $\psi^\ast\colon
  \CTg_\pM\Md\to\CTg_0\D_d$ satisfying $\psi^\ast\circ \phi^\ast =
  (1-\epsilon)\ID$.  Since $\CTg_0\D_d\cong\MIN(\ell^2_d)$, $\rho\circ
  \psi^\ast$ is a contractive representation of $\CTg_\pM\Md$.  It is not
  \PrMn{2}contractive because $\rho\circ (\psi^\ast\circ \phi^\ast)$ is not
  \PrMn{2}contractive.
\end{proof}

\begin{example}  \label{exa:ellOne}
  $\D^2(0)$ has no isometric representation on a finite dimensional Hilbert
  space.
\end{example}

\begin{proof}
  It suffices to show this for $\CTg_0\D^2$, which is isometric to the
  dual~$\ell_2^1$ of~$\ell_2^\infty$ by Theorem~\ref{the:ACQClassify}.  Assume
  that $\phi\colon \ell_2^1\to\Mat$ is a finite dimensional isometric
  representation.  Let $e_1,e_2$ be the standard basis of~$\ell_2^1$, so that
  $\|a_1e_1+a_2e_2\|=|a_1|+|a_2|$.  Let $T_1=\phi(e_1)$, $T_2=\phi(e_2)$.  Then
  $\|T_j\|=\|e_j\|=1$.  Moreover, $\|T_1+\lambda T_2\|=2$ for all
  $\lambda\in\bdD$.  Hence there are vectors $x_\lambda\in \C^n$ with
  $\|x_\lambda\|=1$ and $\|(T_1+\lambda T_2)x_\lambda\|=2$, because the unit
  ball of~$\C^n$ is compact.

  Since the~$T_j$ are contractions, $\|T_1x_\lambda\|=1$ and
  $\|T_2x_\lambda\|=1$.  Moreover, $\lambda T_2 x_\lambda=T_1x_\lambda$ because
  otherwise $\|(T_1+\lambda T_2)x_\lambda\|<2$.  Let $\HilS\subset\C^n$ be the
  \PrMn{1}eigenspace of $T_1^\ast T_1$.  All the vectors~$x_\lambda$ lie
  in~$\HilS$ because $\|T_1 x_\lambda\|=1$ and $\|T_1\|=1$.  Furthermore,
  $x_\lambda\in\HilS$ and $\lambda T_2x_\lambda=T_1 x_\lambda$ imply $\lambda
  T_1^\ast T_2 x_\lambda = T_1^\ast T_1 x_\lambda = x_\lambda$.  Thus
  $x_\lambda$ is an eigenvalue of the operator $T_1^\ast T_2$ with
  eigenvalue~$\conj{\lambda}$.  However, $T_1^\ast T_2\in\Mat$ has only
  finitely many eigenvalues.  Thus~$x_\lambda$ cannot exist for all
  $\lambda\in\bdD$, contradiction.
\end{proof}

\section{Q\PrM{}algebras in complex analysis}
\label{sec:ComplexAnalysis}

Let $\FA\subset\NBC(\Omega)$ be a function algebra.  If~$\NS$ is any normed
vector space, then there is a natural norm on
$\FA\otimes\NS\subset\NBC(\Omega)\otimes\NS$: View elements as functions
from~$\Omega$ to~$\NS$ and take the supremum norm.  This well\PrM{}known tensor
product norm can actually be defined for tensor products of arbitrary normed
spaces, and it has the property that $\|L_1\otimes L_2\|\le\|L_1\|\cdot\|L_2\|$
for any linear maps $L_1\colon \FA_1\to\FA_2$, $L_2\colon \NS_1\to\NS_2$.  In
particular, any contractive linear map $L\colon \FA_1\to\FA_2$ induces
contractions $L_{\NS}=L\otimes\ID[\NS]$.

If $\QA=\FA/\ide$ for some closed ideal~$\ide$, then $\QA\otimes\NS$ carries a
natural quotient norm from the identification $(\FA/\ide)\otimes\NS\cong
(\FA\otimes\NS)/(\ide\otimes\NS)$.  The collection of these norms for all
finite dimensional normed spaces is called the \emph{ac-normed structure}
of~$\QA$ in~\cite{Meyer:97a}.  It can be defined more generally for a quotient
of a normed space by a closed ideal.  Obviously, if a linear map
$\FA_2\to\FA/\ide$ can be lifted to a contractive linear map $\FA_2\to\FA$,
then it must be \emph{ac-contractive}, i.e.\ its tensor product with $\ID[\NS]$
must be contractive for all normed spaces~$\NS$.  Conversely, such a
contractive linear lifting exists for ac-contractive maps if the range is a
dual Banach space and~$\ide$ is weakly closed \cite{Meyer:97a}.

Although it looks formally quite similar to the matrix normed structure for
operator algebras, this ac-normed structure is considerably finer.  This can be
seen most easily from the following example.

Let $\Bd\subset\C^d$, $d>1$, be an absolutely convex domain.  Define~$\UA(\Bd)$
to consist of all bounded, not necessarily continuous functions
$f\colon\Bd\to\C$ such that $f|_{P\cap\Bd}$ is holomorphic for all
\PrMn{1}dimensional planes~$P$ through the origin and $f(z) = a + l(z) +
O(z^2)$ for $z\to0$, with some \PrMn{\C}linear functional~$l$ and $a\in\C$.
This is a function algebra on some compact \Hausdorff space, e.g.\ the
Stone\PrM{}\Cech compactification of $(\Bd,\text{discrete topology})$.

Clearly, $\BH(\Bd)\subset\UA(\Bd)$.  This induces a homomorphism $\theta\colon
\Bd(0)\to\UA(0)$.  Indeed, it is easy to check that this is an isomorphism, the
inverse being
$$
	\UA(0)\ni [a + l(z) + O(z^2)]\mapsto [a+l(z)]\in \Bd(0).
$$

\begin{proposition}  \label{pro:BHUA}
  $\theta\colon \Bd(0)\to\UA(0)$ is a complete equivalence.
\end{proposition}

\begin{proof}
  $\theta$ is obviously completely contractive because it can be lifted to the
  inclusion map $\BH(\Bd)\to\UA(\Bd)$.  To show that~$\theta^{-1}$ is
  completely contractive, it suffices to check this for the restriction to
  $\CTg_0\UA(\Bd)$ by Theorem~\ref{the:unitizationWellDefined}.

  Therefore, take $f\in\Ball(\CTg_0\UA(\Bd)_{(n)})$, i.e., $f\colon
  \Bd\to\Ball(\Mat)$ and $f(0)=0$.  Write $f(z)=l(z)+O(z^2)$ for some linear
  functional~$l$.  An application of the Schwarz lemma to the restriction
  of~$f$ to each \PrMn{1}dimensional plane through the origin yields that
  $\|l(z)\|<1$ for all $z\in\Bd$.  Thus $l\in\HO\bigl(\Bd,\Ball(\Mat)\bigr)$
  is a representative for $\theta^{-1}_{(n)}(f)$ in $\CTg_0\Bd$, so
  that~$\theta^{-1}$ is completely contractive.
\end{proof}

Although the algebra~$\UA(\Bd)$ is quite pathological from the point of view of
complex analysis, the Q\PrM{}algebras $\UA(0)$ and $\Bd(0)$ cannot be
distinguished by dilation theory.  However, they can be distinguished by their
ac-normed structure:

\begin{proposition}  \label{pro:BHUAac}
  $\theta$ is not ac-isometric for $\Bd=\D_2$.
\end{proposition}

\begin{proof}
  Let $(\Md,\pM)$, $\Md\subset\C^2$, be a pointed convex domain as in
  Example~\ref{exa:IsoNotCI} and let $L\colon \Tg_0\D_2\to\Tg_\pM\Md$ be an
  isometry for the \Caratheodory-\Reiffen metric.  For each plane~$P$
  through~$0$ in $\Tg_\pM\Md$, choose an element $X\in P$ with
  $\gamma_\Md(\pM;X)=1$, and choose a complex geodesic $\phi_X\in\HO(\D,\Md)$
  for $(\pM;X)$ according to \Lempert's Theorem \cite{Jarnicki-Pflug:93}.
  Via~$L$, the plane~$P$ corresponds to a plane in $\Tg_0\D_2$, and $\D_2\cap
  L^{-1}(P)=\D$.  View the complex geodesic~$\phi_X$ as a holomorphic map
  $L^{-1}(P)\cap\D_2\to\Md$ with derivative $L|_{L^{-1}(P)}$ at~$0$.  These
  functions can be pieced together to a map $\phi\colon \D_2\to\Md$ whose
  components lie in~$\UA(\D_2)$.

  If~$\theta$ were ac-isometric, there would be a contractive linear map
  $\hat\theta^{-1}\colon \UA(\D_2)\to\BH(\D_2)$ that lifts~$\theta^{-1}$
  by~\cite{Meyer:97a}.  Then $\tilde\phi=\hat\theta^{-1}_\ast(\phi)$ is a
  function in $\HO(\D_2,\Md)$ because~$\Md$ is bounded and convex
  \cite{Meyer:97a}.  Moreover, $\tilde\phi(0)=\phi(0)=\pM$ and
  $D\tilde\phi(0)=D\phi(0)=L$.  Hence~$\tilde\phi$ is biholomorphic by
  Proposition~8.7.2 of~\cite{Jarnicki-Pflug:93}, contradicting the choice
  of~$\Md$.
\end{proof}

The proof Proposition~\ref{pro:BHUAac} goes through for an absolutely convex
domain~$\Bd$ whenever there is a pointed, bounded, convex domain $(\Md,\pM)$
not biholomorphic to~$\Bd$ whose unit ball with respect to the quotient metric
is~$\Bd$.  By~\cite{Lempert:88}, this is the case for all two\PrM{}dimensional
absolutely convex domains with smooth boundary.

The reason why the ac-normed structure contains more information than the
matrix normed structure is Theorem~\ref{the:unitizationWellDefined}: The vector
spaces $\CTg_0\UA(\Bd)$ and $\CTg_0\Bd$ \emph{are} ac-isometric by an obvious
generalization of the proof of Proposition~\ref{pro:BHUA}.  But this does not
imply that the unitizations are ac-isometric.  From the point of view of
complex analysis, this corresponds to the fact that the unit balls of~$\Mat$
have a transitive automorphism group, but general absolutely convex domains do
not.  Thus interpolation with values in arbitrary absolutely convex domains is
more general than interpolation with values in the domains $\Ball(\Mat)$.

From the point of view of complex analysis, restricting attention to
interpolation with values in the domains $\Ball(\Mat)$ does not seem very
fruitful.  On the one hand, this special case is too symmetric to be
``generic''.  It does not contain enough information to tackle more general
interpolation problems whose range is not a symmetric domain.  On the other
hand, this special case is still too complicated for an interesting theory.
This is exemplified by the counterexamples in Section~\ref{sec:Cotangent}.

\section{Arveson's model theory for $d$-contractions}
\label{sec:modelD}

In~\cite{Arveson:97}, Arveson develops a model theory for
\emph{\PrMn{d}contractions}, which are commuting \PrMn{d}tuples of operators
$\mO{T}=(T_1,\dots,T_d)$ on a Hilbert space~$\HilS$
satisfying~\eqref{equ:dcontraction}.  He defines a particular
\PrMn{d}contraction, the \emph{\PrMn{d}shift} $\mO{S}=(S_1,\dots,S_d)$ acting
on a Hilbert space~$H^2_d$ which can be viewed as a variant of Bosonic Fock
space, or as the closure of $\Pol(\cl{\D_d})$ under a somewhat strange norm.
All elements of~$H^2_d$ can be viewed as continuous functions on the closed
Euclidean ball~$\cl{\D_d}$ holomorphic on the interior, but not all such
functions arise.  Let $u_x(z)=(1-\5{z}{x})^{-1}$, where
$\5{z}{x}=z_1\conj{x_1}+\dots +z_d\conj{x_d}$ for $z,x\in\D_d$.  Then $u_x\in
H^2_d$ for all $x\in\D_d$, and
\begin{equation}  \label{equ:fux}
	\5{f}{u_x} = f(x)
\end{equation}
for all $f\in H^2_d$.  Especially,
\begin{equation}  \label{equ:uxuy}
        \5{u_x}{u_y} = (1-\5{y}{x})^{-1}.
\end{equation}
Moreover, the vectors~$\{u_x\}$ span a dense subset of~$H^2_d$.  As remarked by
Arveson, the inner product on~$H^2_d$ does not come from any measure on~$\C^d$
and the \PrMn{d}shift is \emph{not} subnormal.

The \emph{\PrMn{d}Toeplitz algebra}~$\Toep_d$ is the \CstarAlgebra generated by
the \PrMn{d}shift.  The closed subalgebra of~$\Toep_d$ generated by~$S$ is
called $\Shift_d$.  Elements of $\Shift_d$ are often viewed as functions~$f$
on~$\cl{\D_d}$.  Then~$M_f$ denotes the corresponding operator on~$H^2_d$.  It
is shown in~\cite{Arveson:97} that the \PrMn{d}shift has joint spectrum
$\cl{\D_d}$.  Hence the functional calculus yields
$\HO(\cl{\D_d})\subset\Shift_d$.  Indeed, the following inclusions hold:
$$
\HO(\cl{\D_d})\subset \Shift_d\subset H^2_d\subset \clHO(\cl{\D_d}).
$$
The \PrMn{d}shift is important because it is the ``universal''
\PrMn{d}contraction in the following sense:

\begin{theorem}[Arveson \cite{Arveson:97}]  \label{the:Arv1}
  Let $\mO{T}=(T_1,\dots,T_d)$ be a \PrMn{d}contraction on~$\HilS$ and
  let~$\mO{S}$ be the \PrMn{d}shift.  Then $S_j\mapsto T_j$ defines a unital,
  completely contractive representation of~$\Shift_d$ on $\HilS$.  Conversely,
  if $\rho\colon \Shift_d\to\Bound(\HilS)$ is a unital, completely contractive
  representation, then $\rho(\mO{S})=\bigl(\rho(S_1),\dots,\rho(S_d)\bigr)$ is
  a \PrMn{d}contraction.
\end{theorem}

\begin{theorem}[Arveson \cite{Arveson:97}]  \label{the:Arv2}
  Let $d=1,2,\dots$, let $\mO{T}=(T_1,\dots,T_d)$ be a \PrMn{d}contraction
  acting on a separable Hilbert space and let $\mO{S}=(S_1,\dots,S_d)$ be the
  \PrMn{d}shift.  Then there is a triple $(n,\mO{Z},\HilS[K])$ consisting of an
  integer $n=0,1,\dots,\infty$, a spherical operator~$\mO{Z}$ and a full
  co-invariant subspace~$\HilS[K]$ for the operator $n\cdot \mO{S} \oplus
  \mO{Z}$ such that~$\mO{T}$ is unitarily equivalent to the compression of
  $n\cdot \mO{S}\oplus \mO{Z}$ to~$\HilS[K]$.
\end{theorem}

Here a \emph{spherical operator} is a \PrMn{d}tuple $\mO{Z}=(Z_1,\dots,Z_d)$ of
commuting normal operators acting on a common Hilbert space with joint spectrum
$\bd\D_d$, i.e.\ $Z_1Z_1^\ast + \dots +Z_dZ_d^\ast =1$.  $n\cdot\mO{S}$ denotes
the direct sum of~$n$ copies of~$\mO{S}$ acting on $(H^2_d)^n$.  A subspace is
called \emph{co-invariant} if its orthogonal complement is invariant and
\emph{full} (for a collection of operators) if it generates the whole Hilbert
space under the action of the \CstarAlgebra generated by these operators.
Another result of~\cite{Arveson:97} is that this dilation is essentially
unique.  Of course, uniqueness can only hold if one restricts attention to
dilations which involve compression to a full, co-invariant subspace.

It is elementary that every spherical operator is a \PrMn{d}contraction.
Especially, the ``coordinate functions'' $z_j\in\NBC(\bd\D_d)$ form a spherical
operator, and it is easy to check that $S_j\to z_j$ defines a unital
\PrMn{\ast}homomorphism~$\pi$ from~$\Toep_d$ onto $\NBC(\bd\D_d)$.  Its kernel
consists of the compact operators on~$H^2_d$, so that~$\Toep_d$ as a
\CstarAlgebra is not different from more classical higher-dimensional Toeplitz
algebras.  For $d=1$, the \PrMn{1}shift is just the unilateral shift and
$\Shift_1\cong\clHO(\clD)$ is contained in the well\PrM{}known \PrMn{1}Toeplitz
algebra as usual.

The transposed algebra $\Shift_d^t$ is also of interest.  Let
$\mO{S}^t=(S_1^t,\dots,S_d^t)$ be the multioperator $(US_1^\ast U^{-1},\dots,
US_d^\ast U^{-1})$ for any anti-unitary operator $U\colon H^2_d\to H^2_d$.
Then $\Shift_d^t$ is the operator algebra generated by~$\mO{S}^t$.  Call a
multi-operator $\mO{T}=(T_1,\dots,T_d)$ a \emph{transposed \PrMn{d}contraction}
if the matrix
$$
\begin{pmatrix}
 T_1 & 0 & \dots \\
 T_2 & 0 & \dots \\
 \hdotsfor{3} \\
 T_d & 0 & \dots
\end{pmatrix}
$$
is a contraction, i.e.\ $T_1^\ast T_1 + \dots + T_d^\ast T_d\le 1$.
Theorem~\ref{the:transposition} implies that~$\mO{T}$ is a transposed
\PrMn{d}contraction iff~$\mO{T}^t$ is a \PrMn{d}contraction.  The dilation
theory of transposed \PrMn{d}contractions is quite similar to that of ordinary
\PrMn{d}contractions:

\begin{theorem}  \label{the:Arv1transpose}
  Let $\mO{T}=(T_1,\dots,T_d)$ be a transposed \PrMn{d}contraction on~$\HilS$.
  Then $S_j^t\mapsto T_j$ defines a unital, completely contractive
  representation of~$\Shift_d^t$ on $\HilS$.  Conversely, if $\rho\colon
  \Shift_d^t\to\Bound(\HilS)$ is a unital, completely contractive
  representation, then
  $\rho(\mO{S}^t)=\bigl(\rho(S_1^t),\dots,\rho(S_d^t)\bigr)$ is a transposed
  \PrMn{d}contraction.
\end{theorem}

\begin{theorem}  \label{the:Arv2transpose}
  Let $d=1,2,\dots$ and let $\mO{T}=(T_1,\dots,T_d)$ be a transposed
  \PrMn{d}contraction acting on a separable Hilbert space.  Then there is a
  triple $(n,\mO{Z},\HilS[K])$ consisting of an integer $n=0,1,\dots,\infty$, a
  spherical operator~$\mO{Z}$ and a full \emph{invariant} subspace~$\HilS[K]$
  for the operator $n\cdot \mO{S}^t \oplus \mO{Z}$ such that~$\mO{T}$ is
  unitarily equivalent to the compression of $n\cdot \mO{S}^t\oplus \mO{Z}$
  to~$\HilS[K]$.
\end{theorem}

\begin{proof}[Proof of Theorem \ref{the:Arv1transpose}
  and~\ref{the:Arv2transpose}]

If~$\mO{T}$ is a transposed \PrMn{d}contraction, then $S_j\to T_j^t$ defines a
completely contractive representation of $\Shift_d$ by Theorem~\ref{the:Arv1}.
The transpose of this representation is given by $S_j^t\mapsto T_j$.  It is
completely contractive by Lemma~\ref{lem:transposeFunctorial}.  The other half
of Theorem~\ref{the:Arv1transpose} follows similarly.

Theorem~\ref{the:Arv2transpose} follows from Theorem~\ref{the:Arv2} because the
transposes of unitarily equivalent operators are again unitarily equivalent and
because the transpose of a spherical operator is again a spherical operator.
Notice that a co-invariant subspace for~$\mO{S}$ is invariant for~$\mO{S}^\ast$
and hence corresponds under the anti-unitary~$U$ to an invariant subspace
for~$\mO{S}^t$.
\end{proof}

An important question that is not addressed in~\cite{Arveson:97} is whether the
symmetry group of the ball gives rise to symmetries of~$H^2_d$ and thus of
$\Shift_d$.  Recall that~$\D_d$ is a homogeneous domain.  Its symmetry group at
the origin is just the unitary group $\U(d)$.  It is elementary to check that
the unitaries give rise to unitary operators in~$H^2_d$.  In fact, Arveson
shows in~\cite{Arveson:97} that every contraction in~$\Mat[d]$ gives rise to a
contraction in~$H^2_d$ and thus a completely contractive endomorphism of
$\Shift_d$.  In order to see that automorphisms not fixing the origin also give
rise to unitary operators in~$H^2_d$, another characterization of~$H^2_d$ as a
``twisted Bergman space'' is necessary.

\subsection{$H^2_d$ as a twisted Bergman space}
\label{sec:twistedBergman}

The Bergman kernel of the domain~$\D_d$ is \cite{Jarnicki-Pflug:93}
$$
K_{\D_d}(z,w)=\frac{d!}{\pi^d} (1-\5{z}{w})^{-(d+1)}.
$$
Hence, up to a constant factor, $u_z(w)=K_{\D_d}(w,z)^{1/(d+1)}$.
Equation~\eqref{equ:fux} means that $u_z(w)$ is a reproducing kernel for the
Hilbert space~$H^2_d$.  A \emph{twisted Bergman space} is a Hilbert space with
reproducing kernel~$K_{\D_d}^\lambda$ for some $\lambda\in\R$.\footnote{Not
all~$\lambda$ are admissible, i.e.\ give rise to a positive definite kernel.}
Such Hilbert spaces have been studied (also for other symmetric domains) by
harmonic analysts, mainly for the reason that they carry a natural projective
representation of the semi-simple Lie group $\Aut(\D_d)=\PSU(d,1)$:

\begin{theorem}  \label{the:BallAuto}
  Let $h\in\Aut(\D_d)$ be an automorphism of~$\D_d$.  Let
  $$
        \delta(z)
  =     \bigl(\det Dh(z)\bigr)^{1/(d+1)},
  $$
  where any holomorphic branch of the root is chosen and let $(Tf)(z)=
  \delta(z) f\bigl(h(z)\bigr)$ for $f\in H^2_d$, $z\in\D_d$.  Then~$T$ defines
  a unitary operator $H^2_d\to H^2_d$.  This gives rise to a projective
  representation of $\Aut(\D_d)$ on~$H^2_d$.

  Moreover,  $M_{f\circ h}\circ T= T\circ M_f$ for all $f\in\Shift_d$, so that
  $f\mapsto f\circ h$ is a completely isometric automorphism of $\Shift_d$.
\end{theorem}

The proof is based on the behavior of the Bergman kernel under biholomorphic
mappings (\cite{Jarnicki-Pflug:93}, Proposition~6.1.7), which implies
$$
	\bigl(\det Dh(z)\bigr)^{\lambda}
	\conj{\bigl(\det Dh(w)\bigr)^{\lambda}}
	K_{\D_d}\bigl(h(z),h(w)\bigr)^{\lambda}
  =	K_{\D_d}(z,w)^{\lambda}
$$
for all $\lambda\in\R$, $z,w\in\D_d$, and $h\in\Aut(\D_d)$.
See~\cite{Arazy:95}.

However, more general holomorphic mappings $f\colon \D_d\to\D_d$ do
\emph{not} give rise to completely contractive endomorphisms of $\Shift_d$,
although this is the case for linear contractions.  The reason is that any
$f\in\HO(\D_d,\D)$ can be viewed as an endomorphism of~$\D_d$ mapping~$z$ to
$(f(z),0,\dots,0)$.  But if $f\notin\Shift_d$ or $\|f\|_{\Shift_d}>1$, then this
cannot produce a completely contractive endomorphism of $\Shift_d$.

In the article~\cite{Bagchi-Misra:96} by Bagchi and Misra, some interesting
results are proved for the analogue of the \PrMn{d}shift on twisted Bergman
spaces over the symmetric domains $\Ball(\Mat[n,m])$.  They determine when the
generalized shift operator is bounded and they show that its joint spectrum is
$\clBall(\Mat[n,m])$ whenever it is bounded.  Their criterion implies
that~$\mO{S}$ is bounded and that $\Spec \Shift_d=\cl{\D_d}$.  Moreover, they
find necessary and sufficient criteria for the generalized shift operators to
be subnormal.  Their criterion implies that~$\mO{S}$ is not jointly subnormal
and that the inner product on~$H^2_d$ does not come from a measure on~$\C^d$.
However, they fail to notice the special role of the \PrMn{d}shift on~$H^2_d$.
The most important contribution of the twisted Bergman space picture of~$H^2_d$
is the projective representation of $\Aut(\D_d)$.

\subsection{$H^2_d$ and the \Fantappie transform}
\label{sec:Fantappie}

The inner product on~$H^2_d$ is also related to the \emph{\Fantappie
transform}, which is discussed here very briefly, following
\Hoermander~\cite{Hormander:94}.  Let $\Md\subset\C^d$ be a domain.
Endow~$\C^d$ with the usual bilinear form $(x,y)=x_1y_1+\dots +x_dy_d$ and
define $\Md^\dag=\{ x\in\C^d\mid (x,y)\neq1\ \forall y\in\Md\}$.  Notice that
this is a compact set for open~$\Md$ with $0\in\Md$.

If $l\colon \HO(\Md)\to\C$ is a linear functional continuous with respect to
the topology of locally uniform convergence on~$\Md$, its \Fantappie transform
$\Fant{l}\in\HO(\Md^\dag)$ is defined by $\Fant{l}(x)=l(\hat{u}_x)$, where
$\hat{u}_x(y)=\bigl(1-(y,x)\bigr)^{-1}$ for $y\in\Md$.  Actually, $\Fant{l}(x)$
is holomorphic on a neighborhood of~$\Md^\dag$ because the support of~$l$ is a
compact subset of~$\Md$ by continuity.  The main theorem about the
\Fantappie transform is that, if~$\Md$ is \PrMn{\C}convex (especially if~$\Md$
is convex), then~$\Fant$ is a bijection between the space of continuous linear
functionals on $\HO(\Md)$ and the space $\HO(\Md^\dag)$ of functions
holomorphic in a neighborhood of~$\Md^\dag$.  Since $\D_d^\dag=\cl{\D_d}$ and
$\HO(\cl{\D_d})\subset\HO(\D_d)$, this yields a bilinear form on
$\HO(\cl{\D_d})\supset\Rat(\cl{\D_d})\supset\Pol(\cl{\D_d})$ in the special
case $\Md=\D_d$.

However, a sesquilinear form is necessary in order to get a Hilbert space.
Therefore, replace the bilinear form $(,)$ by a sesquilinear form~$\5{}{}$ and
apply the same reasoning.  This yields a ``conjugate \Fantappie'' transform
$\conj{\Fant}\colon l\mapsto l(u_x)$ where~$u_x$ is defined as above, which
maps the dual of $\HO(\Md)$ to the space of conjugate-holomorphic functions on
a neighborhood of the conjugate $(\Md^\dag)^\ast=\{\conj{z}\mid z\in\Md^\dag\}$
of~$\Md^\dag$.  For the ball, this yields a sesquilinear form $B(,)\colon
\HO(\D_d)\times\HO(\cl{\D_d})\to\C$ given by
$B(f,g)=\conj{\Fant}^{-1}(\conj{g})(f)$.

For the simplest linear functionals, the point masses~$\delta_x$ for
$x\in\D_d$, we have $\conj{\Fant}(\delta_x)(z)=u_z(x)=\conj{u_x(z)}$.  Hence
$B(f,u_x)=\delta_x(f)=f(x)$ as for the inner product on~$H^2_d$.  Thus~$B$
coincides with the inner product on~$H^2_d$.

\section{Models for \PrMn{d}contractions with prescribed spectrum}
\label{sec:QuotientMultiplier}

In this section, completely isometric representations of quotients of
$\Shift_d$ are determined.  For $d=1$, this was already done by Arveson
in~\cite{Arveson:69}.

\begin{theorem}  \label{the:modeldContraction}
  Let $\ide\subset\Shift_d$ be a closed ideal.  Let $H^2_d(\ide)= H^2_d \ominus
  \ide\cdot H^2_d\subset H^2_d$ and let
  $\mO{S}(\ide)=\bigl(S_1(\ide),\dots,S_d(\ide)\bigr)$ be the compression of
  the \PrMn{d}shift~$\mO{S}$ to $H^2_d(\ide)$.  Moreover, let~$\SA(\ide)$ (the
  self-adjoint part of~$\ide$) be the quotient of $\NBC(\bd\D_d)$ by the
  closed ideal generated by $\pi(\ide)\subset\NBC(\bd\D_d)$, where $\pi\colon
  \Toep_d\to\NBC(\bd\D_d)$.

  Then~$\SA(\ide)$ is completely isometric to a commutative \CstarAlgebra.  A
  completely isometric representation of $\Shift_d/\ide$ can by obtained as a
  direct sum of a faithful \PrMn{\ast}representation of $\SA(\ide)$ and of the
  representation of $\Shift_d/\ide$ on $H^2_d(\ide)$ induced by the
  \PrMn{d}contraction $\mO{S}(\ide)$.
\end{theorem}

\begin{proof}
  Since any closed ideal in $\NBC(\bd\D_d)$ is necessarily self-adjoint,
  $\SA(\ide)$ is a commutative \CstarAlgebra.

  Let $\rho\colon\Shift_d/\ide\to\Bound(\HilS)$ be any completely isometric
  representation.  Thus $\rho[\mO{S}]=\bigl(\rho[S_1],\dots,\rho[S_d]\bigr)$
  is a \PrMn{d}contraction by Theorem~\ref{the:Arv1}, so that
  Theorem~\ref{the:Arv2} yields certain $(n,\mO{Z},\HilS[K])$.  Call the
  Hilbert space on which~$Z$ acts~$\HilS_{\mO{Z}}$, then $n\cdot \mO{S}\oplus
  \mO{Z}$ acts on $\tilde{\HilS}=(H^2_d)^n\oplus \HilS_{\mO{Z}}$.  Let
  $\hat\rho\colon \Toep_d\to\Bound(\tilde{\HilS})$ be the corresponding
  \PrMn{\ast}representation given by $\hat\rho(S_j)=n\cdot S_j \oplus Z_j$.
  Since~$\HilS[K]$ is co-invariant for $n\cdot \mO{S}\oplus \mO{Z}$, its
  orthogonal complement $\HilS[K]^\bot$ is \PrMn{\hat\rho(\Shift_d)}invariant.

  Let $f\in\ide$ and $\xi\in \tilde{\HilS}$.  Write $\xi=\xi_0+\xi^\bot$ with
  $\xi_0\in\HilS[K]$ and $\xi^\bot\in\HilS[K]^\bot$.  Then
  $\hat\rho(f)\xi^\bot\in\HilS[K]^\bot$ because~$\HilS[K]^\bot$ is
  \PrMn{\hat\rho(\Shift_d)}invariant.  Moreover, the compression of
  $\hat\rho(f)$ to~$\HilS[K]$ is $\rho[f]=\rho(0)=0$, so that
  $\hat\rho(f)\xi_0\in\HilS[K]^\bot$ as well.  Thus
  $\hat\rho(f)\xi\bot\HilS[K]$.  Let
  $\HilS_2=(\ide\tilde{\HilS})^\bot=\{\hat\rho(f)\xi\mid
  f\in\ide,\xi\in\tilde{\HilS}\}^\bot$, then $\HilS[K]\subset\HilS_2$.

  The representation~$\rho$ is obtained from~$\hat\rho$ by compressing
  to~$\HilS[K]$.  Since $\HilS[K]\subset\HilS_2$, first compressing
  to~$\HilS_2$ and then to~$\HilS[K]$ does not change the result.  If
  $f\in\ide$, then the compression
  $\hat\sigma(f)=P_{\HilS_2}\hat\rho(f)P_{\HilS_2}$ is zero by construction.
  Hence~$\hat\sigma$ induces a completely contractive representation
  $\sigma\colon \Shift_d/\ide\to\Bound(\HilS_2)$.  Since the completely
  isometric representation~$\rho$ is a compression of~$\sigma$ to a subspace,
  $\sigma$~must be completely isometric as well.

  By construction, $\tilde{\HilS}$ decomposes into a direct sum of~$n$ copies
  of~$H^2_d$ and a space~$\HilS_{\mO{Z}}$ on which~$\Toep_d$ acts by the
  spherical operator~$\mO{Z}$.  This yields a direct sum decomposition
  of~$\HilS_2$ into~$n$ copies of~$H^2_d(\ide)$ and the part
  $\HilS_{\mO{Z}}\ominus \ide\HilS_{\mO{Z}} = \HilS_{\mO{Z}}\ominus
  \{f(Z)\xi\mid f\in\ide,\xi\in\HilS_{\mO{Z}}\}$.

  Since~$\mO{Z}$ is normal, $\ide\HilS_{\mO{Z}}$ and therefore
  $\HilS_{\mO{Z}}\ominus \ide\HilS_{\mO{Z}}$ is invariant under the
  \CstarAlgebra generated by~$\mO{Z}$.  Therefore, the compression of~$\mO{Z}$
  to $\HilS_{\mO{Z}}\ominus \ide\HilS_{\mO{Z}}$ is still a normal
  (multi)operator with spectrum contained in~$\bd\D_d$, i.e.\ a spherical
  operator.  Thus the representation of $\Shift_d$ on
  $\HilS_{\mO{Z}}\ominus\ide(\mO{Z})\HilS_{\mO{Z}}$ extends to a
  \PrMn{\ast}representation of $\NBC(\bd\D_d)$.  The kernel of this extension
  is a closed ideal of $\NBC(\bd\D_d)$.  It must contain~$\ide$, so
  that~$\mO{Z}$ comes from a \PrMn{\ast}representation of $\SA(\ide)$.

  Hence some completely isometric representation of $\Shift_d/\ide$ can be
  obtained as a direct sum of $n$~copies of $\mO{S}(\ide)$ and a
  \PrMn{\ast}representation of~$\SA(\ide)$.  The representation of~$\SA(\ide)$
  need not be faithful, but replacing it by a faithful representation can only
  increase norms and thus still gives a completely isometric representation.
  Moreover, replacing the~$n$ copies of~$\mO{S}(\ide)$ with just one does not
  change matrix norms either.
\end{proof}

It is not always necessary to add a faithful representation of $\SA(\ide)$.
For example, if $\ide=\{0\}$, then $\mO{S}(\ide)=\mO{S}$ and~$\mO{Z}$ can be
omitted, although $\SA(\ide)=\NBC(\bd\D_d)$.  Indeed, nothing really
interesting happens in the boundary part coming from~$\mO{Z}$.  In the finite
dimensional case, it only contributes a direct sum of several copies of~$\C$ to
the quotient algebra.  To see this, the following lemma is necessary:

\begin{lemma}  \label{lem:quotShiftBd}
  Let $\omega\in\bd\D_d$.  Then $c_{\Shift_d}^\ast(\omega,\omega_2)=1$ for all
  $\omega_2\in\cl{\D_d}\setminus\{\omega\}$ and $\Tg_\omega\Shift_d=\{0\}$.
\end{lemma}

\begin{proof}
  Since rotation by unitaries operates transitively on $\bd\D_d$, we can assume
  without loss of generality that $\omega=(1,0,\dots,0)$.  Since~$\mO{S}$ is a
  \PrMn{d}contraction, $\|S_1\|\le1$.  Moreover, $S_1(\omega)=1$ and
  $S_1(\omega_2)\neq1$ for any $\omega_2\neq\omega$.  Hence
  Corollary~\ref{cor:boundaryTrivial} implies
  $c_{\Shift_d}^\ast(\omega,\omega_2)=1$ for all
  $\omega_2\in\cl{\D_d}\setminus\{\omega\}$ and $\delta(S_1)=0$ for any
  derivation~$\delta$ of $\Shift_d$ at~$\omega$.

  In order to get $\Tg_\omega\Shift_d=\{0\}$, it remains to show that,
  if~$\delta$ is a derivation at~$\omega$, then $\delta(S_j)=0$ also for
  $j=2,\dots,d$.  For then~$\delta$ vanishes on the polynomial algebra which is
  dense in~$\Shift_d$.  If~$\delta$ were a non-zero derivation at~$\omega$,
  then, without loss of generality, $\|\delta\|=1$.  Then the representation
  $$
  \rho\colon f\mapsto
  \begin{pmatrix} f(\omega) & \delta(f) \\ 0 & f(\omega) \end{pmatrix}
  $$
  is a completely contractive representation of~$\Shift_d$.  Indeed, it is a
  completely isometric representation of the quotient $\Shift_d/\Ker\rho$ by
  \PrM{2\times2}matrices.  Hence the matrix
  $$
	\rho\begin{pmatrix} S_1 & S_2 & \cdots & S_d \end{pmatrix}
  =	\begin{pmatrix} 1 & \delta(S_1) & 0 & \delta(S_2) & \cdots
							  & 0 & \delta(S_d) \\
			0 & 1 & 0 & 0 & \cdots & 0 & 0
	\end{pmatrix}
  $$
  is a contraction because $\rho(\mO{S})$ is a \PrMn{d}contraction.  This
  implies $\delta(S_1)=\dots=\delta(S_d)=0$ as desired.  Hence
  $\Tg_\omega\Shift_d=\{0\}$.
\end{proof}

\begin{proposition}  \label{pro:boundaryOrthogonal}
  Let $\ide\subset\Shift_d$ be a closed, finite codimensional ideal.  Then
  $\QA=\Shift_d/\ide$ is completely equivalent to an orthogonal direct sum
  $\Shift_d/\tilde\ide \oplus \C \oplus \dots \oplus \C$, where
  $\tilde\ide\subset\Shift_d$ is a closed ideal with $\SA(\tilde\ide)=\{0\}$ and
  the number of copies of~$\C$ occurring is the dimension of $\SA(\ide)$.
\end{proposition}

\begin{proof}
  The spectrum of $\QA$ can be viewed as a subset of~$\cl{\D_d}$.  By
  Lemma~\ref{lem:quotShiftBd}, a point in the boundary is \PrMn{\QA}related
  only to itself and $\ide(\omega)^\infty=\ide(\omega)$ for $\omega\in\bd\D_d$.
  Hence the assertion follows from Theorem~\ref{the:OrthDirectSum}.
\end{proof}

\begin{corollary}  \label{cor:repDimension}
  Let $\ide\subset\Shift_d$ be a closed ideal of finite codimension~$r$.  Then
  $\Shift_d/\ide$ has a completely isometric representation by \PrM{r\times
  r}matrices.
\end{corollary}

\begin{proof}
  We may assume without loss of generality that $\SA(\ide)=\{0\}$ by
  Proposition~\ref{pro:boundaryOrthogonal}.  By
  Theorem~\ref{the:modeldContraction}, $\Shift_d/\ide$ can be represented
  completely isometrically on $H^2_d\ominus \ide H^2_d$.  Since $\Shift_d$ is
  dense in~$H^2_d$, there is a map $\Shift_d/\ide\to H^2_d\ominus \ide H^2_d$
  whose image is dense.  Since $\Shift_d/\ide$ has finite dimension~$r$, the
  image must be all of $H^2_d\ominus \ide H^2_d$ and this Hilbert space has
  dimension~$r$.
\end{proof}

\begin{theorem}  \label{the:NevPickExplicit}
  Let $x_1,\dots,x_m\in\D_d$ and let $y_1,\dots,y_m\in\Mat$.  Then there
  exists $F\in\Ball\bigl((\Shift_d)_{(n)}\bigr)$ with $F(x_j)=y_j$ for all
  $j=1,\dots,m$ if and only if the block matrix $A\in\Mat[m]\otimes\Mat[n]$
  with entries
  $$
  \frac{1- y_iy_j^\ast}{1-\5{x_i}{x_j}}\in\Mat[n]
  $$
  is positive definite and invertible.  There exists
  $F\in\Cone\bigl((\Shift_d)_{(n)}\bigr)$ with $F(x_j)=y_j$ for all
  $j=1,\dots,m$ if and only if the block matrix $B\in\Mat[m]\otimes\Mat[n]$
  with entries
  $$
  \frac{y_i+y_j^\ast}{1-\5{x_i}{x_j}}\in\Mat[n]
  $$
  is positive definite and invertible.  Moreover, solutions~$F$ can be chosen
  to have polynomial entries.
\end{theorem}

\begin{proof}
  The theorem amounts to a computation of the matrix normed structure of the
  quotient $\QA=\Shift_d/\ide(x_1,\dots,x_n)$.  By
  Theorem~\ref{the:modeldContraction}, $\QA$~has a completely isometric
  representation on the Hilbert space $\HilS=H^2_d\ominus \ide(x_1,\dots,x_m)
  H^2_d$.

  The first step is to write down a (non-orthogonal) basis of~$\HilS$.  We
  claim that the vectors $e_j=u_{x_j}$ for $j=1,\dots,m$ form such a basis.
  Viewing elements of~$H^2_d$ as functions on~$\D_d$, the relation $\5{f}{e_j}=
  f(x_j)$ for all $f\in H^2_d$ shows that the vectors~$e_j$ are all orthogonal
  to the subspace $\ide(x_1,\dots,x_m)H^2_d$ and hence lie in~$\HilS$.
  Moreover, they are linearly independent.  Since $\dim\HilS=m$ by
  Proposition~\ref{cor:repDimension}, they span~$\HilS$.

  Now $\Cone\bigl(\Shift_d/\ide(x_1,\dots,x_n)\bigr)_{(n)}$ can be computed
  using Proposition~\ref{pro:frame}.
  Let $[F]\in\bigl(\Shift_d/\ide(x_1,\dots,x_n)\bigr)_{(n)}$ be given by
  $[F](x_j)=y_j$ for $j=1,\dots,m$.  Let $v_1,\dots,v_n$ be the standard basis
  of~$\C^n$.  Then $\{e_i\otimes v_\mu\}$ is a frame for $\HilS\otimes\C^n$ and
  $$
	\5{[F](e_j\otimes v_\nu)}{e_i\otimes v_\mu}
  =	\5{[F]_{\mu\nu} e_j}{e_i}
  =	[F]_{\mu\nu}(x_i) u_{x_j}(x_i)
  =	(y_i)_{\mu\nu}\cdot (1-\5{x_i}{x_j})^{-1}.
  $$
  This proves the correctness of the criterion for $[F]\in\Cone(\cdots)$.  In
  order to compute the unit ball, however, the action of~$[F]$ must be
  determined in an orthonormal basis.  This will also give a new proof of the
  criterion for $[F]\in\Cone(\cdots)$.

  The inner products between the basis elements are given by the matrix~$B$
  with entries
  $$
  \beta_{i,j}=\5{e_j}{e_i}=(1-\5{x_i}{x_j})^{-1}
  $$
  by~\eqref{equ:uxuy}.  Since the inner product is positive definite, the
  matrix~$B$ is positive and invertible.  Hence the vectors
  $$
	\tilde{e}_j
  =	B^{-1/2}e_j
  =	\sum_{k=1}^m (B^{-1/2})_{kj}e_k
  $$
  are well\PrM{}defined.  It is easy to check that they form an orthonormal basis.
  Moreover, the operator $B\colon e_j\to \sum \beta_{kj}e_k$ still has the
  matrix~$(\beta_{ij})$ in the basis~$(\tilde{e}_j)$ because~$B$ and~$B^{-1/2}$
  commute.

  For $f\in\Shift_d$, the relation $\5{M_f e_i}{e_j} =f(x_j)\5{e_i}{e_j}$ shows
  that the action of~$M_f^\ast$ is given by $M_f^\ast e_j =\conj{f(x_j)}e_j$,
  i.e.\ the basis~$e_j$ is a joint eigenbasis for these adjoints.  Thus the
  action of the compression of~$M_f^\ast$ to~$\HilS$ is given by the matrix
  $$
  B^{1/2}\diag\bigl(\conj{f(x_1)},\dots,\conj{f(x_m)}\bigr) B^{-1/2}
  $$
  in the basis~$(\tilde{e_j})$.  Hence the action of~$\Shift_d$ on~$\HilS$ in
  the orthonormal basis~$(\tilde{e_j})$ is given by
  $$
  f \mapsto B^{-1/2} \diag\bigl( f(x_1),\dots,f(x_m)\bigr)B^{1/2}
  $$
  for $f\in\Shift_d$.

  $[F]$~is represented by $(B\otimes\ID)^{-1/2} \diag(y_1,\dots,y_n)
  (B\otimes\ID)^{1/2}$ on $\HilS\otimes\C^n$.  This matrix has norm less
  than~$1$ iff
  $$
  1- (B\otimes\ID)^{-1/2}\diag(y_j)(B\otimes\ID)\diag(y_j^\ast)
	(B\otimes\ID)^{-1/2}
  $$
  is positive and invertible.  Since $B\otimes\ID$ is invertible, this is
  equivalent to $B\otimes\ID - \diag(y_j)(B\otimes\ID)\diag(y_j^\ast)$ being
  positive and invertible.  This is just the matrix in the statement of the
  theorem, proving the first assertion.

  The real part of~$[F]$ is represented by the matrix
  $$
  (1/2) \cdot (B\otimes\ID)^{-1/2}\diag(y_j)(B\otimes\ID)^{1/2} +
  (1/2) \cdot (B\otimes\ID)^{1/2}\diag(y_j^\ast)(B\otimes\ID)^{-1/2}
  $$
  This is positive and invertible iff the matrix
  $$
  \diag(y_j)(B\otimes\ID) +
  (B\otimes\ID)\diag(y_j^\ast)
  $$
  is positive and invertible.  Again, this is the matrix occurring in the
  statement of the theorem, proving again the second assertion.
\end{proof}

\begin{theorem}  \label{the:NevPickExplicitTranspose}
  Let $x_1,\dots,x_m\in\D_d$ and let $y_1,\dots,y_m\in\Mat$.  Then there
  exists $F\in\Ball\bigl((\Shift_d^t)_{(n)}\bigr)$ with $F(x_j)=y_j$ for all
  $j=1,\dots,m$ if and only if the block matrix $A\in\Mat[m]\otimes\Mat[n]$
  with entries
  $$
  \frac{1- y_i^\ast y_j}{1-\5{x_j}{x_i}}\in\Mat[n]
  $$
  is positive definite and invertible.  There exists
  $F\in\Cone\bigl((\Shift_d^t)_{(n)}\bigr)$ with $F(x_j)=y_j$ for all
  $j=1,\dots,m$ if and only if the block matrix $B\in\Mat[m]\otimes\Mat[n]$
  with entries
  $$
  \frac{y_i^\ast +y_j}{1-\5{x_j}{x_i}}\in\Mat[n]
  $$
  is positive definite and invertible.  Moreover, solutions~$F$ can be chosen
  to have polynomial entries.
\end{theorem}

\begin{proof}
  Since such an $F\in(\Shift_d^t)_{(n)}$ exists iff there is
  $F^t\in(\Shift_d)_{(n)}$ with $F^t(x_j)=y_j^t$ for $j=1,\dots,m$, this is an
  immediate consequence of Theorem \ref{the:NevPickExplicit}.
\end{proof}

For $d=1$, Theorem~\ref{the:NevPickExplicit} is equivalent to the existence
part of Nevanlinna-Pick theory.  However, the proof above is not constructive.

\begin{problem}
  Find a constructive proof of Theorem~\ref{the:NevPickExplicit}.  Is there a
  criterion for the existence of $F\in\Shift_d\otimes\Mat$ with
  $\|F\|_{(n)}\le1$ and prescribed values in finitely many points of~$\D_d$?
\end{problem}

Since $\Shift_1\cong\clPol(\clD)$ is a uniform algebra,
$\Shift_1^t\cong\Shift_1$.  Thus quotients of $\Shift_1^t$ and $\Shift_1$ are
completely equivalent.  This is not so clear, however, from the formulas in
Theorem \ref{the:NevPickExplicit} and~\ref{the:NevPickExplicitTranspose}.

If $\Shift_d/\ide$ has non-trivial tangent space, the orthogonal complement of
$\ide\cdot H^2_d$ can be computed most easily using the characterization of the
inner product by the \Fantappie transform.  There are additional conditions
$l(f)=0$, where the linear functional~$l$ is a differential operator $l=\sum
c_\alpha \partial^\alpha|_{a}$ at some $a\in\D_d$.  Since~$l$ is a continuous
linear functional on $\HO(\D_d)$, the \Fantappie transform yields
\begin{equation}  \label{equ:lf}
        l(f)
  =     \conj{\Fant}^{-1}\conj{\Fant}(l)(f)
  =     \5{f}{\conj{\conj{\Fant}(l)}}_{H^2_d}.
\end{equation}
Hence the function $z\mapsto \conj{\conj{\Fant}(l)(z)}= \conj{l(u_z)}$ lies in
the orthogonal complement of~$\ide$.  If~$l$ is some differential operator as
above, this function can easily be computed.  Appropriate differential
operators~$l$ provide a basis of $H^2_d\ominus \ide H^2_d$.  The inner products
between these vectors can be computed using~\eqref{equ:lf}.  More generally,
\begin{equation}  \label{equ:ll}
	\5{M_f \conj{l_1(u_z)}}{\conj{l_2(u_z)}}
  =	\5{M_f \conj{l_1(u_z)}}{\conj{\conj{\Fant}(l_2)}}
  =	l_2\bigl( f\cdot \conj{l_1(u_z)}\bigr),
\end{equation}
for any continuous linear functionals $l_1,l_2$ on~$\HO(\D_d)$ and
$f\in\Shift_d$.

However, the method of the proof of Theorem~\ref{the:NevPickExplicit} does not
apply to this situation; there is no natural basis for $H^2_d\ominus \ide
H^2_d$ in which the action of~$M_f^\ast$ can be computed easily.  However, the
following recipe still works.  Let $\ide\subset\Shift_d$ be a closed ideal of
finite codimension~$r$ such that $\SA(\ide)=\{0\}$.  Let $\QA=\Shift_d/\ide$
and $F\in\QA_{(n)}$.  In order to determine whether $F\in\Ball(\QA_{(n)})$, do
the following:

\begin{enumerate}
\item Choose a basis $l_1,\dots,l_r$ for the vector space of differential
  operators annihilating~$\ide$ and $g_1,\dots,g_r\in\Shift_d$ representing a
  basis of~$\QA$.  Write~$F$ as a matrix with entries
  $$
	F_{\mu\nu}
  =	\sum_{k=1}^r F_{\mu\nu}^k [g_k].
  $$

\item Compute $\lambda_j=\conj{\conj{\Fant}(l_j)}$ for $j=1,\dots,r$; these
  functions form a basis for $H^2_d\ominus \ide H^2_d$.

\item Using~\eqref{equ:ll}, compute the inner products
  \begin{align*}
	\beta_{ij}
  &=	\5{\lambda_j}{\lambda_i}
   =	l_i(\lambda_j),\\
	\gamma_{kij}
  &=	\5{g_k\lambda_j}{\lambda_i}
   =	l_i(g_k\lambda_j).
  \end{align*}

\item Let $M(F)$ be the block matrix with $\mu,\nu$th entry $\left(\sum_k
  \gamma_{kij} F_{\mu\nu}^k\right)_{ij}\in \Mat[r]$ and let $B=(\beta_{ij})$.

\item Check whether the matrix $B\otimes\ID - M(F) (B^{-1}\otimes\ID)
  M(F)^\ast$ is positive definite and invertible.  This happens iff
  $F\in\Ball(\QA_{(n)})$.

\end{enumerate}

In order to determine whether $F\in\Cone(\QA_{(n)})$, it suffices to compute
the matrix $M(F)$, the matrix~$B$ is not necessary: $F\in\Cone(\QA_{(n)})$ iff
$\RE M(F)$ is positive and invertible.

The proof that the above algorithm works is left to the reader.  It is also
left to the reader to check that it gives the same answer in the special case
of Theorem~\ref{the:NevPickExplicit}.  It is essential to flip the indices
$i,j$ in the definition of $\beta_{ij}$ and~$\gamma_{kij}$ in order to get the
\emph{matrix} normed structure right.  This is because if the matrix~$A$ has
entries~$A_{ij}$ in the orthonormal basis~$\{E_j\}$, then
$\5{AE_i}{E_j}=A_{ji}$.

It is of special interest to compute the two\PrM{}dimensional quotients of the
multiplier algebra.  First look at the quotient by $\ide(x,y)$ with
$x,y\in\cl{\D_d}$.  If one of the points lies in the boundary, then
$c_{\Shift_d}^\ast(x,y)=1$ by Lemma~\ref{lem:quotShiftBd}.  If $x,y\in\D_d$,
the quotient is represented by \PrM{2\times2}matrices.  The adjoints of the
representing matrices have eigenvectors $u_{x}$ and~$u_{y}$.  The angle between
these two vectors is
$$
        \frac{|\5{u_x}{u_y}|}{\|u_x\|\cdot \|u_y\|}
  =     \frac{\sqrt{1-\|x\|^2}\sqrt{1-\|y\|^2}}{|1-\5{x}{y}|}.
$$
Comparing this with the angle between the eigenvectors of the
matrices~$T_c^\ast$ of Section~\ref{sec:twodimOpalg} shows that this number is
$\sqrt{1-c_{\Shift_d}^\ast(x,y)^2}$, so that
\begin{equation}  \label{equ:CaraDistMult}
        c_{\Shift_d}^\ast(x,y)
  =     \left( 1 - \frac{(1-\|x\|^2)(1-\|y\|^2)}{|1-\5{x}{y}|^2} \right)^{1/2}.
\end{equation}
This coincides with the classical \Caratheodory{}\Ast distance for the unit
ball, computed using $\BH(\D_d)$ instead of $\Shift_d$ \cite{Jarnicki-Pflug:93}.
This is not too surprising because $\Shift_d$ is constructed to model
\PrMn{d}contractions, and the quotient distance is constructed to describe
two\PrM{}dimensional quotients.  Multi\PrM{}operators that generate
two\PrM{}dimensional, unital operator algebras can be modeled just as well by
$\BH(\D_d)$.  Only for higher dimensional quotients does $\BH(\D_d)$ fail to
give a satisfactory theory.  Note also that the above distance is invariant
under $\Aut(\D_d)$ because the classical \Caratheodory{}\Ast distance is
invariant.  This must happen because of Theorem~\ref{the:BallAuto}.

Now let us compute the quotient metric for $\Shift_d$.  For this purpose, it is
convenient to use Theorem~\ref{the:BallAuto}, which implies that all
automorphisms of~$\D_d$ induce isometries for the quotient metric.  It is easy
to see that $\Tg_a\Shift_d\cong\Tg_a\D_d\cong\C^d$ for $a\in\D_d$.  Moreover,
Lemma~\ref{lem:quotShiftBd} yields $\Tg_a\Shift_d=\{0\}$ for $a\in\bd\D_d$.
Let $X=\partial/\partial z_1\in\Tg_0\D_d$, then every element of $\Tg\D_d$ is
mapped to a multiple of~$X$ by some automorphism of~$\D_d$ (first map the base
point to the origin and then rotate by a unitary).  Hence the value of
$\gamma_{\Shift_d}(0;X)$ determines the quotient metric completely.

The orthogonal complement of $\ide(0;X) H^2_d$ is spanned by the functions
$1,z_1\in H^2_d$.  Indeed, $\{1,z_1\}$ is an orthonormal basis of $H^2_d\ominus
\ide(0;X)H^2_d$.  Moreover, the compression of~$M_{z_1}$ is the matrix~$T_0$ of
Section~\ref{sec:twodimOpalg} in the basis $\{z_1,1\}$.  Thus the derivation
$[f]\mapsto Df(0;X)$ corresponds to the derivation $T_0\mapsto 1$ and therefore
has norm~$1$.  Consequently, $\gamma_{\Shift_d}(0;X)=1$.
Lemma~\ref{lem:balCRmetric} gives the same result for the classical
\Caratheodory-\Reiffen metric $\gamma_{\BH(\D_d)}(0;X)$.  Since both
$\gamma_{\BH(\D_d)}(a;X)$ and $\gamma_{\Shift_d}(a;X)$ are invariant under
automorphisms of~$\D_d$, these functions coincide everywhere.  Hence the
formula in~\cite{Jarnicki-Pflug:93} for the classical \Caratheodory-\Reiffen
metric on~$\D_d$ can be copied:
$$
        \gamma_{\Shift_d}(a;X)
  =     \left( \frac{\|X\|^2}{1-\|a\|^2} + \frac{|\5{a}{X}|^2}{(1-\|a\|^2)^2}
        \right)^{1/2}.
$$
Here~$\|X\|$ stands for the norm of $X\in\ell^2_d$, not for the norm of the
associated linear functional on~$\Shift_d$, of course.  However,
$\CTg_0\Shift_d$ and $\CTg_0\D_d$ are no longer \emph{completely} isometric for
$d\ge 2$.

\begin{example}  \label{exa:ShiftAtZero}
  The representation of $\Shift_d(0)$ on $H^2_d\ominus \ide(0)^2
  H^2_d$ is completely isometric by Theorem~\ref{the:modeldContraction}.
  Clearly, this orthogonal complement is spanned by the constant function~$1$
  and the linear functions $z_1,\dots,z_d$.  Thus $\CTg_0\Shift_d$ is
  represented by the \PrMn{d}dimensional subspace
  $$
  \begin{pmatrix}
	0    & 0 & 0 & \dots \\
	\ast & 0 & 0 & \dots \\
	\ast & 0 & 0 & \dots \\
	\hdotsfor{4}
  \end{pmatrix}
  $$
  of $\Mat[d+1]$.  Especially, the algebra called $\Shift_2(0)$ in
  Section~\ref{sec:Cotangent} is indeed completely equivalent to
  $\Shift_2/\ide(0)^2$.

  The proof of Theorem~\ref{the:transpositionBad} can easily be generalized to
  show that transposition $\Shift_d(0)\to\Shift_d(0)^t$ is not
  \PrMn{2}isometric.  Thus $\Shift_d(0)$ is not a Q\PrM{}algebra by
  Corollary~\ref{cor:Qtranspose}.  Therefore, $\Shift_d$ cannot be a function
  algebra.  Thus the \PrMn{d}shift cannot be subnormal.
\end{example}

\begin{theorem}  \label{the:MultQutDim}
  Let~$\QA$ be a \PrMn{d'}dimensional quotient of $\Shift_d$.  Then~$\QA$ is
  completely equivalent to quotients of $\Shift_e$ for all $e\ge
  \min\{d,d'-1\}$.
\end{theorem}

\begin{proof}
  Viewing $(S_1,\dots,S_d,0)$ as a \PrMn{d+1}contraction on~$H^2_d$ gives a
  completely contractive unital homomorphism $i_{\downarrow}\colon
  \Shift_{d+1}\to\Shift_d$.  Viewing $(S_1,\dots,S_d)$ as a \PrMn{d}contraction
  on~$H^2_{d+1}$ produces a completely contractive unital homomorphism
  $i_{\uparrow}\colon \Shift_d\to\Shift_{d+1}$.  Clearly, $i_{\downarrow}\circ
  i_{\uparrow}=\ID[\Shift_d]$.  Thus~$i_\uparrow$ is a complete isometry
  and~$i_{\downarrow}$ is a complete quotient map.  Consequently, every
  quotient of $\Shift_d$ is also a quotient of $\Shift_{d+1}$, so that it remains
  to show that~$\QA$ is completely equivalent to a quotient of $\Shift_{d'-1}$
  if $d'\le d<\infty$.

  Since arbitrary finite dimensional quotients of $\NBC(\bd\D_d)$ are
  also quotients of $\Shift_1$, it can be assumed that the spectrum of~$\QA$ is
  not contained in~$\bd\D_d$.  Here view
  $\Spec\QA\subset\Spec(\Shift_d)\cong\cl{\D_d}$ via the map induced by the
  quotient map $\rho\colon \Shift_d\to\QA$.  Using an automorphism of $\Shift_d$
  according to Theorem~\ref{the:BallAuto}, it can be achieved that
  $0\in\Spec\QA$.  Let $\delta=d-d'+1$.

  The space spanned by $1,\rho(S_1),\dots,\rho(S_d)$ in~$\QA$ has dimension at
  most~$d'$, hence there must be at least $d-d'$ independent linear dependence
  relations among these elements.  Since $0\in\Spec\QA$ and $S_j(0)=0$ for all
  $j=1,\dots,d$, these relations cannot involve~$1$ with a nonzero
  coefficient.  Hence there are $d-d'$ linearly independent functionals~$l_j$
  on~$\C^d$ with $l_j\bigl(\rho(S_1),\dots,\rho(S_d)\bigr)=0$ for all
  $j=1,\dots,d-d'$.  By a unitary transformation (and thus a complete isometry
  of $\Shift_d$) it can be achieved that the kernel of these linear functionals
  is spanned by $e_d,e_{d-1},\dots,e_{d-d'}$.

  Hence it can be assumed without loss of generality that
  $\rho(S_{d-d'})=\dots=\rho(S_{d})=0$, so that $\rho=\rho\circ
  i_\uparrow^\delta \circ i_\downarrow^\delta$.  Here, of course,
  $i_\downarrow^\delta= i_\downarrow\circ \dots \circ i_\downarrow$ with
  $\delta$~factors, so that
  $$
        i_\uparrow^\delta\circ i_\downarrow^\delta(S_1,\dots,S_d)
  =     (S_1,\dots,S_{d'-1},0,\dots,0).
  $$
  Now define $\rho_\downarrow=\rho\circ i_\uparrow^\delta\colon
  \Shift_{d'-1}\to\QA$.  By definition, this is completely contractive, and
  $\rho=\rho\circ i_\uparrow^\delta \circ i_\downarrow^\delta = \rho_\downarrow
  \circ i_\downarrow^\delta$.

  If $l_1,l_2$ are contractive linear maps and $l_1\circ l_2$ is an isometry,
  then~$l_2$ must be an isometry.  The dual statement of this is that if
  $l_1\circ l_2$ is a quotient map, then~$l_1$ must be a quotient map.  This
  dual statement can be proved either directly or by realizing that $l\colon
  \NS_1\to\NS_2$ is a quotient map iff $l^\ast\colon \NS_2'\to\NS_2'$ is an
  isometry.  Of course, these statements remain true if ``contractive'',
  ``isometry'', and ``quotient map'' are replaced by ``completely
  contractive'', ``completely isometric'', and ``complete quotient map''
  everywhere.  Especially, this can be applied to $\rho=\rho_\downarrow\circ
  i_\downarrow^\delta$ to obtain that~$\rho_\downarrow$ is a complete quotient
  map.
\end{proof}

\section{The quotient complexity of a commutative operator algebra}
\label{sec:Complexity}

\begin{definition}
  Let~$\OA$ be a unital, commutative operator algebra.  Then~$\OA$ is said to
  have \emph{minimal (quotient) complexity} if every \PrMn{r}dimensional
  quotient of~$\OA$ has a completely isometric representation by \PrM{r\times
  r}matrices.
\end{definition}

The idea behind the concept of quotient complexity is that a subalgebra
of~$\Mat[r']$ for $r'>r$ can have a more complicated matrix normed structure
than an isomorphic subalgebra of~$\Mat[r]$.  Moreover, an \PrMn{r}dimensional,
algebraically semisimple, commutative operator algebra (i.e.\ a direct sum
of~$r$ copies of~$\C$) cannot be represented faithfully in~$\Mat[r-1]$.
However, it is easy to find nilpotent commutative subalgebras of~$\Mat[r]$ of
dimension greater than~$r$:

\begin{example}  \label{exa:complexityMatr}
  Consider~$\Mat[r]$ as an operator space, then its trivial unitization is by
  definition a subalgebra of $\Mat[r+1]$ of dimension $r^2+1$.  However, it
  follows from Lemma~\ref{lem:ellTwoSpecial} that the unitization of~$\Mat[r]$
  does \emph{not} have minimal quotient complexity.
\end{example}

Nevertheless, the term \emph{minimal} seems appropriate, especially if the
following conjecture should turn out to be true: 

\begin{conjecture}
  A unital, commutative subalgebra of~$\Mat[r]$ of dimension greater than~$r$
  does \emph{not} have minimal quotient complexity.
\end{conjecture}

Theorem~\ref{the:cotangentMinimalComplexity} proves this conjecture for trivial
unitizations.

\begin{remark}
  Let $\OA\subset\Mat[n-1]$ be a unital, commutative operator algebra.
  Then~$\OA$ has a unital, completely isometric representation by \PrM{n\times
  n}matrices.
\end{remark}

\begin{proof}
  It is trivial to find a non-unital completely isometric representation
  $\OA\to\Mat$.  Restrict this to any maximal ideal and take the unitization of
  this representation.  By Theorem~\ref{the:unitizationWellDefined}, this
  defines a completely isometric, unital representation of~$\OA$.
\end{proof}

\begin{proposition}  \label{pro:complexitySum}
  Let $\OA_1$ and~$\OA_2$ be unital, commutative operator algebras.  Then
  $\OA_1\oplus\OA_2$ has minimal complexity if and only if both $\OA_1$
  and~$\OA_2$ have minimal complexity.  Moreover, $\OA^t$ has minimal
  complexity if and only if~$\OA$ has minimal complexity.
\end{proposition}

\begin{proof}
  Let $e_j\in\OA_j$ be the identity elements.  Any ideal
  $\ide\subset\OA_1\oplus\OA_2$ is of the form $\ide=\ide_1\oplus\ide_2$ with
  ideals $\ide_j\subset\OA_j$.  Hence $(\OA_1\oplus\OA_2)/\ide\cong
  (\OA_1/\ide_1)\oplus(\OA_2/\ide_2)$.  If $\rho_j\colon
  \OA_j/\ide_j\to\Mat[n_j]$, $j=1,2$, are completely isometric representations,
  then $\rho_1\oplus\rho_2\colon (\OA_1\oplus\OA_2)/\ide\to\Mat[n_1+n_2]$ is a
  completely isometric representation.  Thus if $\OA_1$ and~$\OA_2$ have
  minimal complexity, so has their direct sum.  The reverse implication is
  trivial because every quotient of~$\OA_j$, $j=1,2$, is completely equivalent
  to a quotient of $\OA_1\oplus\OA_2$.

  By Lemma~\ref{lem:transposeFunctorial}, every quotient of~$\OA^t$ is the
  transpose of a quotient of~$\OA$.  Hence it has low-dimensional
  representations iff this is true for the corresponding quotient of~$\OA$ by
  the concrete description of the transpose operation.
\end{proof}

Thus it suffices to study indecomposable algebras of minimal complexity.
Corollary~\ref{cor:repDimension} can be rephrased as follows: All quotients of
$\Shift_d$ have minimal quotient complexity.  The proof uses only formal
properties of the dilation theory for \PrMn{d}contractions.  Hence in order to
find other model theories with similar formal properties, it is interesting to
know whether there are more algebras of minimal complexity.  The only finite
dimensional, indecomposable algebras of minimal complexity known at the moment
are the quotients of $\Shift_d$ and $\Shift_d^t$.  It is quite conceivable that
there are no other examples (Conjecture~\ref{conj:classifyMC}).

Other infinite dimensional operator algebras of minimal quotient complexity are
easy to obtain.  An obvious candidate is the injective limit of the algebras
$\Shift_d$ for $d\to\infty$.  More generally, if~$M$ is a set and~$\Lambda$ is
the net of finite dimensional subsets of~$M$, we can associate to it an
injective system $S\mapsto \Shift_{\# S}$ for $S\in\Lambda$ with the obvious
maps used already in the proof of Theorem~\ref{the:MultQutDim}.  If~$M$ is
uncountable, then the corresponding injective limit will be an operator algebra
of minimal quotient complexity that is not separable.  Indeed, any finite
dimensional quotient of this injective limit is also a quotient of some
$\Shift_d$.

\begin{example}  \label{exa:Qzero}
  The \PrMn{4}dimensional, unital operator algebra
  $\QA_0\otimes\QA_0\subset\Mat[4]$ does not have minimal quotient complexity:
  Its unique \PrMn{2}dimensional cotangent space is invariant under
  transposition.  Hence it cannot be completely isometric to $\Shift_2(0)$ or
  $\Shift_2(0)^t$.
\end{example}

\begin{remark}
  The previous example shows that the spatial tensor product does \emph{not}
  preserve quotient complexity.  Another example is
  $\Pol(\clD^2)\cong\Pol(\clD)\otimes\Pol(\clD)$.  This function algebra has a
  quotient with no finite dimensional completely isometric representations by
  Example~\ref{exa:ellOne}, whereas $\Pol(\clD)$ has minimal quotient
  complexity because it is the algebra generated by the \PrMn{1}shift.  This
  example shows that the maximal tensor product does not preserve quotient
  complexity either.

  However, if $\OA_j\subset\Mat[n_j]$, $j=1,2$, are \PrMn{n_j}dimensional
  subalgebras, then $\OA_1\otimes\OA_2\subset\Mat[n_1\cdot n_2]$ is an
  \PrMn{n_1\cdot n_2}dimensional subalgebra.  Hence the reason for the above
  problem is that taking tensor products is not well\PrM{}behaved with respect
  to quotients.
\end{remark}

\begin{example}
  I expect that the spatial tensor product $\QA_c\otimes\QA_d$ does \emph{not}
  have minimal quotient complexity if $c,d<1$.  For $c=d=0$, this is shown in
  Example~\ref{exa:Qzero}.  It can be shown that $\QA_c\otimes\QA_d$ is
  indecomposable and not a quotient of $\Shift_1$ if $c,d<1$.  Hence the
  general case would follow from Conjecture~\ref{conj:classifyMC}.  It is
  difficult to compute quotients of $\QA_c\otimes\QA_d$ directly, however.
\end{example}

The best structure theorem for algebras of minimal complexity that we can prove
at the moment is the following:

\begin{theorem}  \label{the:cotangentMinimalComplexity}
  Let~$\OA$ be a closed, commutative, unital operator algebra of minimal
  quotient complexity and let $\omega\in\Spec(\OA)$.  Then $\CTg_\omega\OA$ is
  completely equivalent to $\Bound(\C,\HilS)$ or $\Bound(\HilS,\C)$ for some
  Hilbert space~$\HilS$.  Especially, if $\dim\CTg_\omega\OA=d\in\N$, then
  $\CTg_\omega\OA$ is completely equivalent to $\CTg_0\Shift_d$ or
  $\CTg_0\Shift_d^t$.
\end{theorem}

\begin{proof}
  Let~$\HilS$ be a Hilbert space and choose a unit vector $x\in\HilS$.  This
  induces an isometric embedding $\C\subset\HilS$ and a projection
  $\HilS\to\C$, which turn $\Bound(\C,\HilS)$ and $\Bound(\HilS,\C)$ into
  closed subspaces of $\Bound(\HilS)$.  The resulting abstract operator space
  structure on $\Bound(\HilS,\C)$ and $\Bound(\C,\HilS)$ does not depend on the
  choice of~$x$, of course.  Since $\Bound(\C,\HilS)$ and $\Bound(\HilS,\C)$
  are clearly isometric to a Hilbert space, the first step is to show that
  $\CTg_\omega\OA$ is isometric to a Hilbert space.

\begin{lemma}  \label{lem:ellTwoSpecial}
  Assume that every \PrMn{3}dimensional quotient of the unital, commutative
  operator algebra~$\OA$ has an isometric representation on~$\C^3$.  Then, for
  any $\omega\in\Spec(\OA)$, $\CTg_\omega\OA$ is isometric to a pre-Hilbert
  space.
\end{lemma}

\begin{proof}
  The only \PrMn{3}dimensional subalgebras of~$\Mat[3]$ that are trivial
  unitizations are $\Shift_2(0)$ and $\Shift_2(0)^t$, and both are isometric
  to~$\Unze{(\ell^2_2)}$.  By Theorem~\ref{the:unitizationWellDefined}, there
  is a one\PrM{}to\PrM{}one correspondence between quotients of
  $\OA(\omega)=\Unze{(\CTg_\omega\OA)}$ and quotients of~$\CTg_\omega\OA$.
  Hence any two\PrM{}dimensional quotient of $\CTg_\omega\OA$ is isometric
  to~$\ell^2_2$.

  The quotient map $\CTg_\omega\OA\to\ell_2^2$ dualizes to an isometric
  embedding $\ell_2^2\to(\CTg_\omega\OA)'=\Tg_\omega\OA$.  Clearly, any two
  elements of~$\Tg_\omega\OA$ lie in the image of such a map.  Thus the
  parallelogram identity holds in~$\Tg_\omega\OA$, because it only involves
  vectors in a two\PrM{}dimensional subspace.  Hence the norm on~$\Tg_\omega\OA$
  comes from a pre-inner product.  Another dualization yields
  that~$\CTg_\omega\OA$ is a pre-Hilbert space as well.
\end{proof}

The next step is to study subspaces of~$\Mat[n,m]$ that are isometric to a
Hilbert space.

\begin{theorem}  \label{the:ellTwoMat}
  Let $n,m\ge2$ and let $\HilS\subset\Mat[n,m]$ be a subspace of dimension~$r$
  that is isometric to a Hilbert space.  Then $r\le n+m-2$.
\end{theorem}

\begin{remark}
  The bound $n+m-2$ in Theorem~\ref{the:ellTwoMat} probably is not optimal.
  The only candidates of Hilbert spaces contained in $\Mat[n,m]$ that
  immediately come to mind have dimensions $n$ and~$m$, respectively, so a
  likely conjecture is that even $r\le \max\{n,m\}$.  However, this stronger
  estimate is more difficult and not relevant for our purposes.
\end{remark}

\begin{proof}
  Assume the contrary, then there exists a subspace $\HilS\subset\Mat[n,m]$
  that is isometric to $\ell^2_r$ with $r=n+m-1$.  Since transposition
  $\Mat[n,m]\to\Mat[m,n]$ is isometric, we can assume without loss of
  generality that $m\le n$.

  The proof depends on the singular value decomposition of a matrix.  Every
  $A\in\Mat[n,m]$ can be written as $A=U\diag(a_1,\dots,a_m)V$, where
  $U\in\Mat[n,n]$ and $V\in\Mat[m,m]$ are unitary matrices and $\|A\|=a_1\ge
  a_2\ge \dots \ge a_m\ge0$.  Here $\diag(a_1,\dots,a_m)$ stands for the
  \PrM{n\times m}matrix with $i,j$th entry $\delta_{i,j}a_j$.  The matrices~$U$
  and~$V$ are usually not unique, but the ``singular values'' $a_1,\dots,a_n$
  are.  Indeed, they are the eigenvalues of the matrix $(A^\ast A)^{1/2}$.
  This also shows that~$a_k$ depends continuously on~$A$.  Let $\alpha_k\colon
  \Mat[n,m]\to\R_+$ be the map $A\mapsto a_k$.

  For a normed space~$\NS$, let $S(\NS)=\{ \ns\in\NS\mid \|\ns\|=1 \}$.
  Every element of~$S(\HilS)$ has the singular value~$1$, possibly with
  multiplicity.  Choose $A_0\in S(\HilS)$ for which this multiplicity is
  minimal.

  Since left and right multiplication by unitary matrices induces an isometry
  of $\Mat[n,m]$, we can assume without loss of generality that~$A_0$ is of the
  form
  $$
	A_0=\diag\bigl(1,\alpha_2(A_0),\dots,\alpha_m(A_0)\bigr).
  $$
  Let $N\cong\Mat[n-1,m-1]$ be the subspace of~$\Mat[n,m]$ of matrices with
  zeroes in the first row and first column.

\begin{claim}  \label{claim:a}
  The natural map $p\colon \HilS\to\Mat[n,m]/N$ is a vector space isomorphism.
\end{claim}

\begin{proof}
  Since $\dim (\Mat[n,m]/N)=n+m-1=\dim \HilS$, it suffices to show $\HilS\cap
  N=\{0\}$.  Assume the contrary and take $A\in \HilS\cap N\setminus\{0\}$.
  Let $A'\in\Mat[n-1,m-1]$ be the lower right block of~$A$, i.e.\ forget the
  zero column and row of~$A$.  First consider the case $\alpha_2(A_0)<1$.  Then
  $$
	\|A_0+\lambda A\|
  =	\max\bigl\{1, \bigl\|\diag\bigl(\alpha_2(A_0),\dots,\alpha_m(A_0)\bigr)
	 +\lambda A'\bigr\| \bigr\}
  =	1
  $$
  for all scalars~$\lambda$ in a neighborhood of~$0$.  This, however, cannot
  happen in a Hilbert space.  Now assume
  $$
  1=\alpha_2(A_0)=\dots=\alpha_k(A_0)>\alpha_{k+1}(A_0)
  $$
  for some $k\ge 2$.  By construction of~$A_0$, this means that the singular
  value~$1$ occurs $k$~times for all elements of~$S(\HilS)$.  Since $A\in N$,
  the matrix $A_0+\lambda A$ still has the singular value~$1$ for all
  $\lambda\in\C$.  Furthermore, $\alpha_{k+1}(A_0+\lambda A)<1$ and
  $\|A_0+\lambda A\|\neq 1$ for sufficiently small~$|\lambda|$.  Hence, for
  small $|\lambda|>0$, the matrix $(A_0+\lambda A)/\|A_0+\lambda\|\in S(\HilS)$
  has~$k$ singular values that are at least $\|A_0+\lambda\|^{-1}$.  One of
  them is $\|A_0+\lambda\|^{-1}<1$, so that the singular value~$1$ occurs with
  multiplicity less than~$k$.  This contradicts the choice of~$A_0$.  Thus the
  assumption that $\HilS\cap N\neq\{0\}$ gives rise to a contradiction both if
  $\alpha_2(A_0)<1$ and if $\alpha_2(A_0)=1$.  This proves the claim.
\end{proof}

Indeed, the case $\alpha_2(A_0)=1$ cannot occur at all:

\begin{claim}
  $\alpha_2(A_0)<1$.
\end{claim}

\begin{proof}
  Assume that $\alpha_2(A_0)=1$.  Since $\HilS\to \Mat[n,m]/N$ is an
  isomorphism, there exists $B\in \HilS$ of the form
  $$
  B =
  \begin{pmatrix}
	0 & 1 & \ast & \ast & \dots \\
	1 & b & \ast & \ast & \dots \\
	\ast & \ast & \ast & \ast & \dots \\
	\ast & \ast & \ast & \ast & \dots \\
	\hdotsfor{5}
  \end{pmatrix}
  $$
  Since $\|A_0 + \lambda B\|\ge 1$ for all $\lambda\in\C$, the matrix~$B$ must
  be orthogonal to~$A_0$ in the Hilbert space~$\HilS$.  Therefore, $\|A_0+\lambda
  B\|^2=1+ \const\cdot|\lambda|^2$ for all $\lambda\in\C$.  On the other hand,
  $$
	\|A_0+\lambda B\|^2
  \ge	\left\| \begin{pmatrix}
	1		& \lambda	\\
	\lambda	& 1+b\lambda
	\end{pmatrix} \right\|^2
  $$
  This can further be estimated below by
  \begin{align*}
	\left\| \begin{pmatrix}
	1		& \lambda	\\
	\lambda	& 1+b\lambda
	\end{pmatrix} \right\|^2
  &\ge	\frac{1}{2}
	\left\| \begin{pmatrix}
	1		& \lambda	\\
	\lambda	& 1+b\lambda
	\end{pmatrix}
	\begin{pmatrix} \pm 1 \\ 1 \end{pmatrix} \right\|^2
   =	\frac{1}{2} \bigl( |\pm1+\lambda|^2
	 + |\pm\lambda+1+b\lambda|^2 \bigr)\\
  &=	1 + \RE(\pm2\lambda+ b\lambda) + O(\lambda^2)
  \end{align*}
  for $\lambda \to0$.  For a suitable choice of the sign, the \PrM{\lambda}term
  does not vanish, so that $\|A_0+\lambda B\|^2\ge 1 + \RE(c\lambda) +
  O(\lambda^2)$ for $\lambda\to0$, with some $c\neq0$.  However, this
  contradicts $\|A_0+\lambda B\|^2=1 + \const\cdot |\lambda|^2$.  Hence
  $\alpha_2(A_0)\neq1$.
\end{proof}

Now all the preliminary work is done, the following claim is the main part of
the proof.

\begin{claim}
  Any $A\in S(\HilS)$ with $\alpha_2(A)<1$ has rank at most~$1$.
\end{claim}

\begin{proof}
  Take any such~$A$ and bring it into diagonal form $A=\diag(1,a_2,\dots,a_m)$
  as above.  Let $B\in S(\HilS)$ with $B\bot A$ and entries $B=(b_{jk})$.  Then
  $b_{11}=0$ because otherwise $\|A+ \lambda B\| \ge |1+\lambda b_{11}|$ is not
  $1+O(|\lambda|^2)$.  Let
  $$
	M(\lambda) = (A+\lambda B)^\ast (A+\lambda B) - 1 - |\lambda|^2.
  $$
  Since $\|A+\lambda B\|^2=1+|\lambda|^2$ for all $\lambda\in\C$, the matrix
  $M(\lambda)$ is not invertible for any $\lambda\in\C$.  Hence $\lambda\mapsto
  \det M(\lambda)$ must be the zero function.

  The entries of $M(\lambda)$ are quadratic functions of $\lambda$
  and~$\conj{\lambda}$, so that $\det M(\lambda)$ can be written as a
  polynomial in $\lambda$ and~$\conj{\lambda}$.  Since $M(\lambda)$ is
  Hermitian for all $\lambda\in\C$, the function $\det M(\lambda)$ is
  real-valued.  Thus $\det M(\lambda) = c + 2\RE(d\lambda) + 2\RE (e\lambda^2)
  + f|\lambda|^2 + O(|\lambda|^3)$ with some constants $c,d,e,f\in\C$.  All of
  them have to vanish.  However, $c$, $d$, and~$f$ are irrelevant for this
  proof.

  In order to compute~$e$, it suffices to look at the terms of the form
  $\const\cdot \lambda^2$ in the expansion of $\det M(\lambda)$.  Therefore,
  the summand $(\lambda B)^\ast (A+\lambda B) -|\lambda|^2$ in $M(\lambda)$ can
  be ignored.  It remains to compute the second order term of $\det( A^\ast A
  -1 +\lambda A^\ast B)$.  Since the $1,1$th entry of $A^\ast A-1+\lambda
  A^\ast B$ is zero, all nonzero terms in the determinant expansion involve an
  off-diagonal term in the first column and one off-diagonal term in the first
  row.  Since $A^\ast A-1$ is diagonal, terms involving more than two
  off-diagonal terms are $o(\lambda^2)$ and can be ignored.  Finally, the
  result is $\det (A^\ast A - 1 +\lambda A^\ast B)= e\lambda^2 + o(\lambda^2)$
  with
  $$
	e
  =	-\sum_{k=2}^n \left( a_k b_{k1} b_{1k}
	 \prod_{2\le l\le n,l\neq k} (a_l^2-1) \right)
  =	\prod_{l=2}^n (a_l^2-1) \sum_{k=2}^n \frac{a_kb_{k1}b_{1k}}{1-a_k^2}.
  $$

  Since $a_l^2\neq1$ for $l\ge2$, the assumption $e=0$ implies
  \begin{equation}  \label{equ:abbak}
	 \sum_{k=2}^n \frac{a_kb_{k1}b_{1k}}{1-a_k^2} = 0.
  \end{equation}
  \eqref{equ:abbak} holds for any $B\in S(\HilS)$ with $B\bot A$.  But then it
  must hold for all $B\in\HilS$.  However, an application of
  Claim~\ref{claim:a} shows that the constants $b_{k1},b_{1k}$ can be
  prescribed arbitrarily.  Thus $a_k=0$ for all $k=2,\dots,n$.  But that means
  that~$A$ has rank~$1$.
\end{proof}

  Now the proof of Theorem~\ref{the:ellTwoMat} is almost finished.  Pick any
  $A_0\in S(\HilS)$ with $\alpha_2(A)<1$, then~$A_0$ has rank~$1$.  The same
  holds for all $A\in\HilS$ in a suitable neighborhood of~$A_0$
  because~$\alpha_2$ is continuous.  We can assume that
  $A_0=\diag(1,0,\dots,0)$.  By Claim~\ref{claim:a}, there is $B_1\in\HilS$
  with first row zero and first column $(0,1,0,\dots,0)^t$.  Moreover, there is
  $B_2\in\HilS$ with first column zero and first row $(0,1,\dots,0)$.

  The only chance for $A_0+\lambda B_1$ to have rank~$1$ is if $(B_1)_{jk}=0$
  for all $(j,k)\neq(2,1)$, and similarly for~$B_2$.  But then $A_0 + \lambda
  (B_1+B_2)$ has rank~$2$ for all $\lambda\in\C^\ast$.  This
  contradiction proves Theorem~\ref{the:ellTwoMat}
\end{proof}

Now we continue the proof of Theorem~\ref{the:cotangentMinimalComplexity}.
Assume first that $\CTg_\omega\OA$ has finite dimension~$r$.  Let $\rho\colon
\OA(\omega)\to\Mat[r+1]$ be a completely isometric representation.  Let~$K$ be
the intersection of the kernels of all~$\rho(A)$, $A\in\CTg_\omega\OA$, and let
$n=\dim K$.  Since the multiplication on~$\CTg_\omega\OA$ is trivial, $\Ran
\rho(A)\subset K$ for all $A\in\CTg_\omega\OA$.

By a unitary transformation, we can achieve that~$K$ is spanned by the vectors
$e_1,\dots,e_n$.  Then all elements of $\rho(\CTg_\omega\OA)$ are of the form
$$
\begin{pmatrix} 0 & \ast \\ 0 & 0 \end{pmatrix}
$$
Thus ignoring all but the upper right corner gives a (completely) isometric
representation $\CTg_\omega\OA\to\Mat[m,n]$ with $m=r+1-n$.  By
Lemma~\ref{lem:ellTwoSpecial}, the image is isometric to a Hilbert space of
dimension $r=n+m-1$.  By Theorem~\ref{the:ellTwoMat}, this is impossible unless
$m=1$ or $n=1$.  These two cases correspond to $\CTg_\omega\OA=\CTg_0\Shift_r$
and $\CTg_\omega\OA=\CTg_0\Shift_r^t$, respectively.  This proves
Theorem~\ref{the:cotangentMinimalComplexity} in the finite dimensional case.

Now assume that $\CTg_\omega\OA=\OA[B]$ is infinite dimensional.
Lemma~\ref{lem:ellTwoSpecial} implies that~$\OA[B]$ is isometric to a Hilbert
space~$\HilS$ because~$\OA$ and hence~$\OA[B]$ is closed.  We have to show
that~$\OA[B]$ is completely equivalent to $\Bound(\C,\HilS)$ or
$\Bound(\HilS,\C)=\Bound(\C,\HilS)^t$.  

We already know that every \PrMn{r}dimensional quotient of $\OA[B]$ is
completely equivalent to $\CTg_0\Shift_r^t$ or $\CTg_0\Shift_r$.  It is easy to
see that $\CTg_0\Shift_r$ is not a quotient of $\CTg_0\Shift_{r'}^t$ if $r>1$,
and similarly $\CTg_0\Shift_r^t$ cannot be a quotient of $\CTg_0\Shift_{r'}$ for
$r>1$.  Hence either all quotients are of the form $\CTg_0\Shift_r$ or of the
form $\CTg_0\Shift_r^t$, but both types cannot occur at the same time.

Replacing~$\OA[B]$ by~$\OA[B]^t$ if necessary, we can achieve that every
\PrMn{r}dimensional quotient of~$\OA[B]$ is completely isometric to
$\CTg_0\Shift_r$.  Let $\iota\colon\OA[B]\to\HilS\to\Bound(\C,\HilS)$ be the
canonical isometry.  We claim that~$\iota$ is completely isometric.  This
follows if the dual map $\iota'\colon \Bound(\C,\HilS)'\to\OA[B]'$ is a
complete isometry.  Here both dual spaces are equipped with their natural
\PrMn{L^1}matricially normed structure \cite{Effros-Ruan:88}.

If $X\in\OA[B]'_{(n)}$, its entries span a finite dimensional subspace~$L$
of~$\OA[B]'$, which determines a finite dimensional quotient $\OA[B]/L^\bot$
of~$\OA[B]$.  The map $\OA[B]/L^\bot\to \Bound(\C,\HilS)/\iota(L^\bot)$ induced
by~$\iota$ is completely isometric.  Dualizing this shows that
$\iota'|_{\iota(L)}$ is completely isometric.  This implies
$\|\iota'_{(n)}(X)\|_{(n)}=\|X\|_{(n)}$, which is what we need.
\end{proof}

Thus every finite dimensional operator algebra of minimal complexity that is a
trivial unitization is a quotient of $\Shift_d$ or of $\Shift_d^t$.  Notice
that trivial unitizations are automatically indecomposable.  In general, a
direct sum of quotients of $\Shift_d$ need not be a quotient of $\Shift_d^t$
again.  However, the following conjecture has a chance to be true:

\begin{conjecture}  \label{conj:classifyMC}
  Let~$\OA$ be an \PrMn{r}dimensional, indecomposable, (unital, commutative)
  operator algebra of minimal quotient complexity.  Then there exists a closed,
  finite codimensional ideal $\ide\subset\Shift_{r-1}$ with $\SA(\ide)=\{0\}$
  such that $\OA\cong\Shift_{r-1}/\ide$ or $\OA\cong (\Shift_{r-1}/\ide)^t$.
  Moreover, if $\OA\cong \OA^t$, then~$\OA$ is a quotient of $\Shift_1$.
\end{conjecture}

The converse of this conjecture is easy.  Hence it would provide us with a
complete classification of finite dimensional operator algebras of minimal
quotient complexity.

\begin{remark}
  Every unital, commutative subalgebra of~$\Mat[3]$ has minimal quotient
  complexity.  Indeed, it can have at most dimension~$3$.  Its non-trivial
  quotients have dimensions $1$ and~$2$, so that the assertion follows from the
  classification of two\PrMn{}dimensional, unital operator algebras.

  For subalgebras of~$\Mat[3]$, I have verified
  Conjecture~\ref{conj:classifyMC} by direct computations.  But the proof is to
  messy to be included here.  Moreover, new features arise in dimension~$4$
  because no longer all \PrMn{4}dimensional, unital, commutative subalgebras
  of~$\Mat[4]$ have minimal quotient complexity.
\end{remark}

\providecommand{\bysame}{\leavevmode\hbox to3em{\hrulefill}\thinspace}

\end{document}